\documentclass{amsart}

\numberwithin{equation}{section}

\usepackage{graphicx,amsthm,enumerate}
\usepackage{verbatim,amsmath,amscd,amssymb,color}
\usepackage[mathscr]{eucal}
\usepackage{tikz}
\usepackage[all]{xy}
\begin{document}
\renewcommand{\labelenumi}{$($\roman{enumi}$)$}
\renewcommand{\labelenumii}{$(${\rm \alph{enumii}}$)$}
\font\germ=eufm10
\newcommand{\cI}{{\mathcal I}}
\newcommand{\cA}{{\mathcal A}}
\newcommand{\cB}{{\mathcal B}}
\newcommand{\cC}{{\mathcal C}}
\newcommand{\cD}{{\mathcal D}}
\newcommand{\cE}{{\mathcal E}}
\newcommand{\cF}{{\mathcal F}}
\newcommand{\cG}{{\mathcal G}}
\newcommand{\cH}{{\mathcal H}}
\newcommand{\cK}{{\mathcal K}}
\newcommand{\cL}{{\mathcal L}}
\newcommand{\cM}{{\mathcal M}}
\newcommand{\cN}{{\mathcal N}}
\newcommand{\cO}{{\mathcal O}}
\newcommand{\cR}{{\mathcal R}}
\newcommand{\cS}{{\mathcal S}}
\newcommand{\cT}{{\mathcal T}}
\newcommand{\cV}{{\mathcal V}}
\newcommand{\cX}{{\mathcal X}}
\newcommand{\cY}{{\mathcal Y}}
\newcommand{\fra}{\mathfrak a}
\newcommand{\frb}{\mathfrak b}
\newcommand{\frc}{\mathfrak c}
\newcommand{\frd}{\mathfrak d}
\newcommand{\fre}{\mathfrak e}
\newcommand{\frf}{\mathfrak f}
\newcommand{\frg}{\mathfrak g}
\newcommand{\frh}{\mathfrak h}
\newcommand{\fri}{\mathfrak i}
\newcommand{\frj}{\mathfrak j}
\newcommand{\frk}{\mathfrak k}
\newcommand{\frI}{\mathfrak I}
\newcommand{\fm}{\mathfrak m}
\newcommand{\frn}{\mathfrak n}
\newcommand{\frp}{\mathfrak p}
\newcommand{\fq}{\mathfrak q}
\newcommand{\frr}{\mathfrak r}
\newcommand{\frs}{\mathfrak s}
\newcommand{\frt}{\mathfrak t}
\newcommand{\fru}{\mathfrak u}
\newcommand{\frA}{\mathfrak A}
\newcommand{\frB}{\mathfrak B}
\newcommand{\frF}{\mathfrak F}
\newcommand{\frG}{\mathfrak G}
\newcommand{\frH}{\mathfrak H}
\newcommand{\frJ}{\mathfrak J}
\newcommand{\frN}{\mathfrak N}
\newcommand{\frP}{\mathfrak P}
\newcommand{\frT}{\mathfrak T}
\newcommand{\frU}{\mathfrak U}
\newcommand{\frV}{\mathfrak V}
\newcommand{\frX}{\mathfrak X}
\newcommand{\frY}{\mathfrak Y}
\newcommand{\fry}{\mathfrak y}
\newcommand{\frZ}{\mathfrak Z}
\newcommand{\frx}{\mathfrak x}
\newcommand{\rA}{\mathrm{A}}
\newcommand{\rC}{\mathrm{C}}
\newcommand{\rd}{\mathrm{d}}
\newcommand{\rB}{\mathrm{B}}
\newcommand{\rD}{\mathrm{D}}
\newcommand{\rE}{\mathrm{E}}
\newcommand{\rH}{\mathrm{H}}
\newcommand{\rK}{\mathrm{K}}
\newcommand{\rL}{\mathrm{L}}
\newcommand{\rM}{\mathrm{M}}
\newcommand{\rN}{\mathrm{N}}
\newcommand{\rR}{\mathrm{R}}
\newcommand{\rT}{\mathrm{T}}
\newcommand{\rZ}{\mathrm{Z}}
\newcommand{\bbA}{\mathbb A}
\newcommand{\bbB}{\mathbb B}
\newcommand{\bbC}{\mathbb C}
\newcommand{\bbG}{\mathbb G}
\newcommand{\bbF}{\mathbb F}
\newcommand{\bbH}{\mathbb H}
\newcommand{\bbP}{\mathbb P}
\newcommand{\bbN}{\mathbb N}
\newcommand{\bbQ}{\mathbb Q}
\newcommand{\bbR}{\mathbb R}
\newcommand{\bbV}{\mathbb V}
\newcommand{\bbX}{\mathbb X}
\newcommand{\bbZ}{\mathbb Z}
\newcommand{\adj}{\operatorname{adj}}
\newcommand{\Ad}{\mathrm{Ad}}
\newcommand{\Ann}{\mathrm{Ann}}
\newcommand{\rcris}{\mathrm{cris}}
\newcommand{\ch}{\mathrm{ch}}
\newcommand{\coker}{\mathrm{coker}}
\newcommand{\diag}{\mathrm{diag}}
\newcommand{\Diff}{\mathrm{Diff}}
\newcommand{\Dist}{\mathrm{Dist}}
\newcommand{\rDR}{\mathrm{DR}}
\newcommand{\ev}{\mathrm{ev}}
\newcommand{\Ext}{\mathrm{Ext}}
\newcommand{\cExt}{\mathcal{E}xt}
\newcommand{\fin}{\mathrm{fin}}
\newcommand{\Frac}{\mathrm{Frac}}
\newcommand{\GL}{\mathrm{GL}}
\newcommand{\Hom}{\mathrm{Hom}}
\newcommand{\hd}{\mathrm{hd}}
\newcommand{\rht}{\mathrm{ht}}
\newcommand{\id}{\mathrm{id}}
\newcommand{\im}{\mathrm{im}}
\newcommand{\inc}{\mathrm{inc}}
\newcommand{\ind}{\mathrm{ind}}
\newcommand{\coind}{\mathrm{coind}}
\newcommand{\Lie}{\mathrm{Lie}}
\newcommand{\Max}{\mathrm{Max}}
\newcommand{\mult}{\mathrm{mult}}
\newcommand{\op}{\mathrm{op}}
\newcommand{\ord}{\mathrm{ord}}
\newcommand{\pt}{\mathrm{pt}}
\newcommand{\qt}{\mathrm{qt}}
\newcommand{\rad}{\mathrm{rad}}
\newcommand{\res}{\mathrm{res}}
\newcommand{\rgt}{\mathrm{rgt}}
\newcommand{\rk}{\mathrm{rk}}
\newcommand{\SL}{\mathrm{SL}}
\newcommand{\soc}{\mathrm{soc}}
\newcommand{\Spec}{\mathrm{Spec}}
\newcommand{\St}{\mathrm{St}}
\newcommand{\supp}{\mathrm{supp}}
\newcommand{\Tor}{\mathrm{Tor}}
\newcommand{\Tr}{\mathrm{Tr}}
\newcommand{\wt}{\mathrm{wt}}
\newcommand{\Ab}{\mathbf{Ab}}
\newcommand{\Alg}{\mathbf{Alg}}
\newcommand{\Grp}{\mathbf{Grp}}
\newcommand{\Mod}{\mathbf{Mod}}
\newcommand{\Sch}{\mathbf{Sch}}\newcommand{\bfmod}{{\bf mod}}
\newcommand{\Qc}{\mathbf{Qc}}
\newcommand{\Rng}{\mathbf{Rng}}
\newcommand{\Top}{\mathbf{Top}}
\newcommand{\Var}{\mathbf{Var}}
\newcommand{\pmbx}{\pmb x}
\newcommand{\bfp}{\mathbf p}
\newcommand{\pmby}{\pmb y}
\newcommand{\gromega}{\langle\omega\rangle}
\newcommand{\lbr}{\begin{bmatrix}}
\newcommand{\rbr}{\end{bmatrix}}
\newcommand{\cd}{commutative diagram }
\newcommand{\SpS}{spectral sequence}
\newcommand\C{\mathbb C}
\newcommand\hh{{\hat{H}}}
\newcommand\eh{{\hat{E}}}
\newcommand\F{\mathbb F}
\newcommand\fh{{\hat{F}}}
\newcommand\Z{{\mathbb Z}}
\newcommand\Zn{\Z_{\geq0}}
\newcommand\et[1]{\tilde{e}_{#1}}
\newcommand\ft[1]{\tilde{f}_{#1}}

\def\ge{\frg}
\def\AA{{\mathcal A}}
\def\al{\alpha}
\def\bq{B_q(\ge)}
\def\bqm{B_q^-(\ge)}
\def\bqz{B_q^0(\ge)}
\def\bqp{B_q^+(\ge)}
\def\beneme{\begin{enumerate}}
\def\beq{\begin{equation}}
\def\beqn{\begin{eqnarray}}
\def\beqnn{\begin{eqnarray*}}
\def\bfi{{\mathbf i}}
\def\bfii0{{\bf i_0}}
\def\bigsl{{\hbox{\fontD \char'54}}}
\def\bbra#1,#2,#3{\left\{\begin{array}{c}\hspace{-5pt}
#1;#2\\ \hspace{-5pt}#3\end{array}\hspace{-5pt}\right\}}
\def\cd{\cdots}
\def\ci(#1,#2){c_{#1}^{(#2)}}
\def\Ci(#1,#2){C_{#1}^{(#2)}}
\def\mpp(#1,#2,#3){#1^{(#2)}_{#3}}
\def\bCi(#1,#2){\ovl C_{#1}^{(#2)}}
\def\ch(#1,#2){c_{#2,#1}^{-h_{#1}}}
\def\cHZ{{\mathcal H}_{{\bf i}_0}}
\def\cc(#1,#2){c_{#2,#1}}
\def\CC{\mathbb{C}}
\def\CBL{\cB_L(\TY(B,1,n+1))}
\def\CBM{\cB_M(\TY(B,1,n+1))}
\def\CVL{\cV_L(\TY(D,1,n+1))}
\def\CVM{\cV_M(\TY(D,1,n+1))}
\def\ddd{\hbox{\germ D}}
\def\del{\delta}
\def\Del{\Delta}
\def\Delr{\Delta^{(r)}}
\def\Dell{\Delta^{(l)}}
\def\Delb{\Delta^{(b)}}
\def\Deli{\Delta^{(i)}}
\def\Delre{\Delta^{\rm re}}
\def\di(#1,#2){D_{#1}^{(#2)}}
\def\dbi(#1,#2){\ovl D_{#1}^{(#2)}}
\def\ei{e_i}
\def\eit{\tilde{e}_i}
\def\eneme{\end{enumerate}}
\def\EEE#1{{\rm E}_{#1}}
\def\F44{{\rm F}_{4}}
\def\G22{{\rm G}_{2}}
\def\AAN{{\rm A}_n}
\def\BBN{{\rm B}_n}
\def\CCN{{\rm C}_n}
\def\DDN{{\rm D}_n}
\def\ep{\epsilon}
\def\eeq{\end{equation}}
\def\eeqn{\end{eqnarray}}
\def\eeqnn{\end{eqnarray*}}
\def\fit{\tilde{f}_i}
\def\fulp{\Phi_{\rm BK}}
\def\FF{{\rm F}}
\def\ft{\tilde{f}_}
\def\gau#1,#2{\left[\begin{array}{c}\hspace{-5pt}#1\\
\hspace{-5pt}#2\end{array}\hspace{-5pt}\right]}
\def\gl{\hbox{\germ gl}}
\def\hom{{\hbox{Hom}}}
\def\HZ{{\mathcal H}_{{\bf i}_0}}
\def\HI{{\mathcal H}_{\bf i}}
\def\CHI{{\bf H}_{\bf i}}
\def\CHZ{{\bf H}_{{\bf i}_0}}
\def\ify{\infty}
\def\io{\iota}
\def\ist{{i^*}}
\def\ji(#1,#2){j_{#1}^{(#2)}}
\def\kp{k^{(+)}}
\def\km{k^{(-)}}
\def\llra{\relbar\joinrel\relbar\joinrel\relbar\joinrel\rightarrow}
\def\lan{\langle}
\def\lar{\longrightarrow}
\def\LL(#1){{#1}^\vee}
\def\lm{\lambda}
\def\Lm{\Lambda}
\def\mapright#1{\smash{\mathop{\longrightarrow}\limits^{#1}}}
\def\Mapright#1{\smash{\mathop{\Longrightarrow}\limits^{#1}}}
\def\mm{{\bf{\rm m}}}
\def\nd{\noindent}
\def\nn{\nonumber}
\def\nnn{\hbox{\germ n}}
\def\catob{{\mathcal O}(B)}
\def\oint{{\mathcal O}_{\rm int}(\ge)}
\def\ot{\otimes}
\def\op{\oplus}
\def\opi{\ovl\pi_{\lm}}
\def\osigma{\ovl\sigma}
\def\ovl{\overline}
\def\plm{\Psi^{(\lm)}_{\io}}
\def\qq{\qquad}
\def\q{\quad}
\def\qed{\hfill\framebox[2mm]{}}
\def\QQ{\mathbb Q}
\def\qi{q_i}
\def\qii{q_i^{-1}}
\def\ra{\rightarrow}
\def\ran{\rangle}
\def\rlm{r_{\lm}}
\def\ssl{\hbox{\germ sl}}
\def\slh{\widehat{\ssl_2}}
\def\syl{\scriptstyle}
\def\ti{t_i}
\def\tii{t_i^{-1}}
\def\til{\tilde}
\def\tm{\times}
\def\tri{\bigtriangleup}
\def\tt{\frt}
\def\TY(#1,#2,#3){#1^{(#2)}_{#3}}
\def\ua{U_{\AA}}
\def\ue{U_{\vep}}
\def\uq{U_q(\ge)}
\def\uqp{U'_q(\ge)}
\def\ufin{U^{\rm fin}_{\vep}}
\def\ufinp{(U^{\rm fin}_{\vep})^+}
\def\ufinm{(U^{\rm fin}_{\vep})^-}
\def\ufinz{(U^{\rm fin}_{\vep})^0}
\def\uqm{U^-_q(\ge)}
\def\uqmq{{U^-_q(\ge)}_{\bf Q}}
\def\uqpm{U^{\pm}_q(\ge)}
\def\uqq{U_{\bf Q}^-(\ge)}
\def\uqz{U^-_{\bf Z}(\ge)}
\def\ures{U^{\rm res}_{\AA}}
\def\urese{U^{\rm res}_{\vep}}
\def\uresez{U^{\rm res}_{\vep,\ZZ}}
\def\util{\widetilde\uq}
\def\uup{U^{\geq}}
\def\ulow{U^{\leq}}
\def\bup{B^{\geq}}
\def\blow{\ovl B^{\leq}}
\def\vep{\varepsilon}
\def\vp{\varphi}
\def\vpi{\varphi^{-1}}
\def\VV{{\mathcal V}}
\def\xii{\xi^{(i)}}
\def\Xiioi{\Xi_{\io}^{(i)}}
\def\xxi(#1,#2,#3){\displaystyle {}^{#1}\Xi^{(#2)}_{#3}}
\def\xsi(#1,#2,#3){\displaystyle {}^{#1}\Sigma^{(#2)}_{#3}}
\def\xE(#1,#2,#3){\displaystyle {}^{#1}E_{#2}[#3]}
\def\xF(#1,#2){\displaystyle {}^{#1}F_{#2}}
\def\xx(#1,#2){\displaystyle {}^{#1}\Xi_{#2}}
\def\W1{W(\varpi_1)}
\def\WW{{\mathcal W}}
\def\wt{{\rm wt}}
\def\wtil{\widetilde}
\def\what{\widehat}
\def\wpi{\widehat\pi_{\lm}}
\def\ZZ{\mathbb Z}
\def\RR{\mathbb R}

\def\m@th{\mathsurround=0pt}
\def\fsquare(#1,#2){
\hbox{\vrule$\hskip-0.4pt\vcenter to #1{\normalbaselines\m@th
\hrule\vfil\hbox to #1{\hfill$\scriptstyle #2$\hfill}\vfil\hrule}$\hskip-0.4pt
\vrule}}

\newtheorem{thm}{Theorem}[section]
\newtheorem{pro}[thm]{Proposition}
\newtheorem{lem}[thm]{Lemma}
\newtheorem{ex}[thm]{Example}
\newtheorem{cor}[thm]{Corollary}
\newtheorem{conj}[thm]{Conjecture}
\newtheorem{rem}[thm]{Remark}
\theoremstyle{definition}
\newtheorem{df}[thm]{Definition}

\newcommand{\cmt}{\marginpar}
\newcommand{\seteq}{\mathbin{:=}}
\newcommand{\cl}{\colon}
\newcommand{\be}{\begin{enumerate}}
\newcommand{\ee}{\end{enumerate}}
\newcommand{\bnum}{\be[{\rm (i)}]}
\newcommand{\enum}{\ee}
\newcommand{\ro}{{\rm(}}
\newcommand{\rf}{{\rm)}}
\newcommand{\set}[2]{\left\{#1\,\vert\,#2\right\}}
\newcommand{\sbigoplus}{{\mbox{\small{$\bigoplus$}}}}
\newcommand{\ba}{\begin{array}}
\newcommand{\ea}{\end{array}}
\newcommand{\on}{\operatorname}
\newcommand{\eq}{\begin{eqnarray}}
\newcommand{\eneq}{\end{eqnarray}}
\newcommand{\hs}{\hspace*}

\title[Half Potential and Cellular Crystals]
{Half Potential on Geometric Crystals and 
Connectedness of Cellular Crystals}

\author{Yuki K\textsc{anakubo}\quad Toshiki N\textsc{akashima}}
\address{Division of Mathematics, 
Sophia University, Kioicho 7-1, Chiyoda-ku, Tokyo 102-8554,
Japan}
\email{Y.K.:j\_chi\_sen\_you\_ky@eagle.sophia.ac.jp,\,\,T.N.:toshiki@sophia.ac.jp}
\thanks{T.N. is supported in part by JSPS Grants 
in Aid for Scientific Research $\sharp 15K04794$.}

\date{}



\begin{abstract}
For any simply connected, simple complex algebraic group,
we define upper/lower half-decorated geometric crystals and show that
their tropicalization will be upper/lower normal Kashiwara's crystals.
In particular, we show that 
the tropicalization of the half-decorated geometric crystal
on the big Bruhat cell(=$B^-_{w_0}:=B^-\cap U\ovl w_0 U$)
is isomorphic to the Langlands dual crystal $B(\infty)$ of
the nilpotent-half subalgebra of quantum group. As an application, 
we shall show that any cellular crystal associated with a reduced word 
is connected in the sense of a crystal graph.
\end{abstract}

\maketitle

\renewcommand{\thesection}{\arabic{section}}
\section{Introduction}
\setcounter{equation}{0}
\renewcommand{\theequation}{\thesection.\arabic{equation}}

The notion of  geometric crystals for any reductive algebraic group $G$
(or $\ge={\rm Lie}(G)$)
was introduced by 
Berenstein and Kazhdan (\cite{BK}) 
as a sort of geometric lifting of 
Kashiwara's crystal base theory(\cite{K0,K1}), which has been generalized to 
the affine/Kac-Moody settings (\cite{KNO,KNO2,N}). 
Applying the tropicalization functor, certain positive 
geometric crystals are transferred 
to some Langlands dual free-Kashiwara's crystals (\cite{BK}).
Furthermore, in \cite{BK2} Berenstein and Kazhdan initiated 
the theory of decorated geometric crystals, which are  geometric crystals
equipped with the rational function $f$ 
called ``decoration'' or ``potential'', 
satisfying the condition
\begin{equation}
 f(e_i^c(x))=f(x)+(c-1)\vp_i(x)+(c^{-1}-1)\vep_i(x)
\q
(i\in I:=\{1,2,\cd,n\},\,c\in \bbC^\times),
\label{intro-f}
\end{equation}
where $e_i$ is the geometric crystal operator and $\vep_i$ and $\vp_i$
are certain structure functions of geometric crystals.
Here, by tropicalizing a positive decorated geometric crystal 
$(X,f)$ to $(\wtil X,\wtil f)$, 
where the crystal $\wtil X=\bbZ^{{\rm dim} X}$ 
is the tropicalization of $X$ and 
the piece-wise linear function $\wtil f$ is the tropicalization of 
the  decoration/potential $f$ on $\wtil X$, they 
define the subset
\[
\wtil B:=\{x\in\wtil X\,|\,\wtil f(x)\geq0\},
\]
and show that $\wtil B$ is a Langlands dual normal Kashiwara's crystal. 
Here a {\it normal Kashiwara's crystal} 
is defined to be a crystal satisfying 
\[
\vep_i(x):=\max\{n\,|\,\eit^n(x)\in\wtil B\},\q
\vp_i(x):=\max\{n\,|\,\fit^n(x)\in\wtil B\},
\]
for any $i\in I$, where $\vep_i$ and $\vp_i$ are the 
structure functions of crystals. 
And then, they realized the decorated geometric crystal on the big
cell  $TB^-_{w_0}:=T(B^-\cap U\ovl w_0 U)$, where 
$w_0$ is the longest Weyl group element, 
$T$ is a maximal torus of $G$, $U$ is a unipotent radical of $B$
and $B^-$ is an opposite Borel subgroup in $G$, 
and they described the potential 
$\Phi_{\rm BK}$ explicitly as 
\[
\Phi_{\rm BK}(g)= \sum_i\frac{\Del_{w_0\Lm_i,s_i\Lm_i}(g)}
{\Del_{w_0\Lm_i,\Lm_i}(g)}+\sum_i
\frac{\Del_{w_0s_i\Lm_i,\Lm_i}(g)}{\Del_{w_0\Lm_i,\Lm_i}(g)},
\]
where  $g\in G$, $\Del_{\gamma,\del}$ is a  generalized
minor associated with weights $\gamma$ and $\del$, 
$s_i$ is the simple reflection of the Weyl group $W$ 
and $\Lm_i$ is the $i$-th fundamental weight. 
Finally, they have shown the following remarkable results:
For the positive decorated geometric crystal
 $(TB^-_{w_0}, \fulp,T\Theta^-_{\bfi})$ 
its tropicalization is defined by:
\[
B_{\fulp,T\Theta^-_\bfi}:=(\wtil B_{\wtil\fulp,\bfi},
\eit, {\rm wt}_i,\wtil\vep_i)\q
(\wtil B_{\wtil\fulp,\bfi}:=\{\til x\in \bbZ^N\,|\,
\wtil\fulp(\til x)\geq0\}).
\]
Then 
$B_{\fulp,T\Theta^-_{\bfi}}$ is a 
{normal Kashiwara's $\ge^\vee$-crystal} and its proper sub-crystal
$B^\lm_{\fulp,T\Theta^-_{\bfi}}$ ($\lm\in {}^LP_+$) is isomorphic to 
the  $\ge^\vee$-crystal base $B(\lm)$ of the 
irreducible highest weight module $V(\lm)$ with the highest weight $\lm$, where
$\ge^\vee$ is the Langlands dual Lie algebra of $\ge$.
To show this result they used the Joseph's characterization theorem
of the crystals $\{B(\lm)\,|\,\lm\in P_+\}$ (\cite{Jo}).

In this paper, first we shall define ``upper/lower 
half-decorated geometric crystals'', which are geometric crystals 
equipped with the function
$f$ called an upper/lower half potential satisfying
\[
 f(e_i^c(x))=f(x)+(c^{-1}-1)\vep_i(x),
\q / \q f(e_i^c(x))=f(x)+(c-1)\vp_i(x),
\]
for any $i\in I$. Indeed, 
these formulas are regarded as  ``half'' of the formula \eqref{intro-f}
and then we call them ``half decorations/potentials''.
Using these properties, the upper/lower normality 
of the tropicalized decorated geometric crystal $\wtil B$ will be 
shown.

Next, for the big Bruhat cell $B^-_{w_0}$, we shall define 
a certain positive geometric crystal structure
explicitly and furthermore, the upper/lower half potential 
$\Phi^{(\pm)}$ on $B^-_{w_0}$
will be given as a half of $\Phi_{\rm BK}$:
\[
\Phi^{(+)}:=\sum_i\frac{\Del_{w_0\Lm_i,s_i\Lm_i}(g)}
{\Del_{w_0\Lm_i,\Lm_i}(g)},\qq
\Phi^{( - )}:=\sum_i
\frac{\Del_{w_0s_i\Lm_i,\Lm_i}(g)}{\Del_{w_0\Lm_i,\Lm_i}(g)}.
\]
Tropicalizing 
the positive upper half decorated geometric crystal $\bbB$
on $B^-_{w_0}$, we have
the free Kashiwara crystal $\wtil\bbB=\bbZ^{l(w_0)}$ 
and a piece-wise linear function
$\wtil \Phi^{(+)}$. Then, we obtain the upper normal crystal 
\[
\cB:=\{x\in \wtil\bbB\,|\,\wtil \Phi^{(+)}(x)\geq0\}.
\]
We shall show that this is isomorphic to the Langlands dual crystal 
$B(\ify)$ of the nilpotent quantum algebra $U_q(\ge^\vee)$ (\cite{K1,L}). 
To show that the crystal $\cB$ is isomorphic to
the $\ge^\vee$-crystal  $B(\ify)$ we adapted the 
Kashiwara-Saito characterization theorem \cite{KS}, which will be 
presented as Theorem \ref{KS-char} and Corollary \ref{ks-cor}. 


Finally, we will show the connectedness of the cellular crystal 
$B_{i_1i_2\cd i_k}$ associated with 
a reduced word $i_1i_2\cd i_k$ for a Weyl group element.
Here the cellular crystal 
$B_{i_1i_2,\cd,i_k}$ is defined by
\[
B_{i_1i_2,\cd,i_k}:=B_{i_1}\ot\cd\ot B_{i_k},
\]
where $B_i:=\{(x)_i\,|\,x\in\bbZ\}$ 
is the crystal associated with the simple root $\al_i$
(see Definition \ref{Bi}).
Indeed, to show the connectedness of $B_{i_1\cd i_k}$ it suffices to show 
the connectedness of $B_{\bfi_0}$ for a specific reduced longest word 
$\bfi_0$ due to the braid-type isomorphisms (Proposition \ref{braid}).
First, we define the subset $\HZ$ to be the solution space of some 
system of linear equations in $B_{\bfi_0}=\bbZ^{l(\bfi_0)}$. Then 
we obtain the crystal structure on the subset 
\[
B^H(\ify):=\{x+H\in \bbZ^N\,|\,x\in B(\ify)\}\q(H\in\HZ),
\]
by shifting the crystal $B(\ify)\subset \bbZ^N=B_{\bfi_0}$ along 
an element in $\HZ$.
Next, we regard the condition 
$\CHZ$ which requires that the crystal $B(\ify)$ is defined by a
system of linear inequalities in some specific form
(see Definition \ref{cond-h}).
Then, under the condition $\CHZ$, we shall show that 
$\cup_{H\in\HI}B^H(\ify)=B_{\bfi_0}$ and $B(\ify)\cap B^H(\ify)\ne\emptyset$
for any $H\in \HZ$, 
by which we find that the cellular crystal $B_{\bfi_0}$ is 
connected as a crystal graph and then so is $B_{i_1i_2\cd i_k}$
for an arbitrary reduced word $i_1i_2\cd i_k$.
For all finite-dimensional simple Lie algebras, to show the condition $\CHZ$ 
holds for some specific reduced longest word 
$\bfi_0$, we review the results in \cite{N6} and the
Berenstein-Zelevinsky's 
{\it trail} in \cite{BZ2} to obtain the explicit forms of 
all monomials in $\Phi^{(+)}$. As an application of the connectedness of the 
cellular crystals, we will obtain that 
the geometric crystal $B^-_{\bfi}$ for a reduced word $\bfi=i_1i_2\cd i_k$ 
is prehomogeneous (Corollary \ref{prehom}).

The organization of the article is as follows:
after the introduction, in Sect.2 we will give the basic setting 
of this article, defining objects such as, 
Cartan matrix, root data, weights, Weyl groups 
etc. 
In Sect.3, we review briefly the notion of crystals and 
their fundamental properties.
Furthermore, the Kashiwara embedding of the crystal $B(\ify)$
of the nilpotent subalgebra $\uqm$ of the quantum algebra $\uq$ will be 
introduced. We also review the braid-type isomorphisms and 
the Kashiwara-Saito characterization of $\uqm$. 
In Sect.4, after reviewing the theory of geometric crystals, we shall introduce
the half-decorated geometric crystals and show that 
their tropicalization will be upper/lower normal Kashiwara's crystals.
In Sect.5, 
we shall construct the half-decorated geometric crystal on the 
big Bruhat cell $B^-_{w_0}$ and present the first of the 
main results of the article, which asserts that the tropicalization of this 
geometric crystal will be isomorphic to the Langlands dual crystal $B(\ify)$. 
In Sect.6, we shall review the results in \cite{N6} and some 
explicit forms of the generalized minors $\Del_{w_0\Lm_i,s_i\Lm_i}$
for the types ${ \AAN, \BBN,\CCN,\DDN,\G22}$.
In the last subsection, we compare Theorem \ref{b-ify} in Sect.5, and the
image of Kashiwara embedding and the polyhedral realizations in \cite{NZ}. 
In Sect.7, we review the notion {\it trail} \cite{BZ2} to obtain some 
explicit forms of $\Del_{w_0\Lm_i,s_i\Lm_i}$ for the types 
${\EEE6,\EEE7,\EEE8,\F44}$. 
The last section is devoted to proving
 the connectedness of the cellular crystal $B_{i_1}\ot\cd\ot
B_{i_k}$.

The authors deeply grateful to the referee for his/her 
thorough review for this 
article and well-considered suggestions, in particular, 
he/she pointed out some gap in the proof of Theorem \ref{b-ify} and 
gave an advice to resolve it. Then, the article would 
become much nicer than before.

\renewcommand{\thesection}{\arabic{section}}
\section{Preliminaries and Settings}
\setcounter{equation}{0}
\renewcommand{\theequation}{\thesection.\arabic{equation}}

We see the settings used in this paper.
Let  
 $A=({\bf a}_{ij})_{i,j\in I}$ be 
an incomparable Cartan matrix
with a finite index set $I:=\{1,2,\cd,n\}$
and $(\frt,\{\al_i\}_{i\in I},\{h_i\}_{i\in I})$ 
be the associated
root data 
satisfying $\al_j(h_i)={\bf a}_{ij}$ where 
$\al_i\in \tt^*$ is a simple root and 
$h_i\in \tt$ is a simple coroot.
Let $\ge=\ge(A)=\lan \frt,e_i,f_i(i\in I)\ran$ 
be the finite-dimensional simple Lie algebra associated with $A$
over $\bbC$ and $\Del=\Del_+\sqcup\Del_-$
be the root system associated with 
$\ge=\frn_+\oplus \frt\oplus \frn_-$, where $\Del_{\pm}$ is 
the set of positive/negative roots.
Let $P \subset \tt^*$ be the weight lattice, 
$\lan h,\lm\ran=\lm(h)$ the pairing between $\tt$ and $\tt^*$,
and $(\al, \beta)$ be an inner product on
$\tt^*$ such that $(\al_i,\al_i)\in 2\bbZ_{\geq 0}$ and
$\lan h_i,\lm\ran={{2(\al_i,\lm)}\over{(\al_i,\al_i)}}$
for $\lm\in\tt^*$.
Let us define the 
co-weight lattice $P^*=\{h\in \tt: \lan h,P\ran\subset\ZZ\}$ and
the set of dominant weights 
$P_+:=\{\lm\in P:\lan h_i,\lm\ran\in\ZZ_{\geq 0}\}$.
Let $\{\Lm_i|i\in I\}$ be the set of the fundamental weights satisfying
$\lan h_i,\Lm_j\ran=\del_{i,j}$ which is a $\bbZ$-basis of $P$.

Define the simple reflections $s_i\in{\rm Auto}(\frt)$ $(i\in I)$ by
$s_i(h)\seteq h-\al_i(h)h_i$, which generate the Weyl group $W$.
Let $G$ be the simply connected simple algebraic group 
over $\bbC$ whose Lie algebra ${\rm Lie}(G)=\ge$.
Let $U_{\al}\seteq\exp\ge_{\al}$ $(\al\in \Del)$
be the one-parameter subgroup of $G$ associated to $\al$ and 
$T$ a maximal torus of $G$ 
which has $P$ as its weight lattice such that 
Lie$(T)=\frt$.
The groups $U^\pm$ are generated by 
$\{U_\al|\al\in\Del_{\pm}\}$.
Here $U^\pm$ is a unipotent radical of $G$ and 
${\rm Lie}(U^\pm)=\frn_{\pm}$.
For any $i\in I$, there exists 
a unique group homomorphism
$\phi_i\cl SL_2(\bbC)\rightarrow G$ such that
$\phi_i
\left(\begin{smallmatrix}
1&t\\
0&1
\end{smallmatrix}\right)
=\exp(t e_i)=:x_i(t)$, 
$\phi_i
\left(\begin{smallmatrix}
1&0\\
t&1
\end{smallmatrix}\right)
=\exp(t f_i)=:y_i(t)$, 
$\phi_i
\left(\begin{smallmatrix}
c&0\\
0&c^{-1}
\end{smallmatrix}\right)=c^{h_i}=\al^\vee_i(c)$
$(c\in \bbC^\times,t\in\bbC).$
Set 
$G_i\seteq\phi_i(SL_2(\bbC))$,
$T_i\seteq \alpha_i^\vee(\bbC^\times)$ 
and 
$N_i\seteq N_{G_i}(T_i)$. 
Let 
$B^{\pm}(\supset T,U^{\pm})$ be the Borel subgroup of $G$.
We have the isomorphism
$\phi:W\mapright{\sim}N_G(T)/T$ defined by $\phi(s_i)=N_iT/T$.
For each $s_i\in W\cong N_G(T)/T$, we choose the representative
 $\ovl s_i:=x_i(-1)y_i(1)x_i(-1)\in N_G(T)$, which is known to 
satisfy the braid relations.

For a split algebraic torus $T$ over $\bbC$, let us denote 
its lattice of (multiplicative) characters (resp. co-characters) by $X^*(T)$ 
(resp. $X_*(T)$). By the usual way, we identify $X^*(T)$
(resp. $X_*(T)$) with the weight lattice $P$ (resp. the dual weight
lattice $P^*$). 

For a finite-dimensional simple Lie algebra $\ge$ 
(resp. simple algebraic group $G$) let us denote 
its Langlands dual by $\LL(\ge)$ (resp. $\LL(G)$) whose
corresponding Cartan matrix is the transposed one ${}^t A$.
The Langlands dual of the maximal torus is denoted by $\LL(T)$.
We also denote the Langlands dual of the 
weight lattice (resp. co-weight lattice)
by $\LL(P)$ (resp. $\LL(P^*)$), which is identified with 
$X^*(\LL(T))$ (resp. $X_*(\LL(T))$). 
Now, we write the $\bbZ$-basis of each lattice as follows:
\begin{equation}
\begin{split}
&P=X^*(T)=\bigoplus_i\bbZ \Lm_i,\qq P^*=X_*(T)=\bigoplus_i\bbZ h_i,\\
&\LL(P)=X^*(\LL(T))=\bigoplus_i\bbZ h^\vee_i,\qq
\LL(P^*)=X_*(\LL(T))=\bigoplus_i\bbZ \Lm^\vee_i,
\end{split}\label{PPPP}
\end{equation}
where we have 
\[
\lan \Lm_i,h_j\ran=\lan h^\vee_i,\Lm^\vee_j\ran=\del_{ij}.
\]
Here, we identify the $i$-th simple root of $\LL(\ge)$
with the basis element $h^\vee_i$  of $\LL(P)$.

\renewcommand{\thesection}{\arabic{section}}
\section{Crystals}
\label{sect-crystal}
\setcounter{equation}{0}
\renewcommand{\theequation}{\thesection.\arabic{equation}}

\subsection{Crystal bases and Crystals}

The quantum algebra $\uq$
is an associative
$\QQ(q)$-algebra generated by the $e_i$, $f_i \,\, (i\in I)$,
and $q^h \,\, (h\in P^*)$
satisfying the usual relations (see e.g.,\cite{K1}), 
where we use the same notations for the 
generators $e_i$ and $f_i$ of $\uq$ as the ones for $\ge$.
The algebra $\uqm$ is the subalgebra of $\uq$ generated 
by the $f_i$ $(i\in I)$.
For the irreducible highest weight module of $\uq$
with the highest weight $\lm\in P_+$, we denote it by $V(\lm)$
and we denote its {\it crystal base}  by $(L(\lm),B(\lm))$.
Similarly, for the crystal base of the algebra $\uqm$ we denote 
$(L(\ify),B(\ify))$ (see \cite{K1,L}).
Let $\pi_{\lm}:\uqm\longrightarrow V(\lm)\cong \uqm/{\sum_i\uqm
f_i^{1+\lan h_i,\lm\ran}}$ be the canonical projection and 
$\widehat \pi_{\lm}:L(\ify)/qL(\ify)\longrightarrow L(\lm)/qL(\lm)$
be the induced map from $\pi_{\lm}$. Here note that 
$\widehat \pi_{\lm}(B(\ify))=B(\lm)\sqcup\{0\}$.
Let $*:\uq\to\uq$ be the anti-involution defined by 
$e^*_i=e_i$, $f^*_i=f_i$ and $(q^h)^*=q^{-h}$. Then, 
in \cite{K1,K3} it has been shown 
\begin{equation}
L(\ify)^*=L(\ify),\qq
B(\ify)^*=B(\ify).
\label{*=}
\end{equation}

By the terminology {\it crystal } we mean a combinatorial object 
obtained by abstracting the properties of crystal bases as follows
(see \cite{K3,NZ,N2}):
\begin{df}\label{cryst}
A 6-tuple $(B,\wt, \{\vep_i\},\{\vp_i\}, \{\eit\},\{\fit\})_{i\in I}$ is 
a {\it crystal} if $B$ is a set and $\exists0\not\in B$ and maps :
\begin{eqnarray}
&&{\rm wt}:B\to P,\label{wtp-c}\q
\vep_i:B\to \bbZ\sqcup\{-\ify\},\q\vp_i:B\to\bbZ\sqcup\{-\ify\}\q\,(i\in I)
\\
&&\eit:B\sqcup\{0\}\to B\sqcup\{0\},\q
\fit:B\sqcup\{0\}\to B\sqcup\{0\}\,\,(i\in I),\label{eitfit-c}
\end{eqnarray}
satisfying :
\begin{enumerate}
\item
$\vp_i(b)=\vep_i(b)+\lan h_i,\wt(b)\ran$.
\item
If  $b,\eit b\in B$, then $\wt(\eit b)=\wt(b)+\al_i$, 
$\vep_i(\eit b)=\vep_i(b)-1$, $\vp_i(\eit b)=\vp_i(b)+1$.
\item
If $b,\fit b\in B$, then $\wt(\fit b)=\wt(b)-\al_i$, 
$\vep_i(\fit b)=\vep_i(b)+1$, $\vp_i(\fit b)=\vp_i(b)-1$.
\item
For $b,b'\in B$ and $i\in I$, one has 
$\fit b=b'$ iff $b=\eit b'.$
\item
If $\vp_i(b)=-\ify$ for $b\in B$, then $\eit b=\fit b=0$
and $\eit(0)=\fit(0)=0$.
\end{enumerate}

\noindent
By denoting $\fit b_1=b_2$ ($b_1,b_2\in B$) with $b_1\mapright{i} b_2$, 
the crystal $B$ holds an oriented colored graph structure, which is 
called a {\it crystal graph}.

\end{df}

A crystal $B$ is said to be {\it upper (resp. lower) normal} if 
the function $\vep_i$ (resp. $\vp_i$) satisfies:
\begin{equation}
\vep_i(b)=\max\{n\,|\,\eit^n b\ne 0\}\q({\rm resp.} \,
\vp_i(b)=\max\{n\,|\,\fit^n b\ne0\}).\label{semi-n}
\end{equation}
If $B$ is upper and lower normal, it is said to be 
{\it normal}.
It is well-known that  $B(\ify)$ is an upper normal crystal 
and each $B(\lm)$ is a normal crystal.

\begin{df}\label{morph}
For crystals $B_1,\, B_2$, 
$\Psi$ is a {\it strict embedding} (resp. an {\it isomorphism})
from $B_1$ to $B_2$ if 
$\Psi:B_1\sqcup\{0\}\to B_2\sqcup\{0\}$ is an injective
(resp. a bijective) map satisfying : 
\begin{enumerate}
\item
$\Psi(0)=0$.
\item
For $b\in B_1$, $\wt(\Psi(b))=\wt(b)$, 
$\vep_i(\Psi(b))=\vep_i(b)$ and 
$\vp_i(\Psi(b))=\vp_i(b)$.
\item
$\Psi:B_1\to B_2$ commutes with 
all $\eit$'s and $\fit$'s,.
\end{enumerate}
\end{df}

\subsection{Tensor Product of Crystals}
Let $B_i$ ($i=1,2$) be crystals. 
By defining as follows, we obtain the tensor structure of crystals: 
Set 
\[
B_1\ot B_2=\{b_1\ot b_2: b_1\in B_1 ,\, b_2\in B_2\}.
\]
Then, $B_1\ot B_2$ is a crystal by setting:
\begin{eqnarray}
&&\wt(b_1\ot b_2)=\wt(b_1)+\wt(b_2),\\
&&\vep_i(b_1\ot b_2)={\hbox{max}}(\vep_i(b_1),
  \vep_i(b_2)-\lan h_i,\wt(b_1)\ran),
\label{tensor-vep}\\
&&\vp_i(b_1\ot b_2)={\hbox{max}}(\vp_i(b_2),
  \vp_i(b_1)+\lan h_i,\wt(b_2)\ran),
\label{tensor-vp}\\
&&\eit(b_1\ot b_2)=
\left\{
\begin{array}{ll}
\eit b_1\ot b_2 & {\mbox{ if }}\vp_i(b_1)\geq \vep_i(b_2)\\
b_1\ot\eit b_2  & {\mbox{ if }}\vp_i(b_1)< \vep_i(b_2),
\end{array}
\right.
\label{tensor-e}
\\
&&\fit(b_1\ot b_2)=
\left\{
\begin{array}{ll}
\fit b_1\ot b_2 & {\mbox{ if }}\vp_i(b_1)>\vep_i(b_2)\\
b_1\ot\fit b_2  & {\mbox{ if }}\vp_i(b_1)\leq \vep_i(b_2).
\label{tensor-f}
\end{array}
\right.
\end{eqnarray}
Now, let us introduce the following results for a multi-tensor 
product of crystals:
\begin{thm}[\cite{KN}]\label{thm-KN}
Let $B(j)$ ($1\leq j\leq L$) be crystals. For $i\in I$ and 
$b_j\in B(j)$ ($1\leq j\leq L$) set $a_1:=0$ and 
\begin{eqnarray}
a_k&:=&\vep_i(b_k)-\sum_{1\leq j<k}\lan h_i,\wt(b_j)\ran
\qq (1< j\leq L),
\label{ak}\\
k_e&:=&\min\{k\in[1,L]\,|\,a_k=\max\{a_\nu\,|\,1\leq \nu\leq L\}\},
\label{ke}\\
k_f&:=&\max\{k\in[1,L]\,|\,a_k=\max\{a_\nu\,|\,1\leq \nu\leq L\}\}.
\label{kf}
\end{eqnarray}
By using these, the actions of $\eit$ and $\fit$ on the 
tensor product $b=b_1\ot\cd\ot b_L\in B(1)\ot\cd\ot B(L)$ are given as
\begin{eqnarray}
\eit(b)&=& b_1\ot\cd\ot\eit(b_{k_e})\ot\cd\ot b_L,\label{eit-ke}\\
\fit(b)&=& b_1\ot\cd\ot\fit(b_{k_f})\ot\cd\ot b_L.\label{fit-kf}
\end{eqnarray}
\end{thm}

\subsection{Kashiwar embedding}\label{poly-uqm}

First, define the crystal $B_i$ for $i\in I$: 
\begin{df}\label{Bi}
The crystal $B_i:=\{(x)_i\,|\,x\in \bbZ\}$ ($i\in I$) is defined as follows:
\begin{eqnarray}
&&\wt((x)_i)=x\al_i,\,\,
\vep_i((x)_i)=-x,\,\,
\vp_i((x)_i)=x,\,\,
\vep_j((x)_i)=\vp_j((x)_i)=-\ify\,\,\,(i\ne j),\label{bidata}\\
&&
\eit((x)_i)=(x+1)_i,\q\fit((x)_i)=(x-1)_i,
\q
\til e_j((x)_i)=\til f_j((x)_i)=0\,\,\,(i\ne j).\label{bieitfit}
\end{eqnarray}
\end{df}

\begin{thm}[\cite{K3}]\label{kashiwara-embed-thm}
Let $B_i=\{(x)_i\,|\,x\in\bbZ\}$ ($i\in I$) be the crystal as above.
\begin{enumerate}
\item
For any $i\in I$, we obtain the following strict embedding of crystals:
\begin{eqnarray}
\Psi_i:&&B(\ify)\hookrightarrow B(\ify)\ot B_i,\label{kashiwara-embed}\\
&&u_\ify\mapsto u_\ify\ot (0)_i.\nn
\end{eqnarray}
where $u_\ify:=1\,{\rm mod}\,qL(\ify)$ and $\Psi_i$ 
is called the {\it Kashiwara embedding}.
\item
For $b\in B(\ify)$, suppose $\Psi_i(b)=b_0\ot\fit^m(0)_i$. Then we get 
\[
\vep_i(b^*)=m,\qq\Psi_i(\fit^*b)=b_0\ot\fit^{m+1}(0)_i.
\]
\item
We also have 
\[
{\rm Im}(\Psi_i)=\{b\ot\fit^m(0)_i\,|\,\vep_i(b^*)=0,\,m\geq0\}
\]
\end{enumerate}
\end{thm}
Here note by \eqref{*=} if $b\in B(\ify)$, then $b^*\in B(\ify)$.
For the rest of this section, we fix an infinite sequence of indices
$\io=\cd,i_k\cd i_2i_1$ of $I$ such that
\begin{equation}
{\hbox{
$i_k\ne i_{k+1}$ and $\sharp\{k: i_k=i\}=\ify$ for any $i\in I$.}}
\label{seq-con}
\end{equation}
Iterating the Kashiwara embedding according to $\io$, we obtain the 
embedding : 
\begin{equation*}
\Psi_\io:
B(\ify)\mapright{\Psi_{i_1}}B(\ify)\ot B_{i_1}
\mapright{\Psi_{i_2}}B(\ify)\ot B_{i_2}\ot B_{i_1}
\cd\mapright{\Psi_{i_n}}B(\ify)\ot B_{i_n}\ot\cd \ot B_{i_2}\ot B_{i_1}
\longrightarrow \cd.
\end{equation*}
Set
\begin{equation}
\ZZ^{\ify}
:=\{(\cd,x_k,\cd,x_2,x_1): x_k\in\ZZ
\,\,{\rm and}\,\,x_k=0\,\,{\rm for}\,\,k\gg 0\};
\label{uni-cone}
\end{equation}
For any $\io$ as above, by the tensor product structure of 
$B_i$'s (Theorem \ref{thm-KN}), we can define a crystal structure
on $\ZZ^{\ify}$ and denote it by $\ZZ^{\ify}_{\io}$ 
(\cite[2.4]{NZ}):
For $x=(\cd,x_k,\cd,x_2,x_1)\in\bbZ^\ify$, set 
\begin{equation}
\sigma_k(x):=x_k+\sum_{j>k}\lan h_{i_k},\al_{i_j}\ran x_j
\label{sigma-def}
\end{equation}
Here note that the function $\sigma_k$ is formally an infinite sum, 
but, indeed, well-defined since $x_l=0$ for $l\gg0$.
For $i\in I$ define
\begin{eqnarray}
&&\TY(\wtil\sigma,i,)(x):=\max\{\sigma_k(x)\,|\,k\geq 1\,{\rm and}\,
i_k=i.\},\label{simgak}\\
&&\TY(\wtil M,i,)=\TY(\wtil M,i,)(x):=\{k\,|\,k\geq1,\,i_k=i,\,
\sigma_k(x)=\TY(\wtil\sigma,i,)(x).\},
\label{Mi}
\end{eqnarray}
which is also well-defined and $\TY(\wtil\sigma,i,)(x)\geq0$.
$\TY(M,i,)(x)$ is finite set 
if $\TY(\wtil\sigma,i,)(x)>0$. 

For $x\in\bbZ^\ify$, set 
\begin{equation}
\TY(\wtil m,i,f)=\TY(\wtil m,i,f)(x):=\min\,\TY(\wtil M,i,)(x),\q
\TY(\wtil m,i,e)=\TY(\wtil m,i,e)(x):=\max\,\TY(\wtil M,i,)(x)
\label{mfme}
\end{equation}
Then, $\TY(\wtil m,i,e)$ is finite if $\TY(\wtil\sigma,i,)(x)>0$.
Now, we shall define the Kashiwara operators :
\begin{eqnarray}
&&\fit(x)_k:=x_k+\del_{k,\TY(\wtil m,i,f)},\label{fit-Z-ify}\\
&&\TY(\wtil\sigma,i,)(x)>0\,{\rm then}\,
\eit(x)_k:=x_k-\del_{k,\TY(\wtil m,i,e)},\label{eit-Z-ify}\\
&&\TY(\wtil\sigma,i,)(x)=0\,{\rm then}\,\eit(x)=0.\nn
\end{eqnarray}
Finally, we define the crystal structure on $\bbZ^\ify$: 
\begin{equation}
\wt(x):=-\sum_{k\geq1}x_k\al_{i_k},\q
\vep_i(x):=\TY(\wtil\sigma,i,)(x),\q
\vp_i(x):=\lan h_i,\wt(x)\ran+\vep_i(x).
\label{wt-vep-vp-ify}
\end{equation}

We know that for any $b\in B(\ify)$ there exists $k>0$ and 
$a_k,\cd,a_1\geq0$ such that 
\[
\Psi_{i_k,i_{k-1},\cd,i_2,i_1}(b)=u_\ify\ot f_{i_k}^{a_k}(0)_{i_k}
\ot\cd\ot f_{i_1}^{a_1}(0)_{i_1}
\]
Here identifying $b\leftrightarrow (\cd,0,0,a_k,\cd,a_1)\in\bbZ^\ify_\io$, 
we obtain the strict embedding of crystals:
\begin{pro}[\cite{K3},\cite{NZ}]
\label{emb}
There is a unique strict embedding of crystals
\begin{equation}
\Psi_{\io}:B(\ify)\hookrightarrow \ZZ^{\ify}_{\io},
\label{psi}
\end{equation}
such that
$\Psi_{\io} (u_{\ify}) = (\cd,0,\cd,0,0)$, where 
$u_{\ify}\in B(\ify)$ is the vector corresponding to $1\in \uqm$.
It is also called the {\it  Kashiwara embedding} associated with $\io$.
\end{pro}

\subsection{Braid-Type Isomorphisms}

\begin{pro}[\cite{N2-2}]\label{braid}
There exist the following isomorphisms of crystals 
$\phi_{ij}^{(k)}$ 
($k=0,1,2,3$): 
\begin{align*}
&\phi_{ij}^{(0)}:B_i\ot B_j\mapright{\sim}B_j\ot B_i,&
{\bf a}_{ij}{\bf a}_{ji}=0, \\
&\phi_{ij}^{(1)}:B_i\ot B_j\ot B_i\mapright{\sim}B_j\ot B_i\ot B_j,&
{\bf a}_{ij}{\bf a}_{ji}=1, \\
&\phi_{ij}^{(2)}:(B_i\ot B_j)^{\ot 2}\mapright{\sim}(B_j\ot B_i)^{\ot 2},&
{\bf a}_{ij}{\bf a}_{ji}=2, \\
&\phi_{ij}^{(3)}:(B_i\ot B_j)^{\ot 3}\mapright{\sim}(B_j\ot B_i)^{\ot 3},&
{\bf a}_{ij}{\bf a}_{ji}=3.
\end{align*}
They satisfy $\phi_{ij}^{(k)}\circ \phi_{ji}^{(k)}={\rm id}$.
\end{pro}
These isomorphisms are called the {\it braid-type isomorphisms}.
Indeed, the explicit forms of $\phi_{ij}^{(0)}$ and $\phi_{ij}^{(1)}$ are 
\begin{eqnarray*}
&&\phi_{ij}^{(0)}((x)_i\ot(y)_j)=
(y)_j\ot(x)_i,\\
&&\phi_{ij}^{(1)}((x)_i\ot(y)_j\ot(z)_i)
=(\max(z,y-x))_j\ot(x+z)_i\ot (-\max(-x,z-y))_j.
\end{eqnarray*}
We omit the explicit forms of $\phi_{ij}^{(2)}$ and $\phi_{ij}^{(3)}$
(see \cite{N2-2}).

For a reduced word $\bfi=i_1i_2\cd i_k$ of a Weyl group element $w$, 
let us denote $B_{i_1}\ot B_{i_2}\ot\cd\ot B_{i_k}$ by $B_\bfi$, which 
will be called a {\it cellular crystal} associated with $\bfi$. 
\subsection{Kashiwara-Saito Characterization of the Crystal $B(\ify)$}

The following theorem has been introduced in \cite{KS} to 
characterize the crystal $B(\ify)$:
\begin{thm}[\cite{KS}]\label{KS-char}
Let $B$ be a crystal and $b_0\in B$ an element with weight 0. Assume the 
following seven conditions: 
\begin{enumerate}
\item
$\wt(B):=\{\wt(b)\,|\,b\in B\}\subset Q_-:=\bigoplus_{i\in I}\bbZ_{\leq0}\al_i$.
\item
$b_0$ is a unique element of $B$ with weight 0.
\item
$\vep_i(b_0)=0$ for any $i\in I$.
\item
$\vep_i(b)\in\bbZ$ for any $b\in B$ and $i\in I$.
\item
For any $i\in I$, there exists a strict embedding $\Psi_i:B\hookrightarrow 
B\ot B_i$.
\item
$\Psi_i(B)\subset B\ot\{\fit^n(0)_i\,|\,n\geq0\}$.
\item
For any $b\ne b_0$, there exists $i\in I$ such that 
$\Psi_i(b)=b'\ot \fit^n(0)_i$ with $n>0$.
\end{enumerate}
Then $B$ is isomorphic to $B(\ify)$.
\end{thm}

\begin{cor}\label{ks-cor}
If the 7th condition (vii) in Theorem \ref{KS-char} is replaced by the condition
\begin{enumerate}
\item[(vii')] If $b\ne b_0$, then there exists $i\in I$ such that 
$\eit b\ne0$,
\end{enumerate}
that is, 
under the conditions (i)--(vi) and (vii'), 
then we also obtain $B\cong B(\ify)$.
\end{cor}
{\sl Proof.}
To show this corollary, let us review the proof of Theorem \ref{KS-char}
following \cite{KS}.
First it follows from the condition (i), (ii) and (vi) that 
$\Psi(b_0)=b_0\ot (0)_i$ for any $i\in I$.
Let us show that for $b\ne b_0$, there exists $i$ such that $\eit b\ne0$.
If we assume $b\ne b_0$, by the condition (vi) and (vii) 
there exists $i\in I$ and $b'\in B$ 
such that $\Psi_i(b)=b'\ot \fit^n(0)_i$ with $n>0$.
If $b'= b_0$, then $\eit(\Psi_i(b))=b_0\ot \fit^{n-1}(0)_i\ne0$ 
and then $\eit(b)\ne0$.
If $b'\ne b_0$, then by the induction on the weight, there exists 
$j\in I$ such that $\til e_j(b')\ne0$ and then $\til e_j(b)\ne0$. 
Thus, for any $b\in B$ there exist $i_1,\cd,i_l\in I$ such that
$b=\til f_{i_1}\cd\til f_{i_l}b_0$, which implies that 
the crystal $B$ is connected as a crystal graph.

Assuming the conditions (i)--(vi) and (vii'), let us show the condition (vii).
If $b\ne b_0$, then by (i), (ii) and (vii'), 
 there exists a sequence $i_1,\ldots,i_k$ from 
$I$ such that $b_0=\til e_{i_k}\cd \til e_{i_1}b$
or equivalently, $b=\til f_{i_1}\cd \til f_{i_k}b_0$.
Then, it follows from (v) and (vi) that 
\begin{eqnarray*}
\Psi_{i_k}(b)&=&\til f_{i_1}\cd \til f_{i_k}\Psi_{i_k}(b_0)=
\til f_{i_1}\cd \til f_{i_k}(b_0\ot (0)_{i_k})\\
&=&\til f_{i_1}\cd \til f_{i_{k-1}}(b_0\ot \til f_{i_k}(0)_{i_k})
\in B\ot \{\til f_{i_k}^n(0)_{i_k}\,|\,n>0\},
\end{eqnarray*}
which shows the condition (vii).
%
Note that by the proof of Theorem \ref{KS-char}, we find that
the condition (vii) is only used to show that if $b\ne b_0$, 
then there exists $i\in I$ such that $\eit(b)\ne0$.\qed

\subsection{Monomial Realization of Crystals}
\label{mono-subsec}
Following \cite{K7,Nj}, 
we shall introduce the monomial realization of crystals.
For doubly-indexed variables $\{Y_{m,i}|i\in I, m\in\bbZ\}$, 
define the set of monomials
\[
 \cY:=\{Y=\prod_{m\in\bbZ,i\in I}
Y_{m,i}^{l_{m,i}}|l_{m,i}\in \bbZ\hbox{ and }l_{m,i}=0\text{
for all but finitely many }(m,i)\}.
\]
Fix a set of integers $p=(p_{i,j})_{i,j\in I,i\ne j}$ such that 
$p_{i,j}+p_{j,i}=1$, which we call a {\it sign}. 
Take a sign $p:=(p_{i,j})_{i,j\in I,i\ne j}$ and a
Cartan matrix $({\bf a}_{i,j})_{i,j\in I}$. For $m\in\bbZ,\,i\in I$ 
define the monomial 
\[
 A_{m,i}=Y_{m,i}Y_{m+1,i}\prod_{j\ne i}Y_{m+p_{j,i},j}^{{\bf a}_{j,i}}.
\]
Here, when we emphasize the sign $p$, we shall denote the
monomial
$A_{m,i}$ by $A^{(p)}_{m,i}$.
For any cyclic sequence of the indices
$\io=\cd (i_1i_2\cd i_n)(i_1i_2\cd i_n)\cd$
s.t. $\{i_1,\cd,i_n\}=I$,
we can associate the following sign $(p_{i,j})$ by:
\begin{equation}
 p_{i_a,i_b}=\begin{cases}
1&a<b,\\0&a>b.\end{cases}
\label{iaib}
\end{equation}
For example, if we take $\io=\cd (213)(213)\cd$, then we have
$p_{2,1}=p_{1,3}=p_{2,3}=1$ and $p_{1,2}=p_{3,1}=p_{3,2}=0$.
Thus, we can identify a cyclic sequence $\cd(i_1\cd i_n)(i_1\cd i_n)\cd$ 
with such $(p_{i,j})$.

For a monomial $Y=\prod_{m,i}Y_{m,i}^{l_{m,i}}$, 
set 
\begin{eqnarray*}
&&\hspace{-20pt}\wt(Y)=\sum_{i,m} l_{m,i}\Lm_i, \,\,
\vp_i(Y)=\operatorname{max}_{m\in\bbZ}\{\sum_{k\leq m}l_{k,i}\},\,\,
\vep_i(Y)=\vp_i(Y)-\wt(Y)(h_i)
=\max_{m\in\bbZ}\{-\sum_{k>m}l_{i,k}\},\\
&&\hspace{-20pt}\fit(Y)=\begin{cases}
A_{n_f^{(i)}(Y),i}^{-1}\cdot Y&\text{ if }\vp_i(Y)>0,\\
0&\text{ if }\vp_i(Y)=0,
\end{cases}\q\q
\eit(Y)=\begin{cases}
A_{n_e^{(i)}(Y),i}\cdot Y&\text{ if }\vep_i(Y)>0,\\
0&\text{ if }\vep_i(Y)=0,
\end{cases}\\
&&\text{where }n_f^{(i)}(Y)=\min\{n|\vp_i(Y)=\sum_{k\leq n}l_{k,i}\},\q
n_e^{(i)}(Y)=\max\{n|\vp_i(Y)=\sum_{k\leq n}l_{k,i}\}.
\end{eqnarray*}

\begin{thm}[\cite{K7,Nj}]
\begin{enumerate}
\item
In the above setting, $\cY$ is a crystal, which is denoted by $\cY(p)$.
\item
If $Y\in\cY(p)$ satisfies $\vep_i(Y)=0$ (resp. $\vp_i(Y)=0$) for any $i\in I$, 
then the connected component containing $Y$ is isomorphic to
$B(\wt(Y))$ (resp. $B(w_0\wt(Y))$), 
where we call such monomial $Y$ a highest (resp. lowest) monomial.
\end{enumerate}
\end{thm}
By the above crystal structure of monomials, we know that 
for any $k\in \bbZ,\,\,i\in I$
the monomial $Y_{k,i}$ is a highest monomial with the weight $\Lm_i$. 
Thus, we can define the embedding of crystal 
$\mpp(m,k,i)$ $(i\in I,k\in\bbZ)$ by 
\begin{eqnarray}
 \mpp(m,k,i):B(\Lm_i)&\hookrightarrow& \cY(p)
\label{embed-mki}\\
u_{\Lm_i}&\mapsto& Y_{k,i} \nn
\end{eqnarray}

\renewcommand{\thesection}{\arabic{section}}
\section{Geometric crystals}
\setcounter{equation}{0}
\renewcommand{\theequation}{\thesection.\arabic{equation}}

The basic references for this section are \cite{BK,BK2,N}.
We assume the setting as in Sect.2.
\subsection{Definition of geometric crystals}

\begin{df}
\label{def-gc}
Let $X$ be an irreducible algebraic variety over $\bbC$, 
$\gamma_i$, $\vep_i$ 
$(i\in I)$ rational functions on $X$, and 
$e_i:\bbC^\times\times X\to X$ ($(c,x)\mapsto e_i^c(x)$) 
a rational unital action of multiplicative group $\bbC^\times$
($i\in I$).
A quadruple $\chi=(X,\{e_i\}_{i\in I},\{\gamma_i,\}_{i\in I},
\{\vep_i\}_{i\in I})$ is a 
$\ge$-{\it geometric crystal} 
if
\begin{enumerate}
\item
The rational functions  $\{\gamma_i\}_{i\in I}$ satisfy
$\gamma_j(e^c_i(x))=c^{{\bf a}_{ij}}\gamma_j(x)$ for any $i,j\in I$.
\item
$e_i$ and $e_j$ satisfy the following {\it Verma relations}:
\[
 \begin{array}{lll}
&\hspace{-20pt}e^{c_1}_{i}e^{c_2}_{j}
=e^{c_2}_{j}e^{c_1}_{i}&
{\rm if }\,\,{\bf a}_{ij}={\bf a}_{ji}=0,\\
&\hspace{-20pt} e^{c_1}_{i}e^{c_1c_2}_{j}e^{c_2}_{i}
=e^{c_2}_{j}e^{c_1c_2}_{i}e^{c_1}_{j}&
{\rm if }\,\,{\bf a}_{ij}={\bf a}_{ji}=-1,\\
&\hspace{-20pt}
e^{c_1}_{i}e^{c^2_1c_2}_{j}e^{c_1c_2}_{i}e^{c_2}_{j}
=e^{c_2}_{j}e^{c_1c_2}_{i}e^{c^2_1c_2}_{j}e^{c_1}_{i}&
{\rm if }\,\,{\bf a}_{ij}=-2,\,
{\bf a}_{ji}=-1,\\
&\hspace{-20pt}
e^{c_1}_{i}e^{c^3_1c_2}_{j}e^{c^2_1c_2}_{i}
e^{c^3_1c^2_2}_{j}e^{c_1c_2}_{i}e^{c_2}_{j}
=e^{c_2}_{j}e^{c_1c_2}_{i}e^{c^3_1c^2_2}_{j}e^{c^2_1c_2}_{i}
e^{c^3_1c_2}_je^{c_1}_i&
{\rm if }\,\,{\bf a}_{ij}=-3,\,
{\bf a}_{ji}=-1.
\end{array}
\]
\item
The rational functions $\{\vep_i\}_{i\in I}$ satisfy
$\vep_i(e_i^c(x))=c^{-1}\vep_i(x)$ and 
$\vep_i(e_j^c(x))=\vep_i(x)$ if ${\bf a}_{i,j}={\bf a}_{j,i}=0$.
\end{enumerate}
\end{df}
Here we call $e_i$'s {\it geometric crystal operators}
and define $\vp_i:=\vep_i\cdot\gamma_i$.

\begin{df}
Let $\chi=(X,\{e_i\}_{i\in I},\{\gamma_i,\}_{i\in I},
\{\vep_i\}_{i\in I})$ be a
$\ge$-geometric crystal. 
\begin{enumerate}
\item
A pair $(\chi,f)$ is a 
$\ge$-{\it decorated geometric crystal} if $f$ is a rational function on 
$X$ satisfying
\begin{equation}
f(e_i^c(x))=f(x)+{(c-1)\vp_i(x)}+{(c^{-1}-1)\vep_i(x)},
\label{f}
\end{equation}
for any $i\in I$, $c\in \bbC^\times$ 
and $x\in X$.
We call the function $f$ a {\it decoration} or a {\it potential} on $X$.
\item 
A pair $(\chi,f)$ is a 
$\ge$-{\it upper (resp. lower) half-decorated geometric crystal} 
if $f$ is a rational function on 
$X$ satisfying
\begin{equation}
f(e_i^c(x))=f(x)+{(c^{-1}-1)\vep_i(x)}\q
({\rm resp.\,}f(e_i^c(x))=f(x)+{(c-1)\vp_i(x)})
\label{half-f}
\end{equation}
for any $i\in I$, $c\in\bbC^\times$  and $x\in X$.
We call the function $f$ an 
{\it upper ({\rm resp.} lower) half-decoration} or an
{\it upper ({\rm resp.} lower)  half potential} on $X$.
\end{enumerate}
\end{df}

\nd
\begin{rem}
The definitions of $\vep_i$ and $\vp_i$ are different from the ones in 
e.g., \cite{BK2} since we adopt the definitions following
\cite{KNO,KNO2}. Indeed, if we flip $\vep_i\to \vep_i^{-1}$ and 
$\vp_i\to \vp_i^{-1}$, they coincide with ours.
\end{rem}
%
%
%

\subsection{Positive structure and Tropicalization of Geometric Crystals}
\label{subsec-posi}

For this subsection, see {\it e.g.,} \cite{BK2}, \cite[3.3]{N4}.

First, for a rational function $f(x)\in\bbC(x)$ let us define its degree 
${\rm deg}(f(x)):={\rm deg}(g(x))-{\rm deg}(h(x))$ where 
$g(x),h(x)\in \bbC[x]$ such that $f(x)=\frac{g(x)}{h(x)}$. Furthermore, 
a rational function $f(x)$ is positive if there exist polynomials 
$g(x),h(x)\in\bbR_{\geq 0}[x]$ such that $f(x)=\frac{g(x)}{h(x)}$. 
Next, let us define the valuation 
$\bbV:\bbC(x)\to {}^t\bbZ:=\bbZ\sqcup\{-\ify\}$ by 
\begin{equation}
\bbV(f(x)):=\begin{cases}-{\rm deg}(f(x^{-1}))&f\not\equiv0,\\
-\ify&f\equiv0,
\end{cases}
\label{mapV}
\end{equation}
where ${}^t\bbZ$ is the tropical semi-field with multiplication
``$+$'' and addition ``$\min$''. This valuation satisfies 
$\bbV(f_1\cdot f_2)=\bbV(f_1)+\bbV(f_2)$ and 
$\bbV(f_1/f_2)=\bbV(f_1)-\bbV(f_2)$ for $f_1,f_2\in\bbC(x)$ and also satisfies
$\bbV(f_1+f_2)=\min(\bbV(f_1),\bbV(f_2))$ 
for positive rational functions $f_1,f_2$.

\begin{df}
For algebraic tori $T$, $T'$ and a rational map 
$f:T\to T'$, define 
$\what f:X_*(T)\to X_*(T')$ by 
\begin{equation}
\lan \chi,\what f(\xi)\ran=\bbV(\chi\circ f\circ \xi)\q
(\chi\in X^*(T'),\q \xi\in X_*(T)).
\label{whatf}
\end{equation}
\end{df}

Let $\cT_+$ be the category whose objects are algebraic tori over
$\bbC$ and whose morphisms are positive rational maps.
Then, we define the {\it tropicalization} functor
${\mathcal TR}: \cT_+\longrightarrow {\mathfrak Set}(=$
the category of sets)  by 
\begin{eqnarray*}
{\mathcal TR}&:&\cT_+\longrightarrow {\mathfrak Set}\\
&&T \mapsto X_*(T)\\
&&(f:T\to T')\mapsto (\what f:X_*(T)\to X_*(T'))
\end{eqnarray*}
\begin{rem}
\begin{enumerate}
\item
The functor ${\mathcal TR}$ has been denoted by ${\mathcal ID}$ in 
\cite{N,N3,N6}.
\item
This definition is slightly different from the one in \cite{KN,N,N3,N4}
since here we take a different valuation from the one in those references.
\end{enumerate}
\end{rem}
For an algebraic torus $T$, let 
$\theta:T\rightarrow X$ be a positive structure 
(see \cite{BK2},\cite[3.3]{N4}) on 
a (half) decorated geometric crystal $\chi=(X,\{e_i\}_{i\in I},
\{{\rm wt}_i\}_{i\in I},
\{\vep_i\}_{i\in I},f)$, namely, 
all rational maps
$e_{i,\theta}:=\theta^{-1}\circ e_i\circ\theta
:\bbC^\tm \tm T\rightarrow T$ and
$f\circ\theta,\gamma_i\circ \theta,\vep_i\circ\theta:T\ra \bbC$
are positive maps (see \cite{BK2},\cite[3.3]{N4}).
Applying the functor ${\mathcal TR}$ 
to these positive rational morphisms, 
we obtain piece-wise linear maps 
\begin{eqnarray*}
&&\til e_i:={\mathcal TR}(e_{i,\theta}):
\ZZ\tm X_*(T) \rightarrow X_*(T),\qq
{\rm wt}_i:={\mathcal TR}(\gamma_i\circ\theta):
X_*(T)\rightarrow \bbZ,\\
&&\wtil\vep_i:={\mathcal TR}(\vep_i\circ\theta):
X_*(T)\rightarrow \bbZ,\qq
\wtil f:= {\mathcal TR}(f \circ\theta):
X_*(T)\rightarrow \bbZ,
\end{eqnarray*}
where ${\rm wt}_i(b):=\lan h_i,{\rm wt}(b)\ran$ for $b\in X_*(T)$.
%
\begin{thm}[\cite{BK,BK2,N}]
For any $\ge$-geometric crystal 
$\chi=(X,\{e_i\}_{i\in I},\{\gamma_i\}_{i\in I},
\{\vep_i\}_{i\in I})$ and positive structure
$\theta:T\rightarrow X$, the associated crystal 
${\mathcal TR}_{\theta,T}(\chi)=
(X_*(T),\{e_i\}_{i\in I},\{{\rm wt}_i\}_{i\in I},
\{\wtil\vep_i\}_{i\in I})$ 
is a Langlands dual $\ge^\vee$-Kashiwara's crystal.
\end{thm}
{\sl Remark.}
The definition of $\wtil\vep_i$ is different from the one in 
\cite[6.1]{BK2} since our definition of $\vep_i$ corresponds to 
$\vep_i^{-1}$ in \cite{BK2}.

For a positive upper/lower-half decorated geometric crystal 
$\cX=((X,\{e_i\}_{i\in I},\{\gamma_i\}_{i\in I},
\{\vep_i\}_{i\in I},f),\theta,T)$, set 
\begin{equation}
\wtil B_{f,\theta}:=\{\wtil x\in X_*(T)(=\bbZ^{\dim(T)})
|\wtil f(\wtil x)\geq0\},
\label{btil}
\end{equation}
where we define $\eit^n(\til x):={\mathcal TR}(e_{i,\theta})(n,\til x)$
($n\in\bbZ$) and 
if $\eit(\wtil x)\not\in \wtil B_{f,\theta}$, then 
$\eit(\wtil x)=0$. Set 
$\bbB_{f,\theta}:=(\wtil B_{f,\theta},\wtil{\wt}_i|_{\wtil B_{f,\theta}},
\wtil{\vep}_i|_{\wtil B_{f,\theta}},\wtil{e}_i|_{\wtil B_{f,\theta}})_{i\in I}$.
\begin{pro}\label{ul-normal}
For a positive upper/lower half decorated geometric crystal \\
$\cX=((X,\{e_i\}_{i\in I},\{\gamma_i\}_{i\in I},
\{\vep_i\}_{i\in I},f),\theta,T)$, the
 quadruple $\bbB_{f,\theta}$ is a $\ge^\vee$-upper/lower normal crystal.
\end{pro}
{\sl Proof.}
The case where $f$ is a decoration has been shown in \cite{BK2}.
We shall address the case  where $f$ is an upper half potential 
in a similar way.
The lower case is also shown similarly.
For $i\in I$ set $f_0:=f-\vep_i$ and then we get
\begin{equation*}
f(e_i^c x)=f_0(x)+c^{-1}\vep_i(x)
\end{equation*}
Since $e_i^c$ and $f$ are positive and the function $f_0$ does not depend on 
$c$, by tropicalization we get for an integer $n$, 
\begin{equation}
\wtil f(\eit^n b)=\min(\wtil f_0(b),\wtil \vep_i(b)-n).
\label{f0}
\end{equation}
If we take $n=0$, then we have
$\wtil f(b)=\min(\wtil f_0(b),\wtil \vep_i(b))$.
Therefore, we have 
\[
b\in \wtil B_{f,\theta}\Leftrightarrow \wtil f(b)\geq0
\Leftrightarrow \wtil f_0(b)\geq0, \,\,\wtil\vep_i(b)\geq0.
\]
Then we get $\wtil f_0(b)\geq0$ for $b\in\wtil B_{f,\theta}$. Thus 
by \eqref{f0} we obtain 
\begin{equation*}
\eit^n(b)\in \wtil B_{f,\theta}\Leftrightarrow 
\wtil\vep_i(b)\geq n.
\end{equation*}
Hence,
\begin{equation*}
\eit^{n+1}(b)\not\in \wtil B_{f,\theta}\Leftrightarrow 
\wtil\vep_i(b)<n+1.
\end{equation*}
These results imply that $B_{f,\theta}$ is a upper normal crystal, 
that is,
\begin{equation*}
\wtil\vep_i(b)=\max\{n\,|\,\eit^n(b)\in \wtil B_{\theta, f}\}.
\end{equation*}\qed

\renewcommand{\thesection}{\arabic{section}}
\section{Half potentials of the geometric crystal on $B^-_{w_0}$}
\setcounter{equation}{0}
\renewcommand{\theequation}{\thesection.\arabic{equation}}

\subsection{Geometric crystal on $B^-_w$ and the associated cellular crystals}
\label{BBB}
For a Weyl group element $w\in W$, define 
$B^-_w:=B^-\cap U\ovl w U$, which is called here a Bruhat cell of $B^-$ 
associated with $w\in W$.
Let $\gamma_i:B^-_w\to\bbC^\times$ be the rational function defined by 
\begin{equation}
\gamma_i:B^-_w\,\,\hookrightarrow \,\,\,
B^-\,\,\mapright{\sim}\,\,T\times U^-\,\,
\mapright{\rm proj}\,\,\, T\,\,\,\mapright{\al_i^\vee}\,\,\,\bbC^\times.
\label{gammai}
\end{equation}

For any $i\in I$, there exists the natural projection 
$pr_i:B^-\to B^-\cap \phi(SL_2)$. Hence, 
for any $x\in B^-_w$ there exists unique
       $v=\begin{pmatrix}b_{11}&0\\b_{21}&b_{22}\end{pmatrix}
\in SL_2$ such that 
$pr_i(x)=\phi_i(v)$. Using this fact, we define 
the rational function $\vep_i$ on $B^-_w$ as in \cite{N4}:
\begin{equation}
\vep_i(x)=\frac{b_{22}}{b_{21}}\q(x\in B^-_w).
\label{vepi}
\end{equation}
The unital rational $\bbC^\times$-action $e_i$ on $B^-_w$ is defined by
\begin{equation}
e_i^c(x):=x_i\left((c-1)\vp_i(x)\right)\cdot x\cdot
x_i\left((c^{-1}-1)\vep_i(x)\right)\qq
(c\in\bbC^\times,\,\,x\in B^-_w),
\label{ei-action}
\end{equation}
if $\vep_i(x)$ is well-defined, that is, $b_{21}\ne0$, 
and define $e_i^c(x)=x$ if $b_{21}=0$.\\
{\sl Remark.} The definition (\ref{vepi}) is different from the one in 
\cite{BK2}. Indeed, if for \eqref{vepi} we take $\vep_i(x)=b_{21}/b_{22}$, 
then it coincides with
the one in \cite{BK2}.
\begin{pro}[\cite{BK2}]
For any $w\in W$,
the quadruple $\chi:=(B^-_w,\{e_i\}_i,\{\gamma_i\}_i,\{\vep_i\}_i)$
is a geometric crystal, where 
$\gamma_i$, $\vep_i$ and $e_i$ are given in (\ref{gammai}), 
(\ref{vepi}) and (\ref{ei-action}) respectively.
\end{pro}

\def\ld{\ldots}
For a Weyl group element $w\in W$, let 
$\bfi=i_1\ld i_k$ be one of its reduced expressions and 
define the positive structure 
$\Theta^-_\bfi:(\bbC^\times)^k\longrightarrow B^-_w$ by 
\[
 \Theta^-_\bfi(c_1,\cd,c_k):=\pmby_{i_1}(c_1)\cd \pmby_{i_k}(c_k),
\]
where $\pmby_i(c)=y_i(c)\al_i^\vee(c^{-1})$, 
which is different from $Y_i(c)$ in 
\cite{N,N2,KNO,KNO2}. Indeed, $Y_i(c)=\pmby_i(c^{-1})$.

Now, let us  describe the positive geometric crystal
structure on $B^-_{w}$ explicitly following \cite{N4}.
Taking a reduced word $\bfi=i_1i_2\cd i_k$ 
for the  Weyl group element $w$, we set
\[
B^-_\bfi:=\{\pmby_{i_1}(c_1)\cd \pmby_{i_N}(c_k)\,|\,
c_1,\ldots,c_k\in\bbC^\times\}.
\]
where note that $B^-_{w}$ and $B^-_\bfi$ are birationally
isomorphic to each other.
\begin{pro}[\cite{N4}]
The action $e^\xi_i$ $(\xi\in\bbC^\times)$ on 
$\Theta^-_{\bfi}(c_1,\cd,c_k)$ is given by
\[
e_i^\xi(\Theta^-_{\bfi}(c_1,\cd,c_k))
=\Theta^-_{\bfi}(c'_1,\cd,c'_k), 
\]
where for $j=1,2,\cd, k$, 
\begin{equation}
c'_j\seteq 
c_j\cdot \frac{\displaystyle \sum_{1\leq m< j,\,i_m=i}
{\xi\cdot c_1^{{\bf a}_{i_1,i}}\cd c_{m-1}^{{\bf a}_{i_{m-1},i}}c_m}
+\sum_{j\leq m\leq k,\,i_m=i} 
{c_1^{{\bf a}_{i_1,i}}\cd c_{m-1}^{{\bf a}_{i_{m-1},i}}c_m}}
{\displaystyle \sum_{1\leq m\leq j,\,i_m=i} 
{\xi\cdot c_1^{{\bf a}_{i_1,i}}\cd c_{m-1}^{{\bf a}_{i_{m-1},i}}c_m}+
\mathop\sum_{j< m\leq k,\,i_m=i}  
{c_1^{{\bf a}_{i_1,i}}\cd c_{m-1}^{{\bf a}_{i_{m-1},i}}c_m}},
\label{eici}
\end{equation}
where note that in \eqref{eici}, if $i_j\ne  i$, then $c'_j=c_j$.
The explicit forms of 
rational functions $\vep_i$ and $\gamma_i$ are:
\begin{equation}
 \vep_i(\Theta^-_{\bfi}(c))=
\left(\sum_{1\leq m\leq k,\,i_m=i} \frac{1}
{c_mc_{m+1}^{{\bf a}_{i_{m+1},i}}\cd c_{k}^{{\bf a}_{i_{k},i}}}\right)^{-1},\q
\gamma_i(\Theta^-_{\bfi}(c))
=\frac{1}{c_1^{{\bf a}_{i_1,i}}\cd c_k^{{\bf a}_{i_k,i}}}.
\label{th-vep-gamma}
\end{equation}
\end{pro}

Applying the tropicalization functor to this positive geometric crystal
$B^-_w\cong B^-_\bfi$, we shall obtain 
the following Langlands-dual free crystal 
\[
{\mathcal TR}(B^-_\bfi)=\{(x_1,\cd,x_k)\,|\,x_j\in\bbZ\}
\]
 given as follows:
For $x=(x_1,\cd,x_k)\in {\mathcal TR}(B^-_\bfi)$ and $i\in I$, we
 have
\begin{equation}
\eit^n(x)=(x'_1,\cd,x'_k)\q(n\in\bbZ),
\end{equation}
where 
\begin{eqnarray}
&& x'_j=x_j+\min\left(\min_{1\leq m<j,i_m=i}(n+x_m+\sum_{l=1 }^{m-1}
 {\bf a}_{i_l,i}x_l),
\min_{j\leq m\leq k,i_m=i}(x_m+\sum_{l=1 }^{m-1} {\bf a}_{i_l,i}x_l)\right)
\label{eitn} \\
&&\qq -\min\left(\min_{1\leq m\leq j,i_m=i}(n+x_m+\sum_{l=1 }^{m-1}
 {\bf a}_{i_l,i}x_l),
\min_{j<m\leq k,i_m=i}(x_m+\sum_{l=1 }^{m-1} {\bf a}_{i_l,i}x_l)\right),
\nn
\\
&&\wt_i(x)=-\sum_{l=1}^k {\bf a}_{i_l,i}x_l\qq
(\wt(x):=-\sum_{l=1}^k x_l\al_{i_l}),
\label{wt-ud}\\
&&\vep_i(x)=\max_{1\leq m\leq
k,i_m=i}(x_m+\sum_{l=m+1}^{k}{\bf a}_{i_l,i}x_l),
\label{vep-ud}
\end{eqnarray}
where $\al_{i_l}$'s in \eqref{wt-ud} are simple roots of the Langlands dual 
simple Lie algebra ${}^L\ge$.
\begin{rem}
The formula (3.13) in \cite{N4} is incorrect. 
Here the formula \eqref{eitn} is correct.
\end{rem}
Then the crystal ${\mathcal TR}(B^-_\bfi)$
 is isomorphic to the Langlands dual crystal
$(B_{i_k}\ot \cd\ot B_{i_2}\ot B_{i_1})^\vee=
B_{i_k}^\vee\ot\cd\ot B_{i_1}^\vee$ 
by the correspondence
\begin{equation}
\begin{array}{ccc}
{\mathcal TR}(B^-_\bfi)
&\longleftrightarrow &(B_{i_k}\ot \cd\ot B_{i_2}\ot B_{i_1})^\vee,\\
(x_1,x_2,\cd,x_k)&\longleftrightarrow&(-x_k)_{i_k}\ot\cd\ot (-x_2)_{i_2}\ot
(-x_1)_{i_1},
\end{array}
\label{i-i}
\end{equation}
where $B_i^\vee:=\{(x)_i\,|\,x\in\bbZ\}$ is the $\LL(\ge)$-crystal associated with a simple root $\LL(\al_i)$ of $\LL(\ge)$.
We call the crystal $B_\bfi:=B_{i_1}\ot\cd\ot
B_{i_k}$ a {\it cellular crystal} associated with the 
reduced word $\bfi=i_1i_2\cd i_k$.

The following lemma is needed to prove Proposition \ref{optimal}.
\begin{lem}\label{BBiso}
Let $i_1i_2\cd i_k$ be a reduced word for a Weyl group element $w$.
Suppose that applying the braid-type isomorphisms properly to 
$(x_1)_{i_1}\ot (x_2)_{i_2}\ot\cd\ot(x_k)_{i_k}$, one gets 
$(y_1)_{i_1}\ot (y_2)_{i_2}\ot\cd\ot(y_k)_{i_k}$. 
Then we find that $x_1=y_1, x_2=y_2,\cd,x_k=y_k$.
\end{lem}
{\sl Proof.}
Set $Y_i(c):=y_i(c^{-1})\al_i^\vee(c))=\pmby_i(c^{-1})$.
The braid-type isomorphisms are obtained from the following equations
by considering the tropicalization method (\cite{N}):
\begin{equation}
\begin{split}
&Y_i(c_1)Y_j(c_2)=Y_j(d_1)Y_i(d_2)\q({\bf a}_{ij}{\bf a}_{ji}=0),\\
&Y_i(c_1)Y_j(c_2)Y_i(c_3)=Y_j(d_1)Y_i(d_2)Y_j(d_3)\q({\bf a}_{ij}{\bf a}_{ji}=1),\\
&Y_i(c_1)Y_j(c_2)Y_i(c_3)Y_j(c_4)=Y_j(d_1)Y_i(d_2)Y_j(d_3)Y_i(d_4)
\q({\bf a}_{ij}{\bf a}_{ji}=2),\\
&Y_i(c_1)Y_j(c_2)Y_i(c_3)Y_j(c_4)Y_i(c_5)Y_j(c_6)
=Y_j(d_1)Y_i(d_2)Y_j(d_3)Y_i(d_4)Y_j(d_4)Y_i(d_6)
\q({\bf a}_{ij}{\bf a}_{ji}=3),\\
\end{split}
\label{2346move}
\end{equation}
where $d_j$'s are positive rational functions in $c_j$'s.
Here note that the tropicalization process \cite{N} and 
here in \ref{subsec-posi} are different, but the resulting formula
for crystals are correct. Indeed, for example, by the present 
tropicalization of the 3rd formula in \eqref{2346move} we obtain the 
isomorphism of crystals 
\begin{equation}
\begin{split}
B_j\ot B_i\ot B_j\ot B_i&\mapright{\sim}B_i\ot B_j\ot B_i\ot B_j,\\
(-\til c_4)_j\ot (-\til c_3)_i\ot 
(-\til c_2)_j\ot (-\til c_1)_i&\mapsto
(-\til d_4)_i\ot (-\til d_3)_j\ot 
(-\til d_2)_i\ot (-\til d_1)_j,
\end{split}\label{bb}
\end{equation}
where $\til c_j,\til d_j$ are the tropicalization of $c_j,d_j$ respectively.
The formula \eqref{bb}
gives the braid-type isomorphism $\phi^{(3)}_{ji}$ in \eqref{braid}.

Thus if we apply these formula by the same way as in the lemma to
$Y_{i_1}(c_1)Y_{i_2}(c_2)\cd Y_{i_k}(c_k)$ and get 
$Y_{i_1}(d_1)Y_{i_2}(d_2)\cd Y_{i_k}(d_k)$, then they coincide with 
each other, that is, 
\[
Y_{i_1}(c_1)Y_{i_2}(c_2)\cd Y_{i_k}(c_k)
=Y_{i_1}(d_1)Y_{i_2}(d_2)\cd Y_{i_k}(d_k).
\]
Since the factorization 
\begin{equation}
(\bbC^\times)^k\mapright{\sim}B^-_w, 
\q(c_1,\cd,c_k)\mapsto 
Y_{i_1}(c_1)\cd Y_{i_k}(c_k),
\label{12cdk}
\end{equation}
is a birational map, it is injective and then we obtain 
$c_j=d_j$($j=1,2,\cd,k$).
Tropicalizing \eqref{12cdk} and set $x_j:={\mathcal TR}(c_j)$ and 
$y_j:={\mathcal TR}(d_j)$, we have $x_j=y_j$($j=1,2,\cd,k$).\qed

The following proposition will be needed later.
\begin{pro}\label{optimal}
\begin{enumerate}
\item
Take a reduced word $\bfi=i_1i_2\cd i_k$ of a Weyl group element $w$. 
Now for two reduced words ${\bf j}=j_1j_2\cd j_k$ and 
${\bf l}=l_1l_2\cd l_k$ of $w$ satisfying $j_1=l_1$ (resp. $j_k=l_k$). 
We assume that by applying the braid-type isomorphisms to
${\bf x}=(x_1)_{i_1}\ot (x_2)_{i_2}\ot\cd\ot (x_k)_{i_k}\in B_\bfi$ 
properly in 
two ways we obtain 
${\bf y}=(y_1)_{j_1}\ot (y_2)_{j_2}\ot\cd\ot (y_k)_{j_k}
\in B_{\bf j}$ and 
${\bf z}=(z_1)_{l_1}\ot (z_2)_{l_2}\ot\cd\ot (z_k)_{l_k}\in B_{\bf l}$. 
Then we have $y_1=z_1$ (resp. $y_k=z_k$).
Here, let us define the map 
\begin{equation}
\omega_i\,\, ({\rm resp.}\,\,\omega'_i)
:B_{\bf i}\to\bbZ,\q{\rm by}\,\, \omega_i({\bf x})=y_1\q
({\rm resp.}\,\, \omega'_i({\bf x})=y_k).
\label{omegai}
\end{equation}
\item
For any $i\in I$ and any reduced word $\bfi=i_1i_2\cd i_k$ of a
Weyl group element $w\in W$
take another reduced word ${\bf j}=j_1j_2\cd j_l$ of 
$w\in W$ such that $j_1=i$ and  define the following map 
\begin{eqnarray*}
\xi^{(i)}_{\bfi,\bf j;\wtil\phi}:
&&B_{\bfi}\,\,\mapright{\wtil\phi}\,\,B_{\bf j}\,\,
\mapright{r}\,\,B_{\bf j}\,\,\mapright{{\wtil\phi}^{-1}}\,\,
B_{\bfi}\\
&&(x_1)_{i_1}\ot(x_2)_{i_2}\ot\cd\ot (x_k)_{i_k}\mapsto
(y_1)_{i}\ot(y_2)_{j_2}\ot\cd\ot (y_k)_{j_k}\\
&&\mapsto
(0)_{i}\ot(y_2)_{j_2}\ot\cd\ot (y_k)_{j_k}\mapsto
(x'_1)_{i_1}\ot(x'_2)_{i_2}\ot\cd\ot(x'_k)_{i_k},
\end{eqnarray*}
where $\wtil \phi$ is a proper composition of the braid-type isomorphisms
and the map $r:B_{\bf j}\to B_{\bf j}$ is defined by 
$(y_1)_{i}\ot(y_2)_{j_2}\ot\cd\ot (y_k)_{j_k}
\mapsto
(0)_{i}\ot(y_2)_{j_2}\ot\cd\ot (y_k)_{j_k}$.
Then $\xi_{\bf i,j;\wtil\phi}$ 
does not depend on 
the composition of braid-type isomorphisms $\wtil\phi$ 
and the choice of ${\bf j}$.
That is, if we take another reduced word ${\bf l}=l_1l_2\cd l_k$ of $w$
such that $j_1=l_1=i$  and a composition of the braid-type isomorphisms
$\wtil\phi':B_\bfi\to B_{\bf l}$, then we obtain 
$\xi^{(i)}_{\bf i,j;\wtil\phi}=\xi^{(i)}_{\bf i,l;\wtil\phi'}$.
Thus, we denote this map by $\xi^{(i)}_{\bfi}$.
\end{enumerate}
\end{pro}
{\sl Proof.} (i) 
We shall show the case $j_1=l_1$ since the case $j_k=l_k$ is
similar.
Since $j_1=l_1=:j$, one knows that both ${\bf j}':=j_2\cd j_k$ and 
${\bf l}'=l_2\cd l_k$ are
the reduced words for $s_j w$. Then, applying the usual braid relations on 
${\bf j}'$ properly, we obtain ${\bf l}'$. 
Then applying the braid-type isomorphisms to 
${\bf y}':=(y_2)_{j_2}\ot\cd\ot (y_k)_{j_k}$ in the same way, we get 
${\bf z}':=(z_2')_{l_2}\ot\cd\ot (z_k')_{l_k}$. Since 
both 
$(z_1)_{l_1}\ot (z_2)_{l_2}\ot\cd\ot (z_k)_{l_k}$ and 
$(y_1)_{l_1}\ot (z_2')_{l_2}\ot\cd\ot (z_k')_{l_k}$ are obtained from 
${\bf x}=(x_1)_{i_1}\ot (x_2)_{i_2}\ot\cd\ot (x_k)_{i_k}$, 
it follows from Lemma \ref{BBiso} above that 
$z_1=y_1$, $z_m=z_m'$ ($m=2,3,\cd,k$). 

\nd
(ii)
For $c_1,\ldots,c_k\in\bbC^\times$, take $d_j,d'_j\in\bbC$ ($j=1,2,\ldots,k$)
such that
\[
Y_{i_k}(c_k)\cd Y_{i_2}(c_2)Y_{i_1}(c_1)
=Y_{j_k}(d_k)\cd Y_{j_2}(d_2)Y_{i}(d_1)
=Y_{l_k}(d'_k)\cd Y_{l_2}(d'_2)Y_{i}(d'_1),
\]
where note that by (i) one has $d_1=d'_1$ and then 
\[
Y_{j_k}(d_k)\cd Y_{j_2}(d_2)=
Y_{l_k}(d'_k)\cd Y_{l_2}(d'_2).
\]
Furthermore, one has $c'_1,\ldots,c'_k\in\bbC^\times$ such that 
\[
Y_{j_k}(d_k)\cd Y_{j_2}(d_2)Y_i(1)=
Y_{l_k}(d'_k)\cd Y_{l_2}(d'_2)Y_i(1)
=Y_{i_k}(c'_k)\cd Y_{i_2}(c'_2)Y_{i_1}(c'_1).
\]
Here setting $x_j:=-\wtil c_j$, $y_j=-\wtil d_j$, $z_j=-\wtil d'_j$, 
$x'_j=-\wtil c'_j$ ($j=1,2,\ldots, k$), we get 
\begin{eqnarray*}
&&\wtil\phi(x_1,\ldots,x_k)=(y_1,\ldots,y_k),\\
&&\wtil\phi'(x_1,\ldots,x_k)=(z_1,\ldots,z_k),\\
&&\wtil\phi^{-1}(0,y_2,\ldots,y_k)=(x'_1,\ldots,x'_k)
=\wtil{\phi'}^{-1}(0,z_2,\ldots,z_k),
\end{eqnarray*}
which implies that $\xi^{(i)}_{\bf i,j;\wtil\phi}
=\xi^{(i)}_{\bf i,l;\wtil\phi'}$.
\qed

\subsection{Generalized Minors and the potentials}

Let $G$ be a simply connected simple algebraic groups over $\bbC$ and 
$T\subset G$ a maximal torus. 
Let  $X^*(T):=\Hom(T,\bbC^\times)$ and $X_*(T):=\Hom(\bbC^\times,T)$ be
the lattice of characters and co-characters respectively.
We identify $P$ (resp. $P^*$) with $X^*(T)$ 
(resp. $X_*(T)$) as in Sect.2.

\begin{df}[\cite{BZ,BZ2}]
For a dominant weight $\mu\in P_+$, the
{\it principal minor} $\Del_\mu:G\to\bbC$ is defined as
\[
 \Del_\mu(u^-tu^+):=\mu(t)\q(u^\pm\in U^\pm,\,\,t\in T).
\]
Let $\gamma,\del\in P$ be extremal weights such that 
$\gamma=u\mu$ and $\del=v\mu$ for some $u,v\in W$. 
Then the {\it generalized minor} $\Del_{\gamma,\del}$ is defined
by
\[
 \Del_{\gamma,\del}(g):=\Del_\mu(\ovl u^{-1}g\ovl v)
\q(g\in G),
\]
which is a regular function on $G$.
\end{df}
Let us denote the one-dimensional additive group $\bbG_a(=\bbC)$ and 
multiplicative group $\bbG_m(=\bbC^\times)$.
Let $\what U:={\rm Hom}(U,\bbG_a)$ be the set of additive characters 
of $U$.
The {\it elementary character }$\chi_i\in \what U$ and  
the {\it standard regular character} $\chi^{\rm st}\in \what U$ are  
defined as follows (\cite{BK2}):
\[
\chi_i(x_j(c))=\del_{i,j}\cdot c \q(c\in \bbC,\,\, i,j\in I),\qq
\chi^{st}=\sum_{i\in I}\chi_i.
\]
\begin{lem}[\cite{BK2}]
Suppose that $G$ is a simply connected simple algebraic group as above.
\begin{enumerate}
\item
For $u\in U$ and $i\in I$, we have 
$\Del_{\mu,\mu}(u)=1$ and $\chi_i(u)=\Del_{\Lm_i,s_i\Lm_i}(u)$,
where $\Lm_i$ is the $i$-th fundamental weight.
\item
Define the map $\pi^+:B^-\cdot U\to U$ by $\pi^+(bu)=u$ for 
$b\in B^-$ and $u\in U$. For any $g\in G$, we have 
\begin{equation}
 \chi_i(\pi^+(g))=\frac{\Del_{\Lm_i,s_i\Lm_i}(g)}
{\Del_{\Lm_i,\Lm_i}(g)}.
\end{equation}
\end{enumerate}
\end{lem}

\begin{df}[\cite{BK2}]
\begin{enumerate}
\item
Define the rational functions $\Phi^{(\pm)}$ on $G$ by 
\begin{equation}
\Phi^{(+)}(g)=\chi^{st}(\pi^+({\ovl w_0}^{-1}g))=\sum_i\frac{\Del_{w_0\Lm_i,s_i\Lm_i}(g)}{\Del_{w_0\Lm_i,\Lm_i}(g)}
,\qq
\Phi^{(-)}(g)
=\chi^{st}(\pi^+({\ovl w_0}^{-1}\eta(g)))
=\sum_i\frac{\Del_{w_0s_i\Lm_i,\Lm_i}(g)}
{\Del_{w_0\Lm_i,\Lm_i}(g)},
\label{chi-min}
\end{equation}
where $\eta:G\to G$ is the anti-automorphism defined by
\[
\eta(x_i(c))=x_i(c),\q \eta(y_i(c))=y_i(c),\q \eta(t)=t^{-1}.
\]
\item
For $w\in W$ and a character $\chi$, 
we define the function $\Phi_{w,\chi}$ on $G$ by 
\begin{equation}
\Phi_{w,\chi}(g)=\chi(\pi^+({\ovl w}^{-1}g)).
\label{phi-w-chi}
\end{equation}
\end{enumerate}
\end{df}
Note that $\Phi^{(+)}=\Phi_{w_0,\chi^{st}}$.

Let $\omega:\ge\to\ge$ be the anti-involution 
\[
\omega(e_i)=f_i,\q
\omega(f_i)=e_i\,\q\omega(h)=h,
\] and extend it to $G$ by setting
$\omega(x_i(c))=y_i(c)$, $\omega(y_i(c))=x_i(c)$ and $\omega(t)=t$
$(t\in T)$.

There exists a $\ge$(or $G$)-invariant bilinear form on the
finite-dimensional  irreducible
$\ge$(or $G$)-module $V(\lm)$ such that 
\begin{equation}
 \lan au,v\ran=\lan u,\omega(a)v\ran,
\q\q(u,v\in V(\lm),\,\, a\in \ge(\text{or }G)).
\label{inv-bilin}
\end{equation}
For $g\in G$, 
we have the following simple fact:
\[
 \Del_{\Lm_i}(g)=\lan gu_{\Lm_i},u_{\Lm_i}\ran,
\]
where $u_{\Lm_i}$ is a highest weight vector in
the fundamental representation $V(\Lm_i)$ satisfying 
$\lan u_{\Lm_i},u_{\Lm_i}\ran=1$. Hence, for $w,w'\in W$ we have
\begin{equation}
 \Del_{w\Lm_i,w'\Lm_i}(g)=
\Del_{\Lm_i}({\ovl w}^{-1}g\ovl w')=
\lan {\ovl w}^{-1}g\ovl w'\cdot u_{\Lm_i},u_{\Lm_i}\ran
=\lan g\ovl w'\cdot u_{\Lm_i}\, ,\, \ovl{w}\cdot u_{\Lm_i}\ran,
\label{minor-bilin}
\end{equation}
where note that $\omega(\ovl s_i)=\ovl s_i$.

\begin{pro}\label{form-alt}
Let ${\bf i}=i_1\cd i_N$ be a reduced word for the longest Weyl
group element $w_0$. 
For $\Theta_{\bf i}^-(c)\in B^-_{w_0}$, 
we get the following formula.
\begin{equation}
\Phi^{(+)}(\Theta_{\bf i}^-(c))
=\sum_i\Del_{w_0\Lm_i,s_i\Lm_i}(\Theta_{\bf i}^-(c)), \q
\Phi^{( - )}(\Theta_{\bf i}^-(c))=
\sum_i\Del_{w_0s_i\Lm_i,\Lm_i}(\Theta_{\bf i}^-(c)).
\label{fb-th}
\end{equation}
\end{pro}
{\sl Proof.}
We must show that 
\begin{equation}
\Del_{w_0\Lm_i,\Lm_i}(\Theta_{\bf i}^{-}(c))=1.
\label{del-1}
\end{equation}
Since $\Theta_{\bf i}^{-}(c)\in U\ovl w_0 U$, we have 
$\ovl w_0^{-1}\Theta_{\bf i}^{-}(c)\in \ovl w_0^{-1}U\ovl w_0U=U^-\cdot
U$.
So, there exist  $u_1\in U^-$ and $u_2\in U$ such that 
$\ovl w_0^{-1}\Theta_{\bf i}^{-}(c)=u_1u_2$.
Thus, it follows from \eqref{minor-bilin} and the fact $\omega(u_1)\in U$
that 
\[
\Del_{w_0\Lm_i,\Lm_i}(\Theta_{\bf i}^{-}(c))=
\lan \ovl w_0^{-1}\Theta_{\bf i}^-(c)u_{\Lm_i},u_{\Lm_i}\ran
=\lan u_1u_2u_{\Lm_i},u_{\Lm_i}\ran
=\lan u_2u_{\Lm_i},\omega(u_1)u_{\Lm_i}\ran=\lan
u_{\Lm_i},u_{\Lm_i}\ran=1.
\]
\qq \qed

Set $\Phi_{\rm BK}:=\Phi^{(+)}+\Phi^{( - )}$ and then it is shown in 
\cite{BK2} that $\Phi_{\rm BK}$ is a potential on $TB^-_{w_0}$, 
which deduces the following result:

\begin{pro}\label{ul-dec-gc}
The function $\Phi^{(\pm)}$ is an upper/lower half-potential on the 
geometric crystal $B^-_{w_0}$ and then $(B^-_{w_0},\Phi^{(\pm)})$
is an upper/lower half-decorated geometric crystal.
\end{pro}
{\sl Proof. }
As we have seen that $B^-_{w_0}$ holds the positive 
geometric crystal structure, it suffices to show that 
the function $\Phi^{(\pm)}$ is an upper/lower half-potential on $B^-_{w_0}$. 
We only see $\Phi^{(+)}$ since it will be done similarly for $\Phi^{( - )}$.
The action $e^c_i$ on an element in $B^-_{w_0}$ is given in \eqref{ei-action}:
\[
e^c_i(x)=x_i((c-1)\vp_i(x))\cdot x\cdot x_i((c^{-1}-1)\vep_i(x)).
\]
Then applying \eqref{minor-bilin}, for $g\in B^-_{w_0}$ we obtain
\begin{eqnarray*}
\Del_{w_0\Lm_j,s_j\Lm_j}(e_i^c(g))&=&
\lan x_i((c-1)\vp_i(g))\cdot g\cdot x_i((c^{-1}-1)\vep_i(g))s_ju_{\Lm_j},
\ovl w_0\cdot u_{\Lm_j}\ran\\
&=&
\lan g\cdot x_i((c^{-1}-1)\vep_i(g))s_ju_{\Lm_j},
y_i((c-1)\vp_i(g))\cdot \ovl w_0\cdot u_{\Lm_j}\ran\\
&=&
\lan g\cdot (1+(c^{-1}-1)\vep_i(g)e_i)f_ju_{\Lm_j},
\ovl w_0\cdot u_{\Lm_j}\ran\\
&=&
\lan g\cdot f_ju_{\Lm_j},
\ovl w_0\cdot u_{\Lm_j}\ran+
\lan g\cdot (c^{-1}-1)\vep_i(g)(f_je_i+\del_{ij}h_i)\cdot u_{\Lm_j},
\ovl w_0\cdot u_{\Lm_j}\ran\\
&=&
\lan g\cdot \ovl s_j\cdot u_{\Lm_j},
\ovl w_0\cdot u_{\Lm_j}\ran+
\del_{ij}(c^{-1}-1)\vep_i(g)\cdot 
\lan g\cdot u_{\Lm_j},
\ovl w_0\cdot u_{\Lm_j}\ran\\
&=&\Del_{w_0\Lm_j,s_j\Lm_j}(g)+\del_{ij}(c^{-1}-1)\vep_i(g),
\end{eqnarray*}
where for the 3rd equality, we use the fact that $\ovl w_0\cdot u_{\Lm_j}$
is the lowest weight vector and $e_i^2f_ju_{\Lm_j}=0$, 
and for the last equality, we use 
$\lan g\cdot u_{\Lm_j},\ovl w_0\cdot u_{\Lm_j}\ran=
\Del_{w_0\Lm_j,\Lm_j}(g)=1$ from the proof of Proposition \ref{form-alt}.
Therefore, we obtain 
\[
\Phi^{(+)}(e_i^c(g))=\sum_j \Del_{w_0\Lm_j,s_j\Lm_j}(e_i^c(g))
=\Phi^{(+)}(g)+(c^{-1}-1)\vep_i(g),
\]
which completes the proof. \qed
\subsection{Tropicalization of the half decorated geometric crystal 
$B^-_{w_0}$ and the crystal $B(\ify)$ }

As we have seen in the last subsection, the geometric crystal $B^-_{w_0}$ is 
equipped with the upper half potential $\Phi^{(+)}$ and then we consider 
the tropicalized crystal 
\[
(\bbB^-_{w_0})_{\Phi^{(+)},\Theta_\bfi}
=((\wtil B^-_{w_0})_{\Phi^{(+)},\Theta_{\bfi}}, 
\tilde {{e_{i}}}_{|(\wtil B^-_{w_0})_{\Phi^{(+)},\Theta_{\bfi}}},
{\wt_i}_{|(\wtil B^-_{w_0})_{\Phi^{(+)},\Theta_{\bfi}}}, 
{\vep_i}_{|(\wtil B^-_{w_0})_{\Phi^{(+)},\Theta_{\bfi}}}),
\]
 where 
$\eit, \wt_i,  \vep_i$ are given in \eqref{eitn}--\eqref{vep-ud} and 
\begin{equation}
(\wtil B^-_{w_0})_{\Phi^{(+)},\Theta_{\bfi}}
:=\{x=(x_1,\cd,x_N)\in{\mathcal TR}(B_{w_0}^-)=
\bbZ^N\,|\,\wtil \Phi^{(+)}(x)\geq0\}. 
\label{tr-bw}
\end{equation}
By Proposition \ref{ul-normal}, one knows that the crystal 
$(\bbB^-_{w_0})_{\Phi^{(+)},\Theta_\bfii0}$ 
is a $\ge^\vee$-upper normal crystal.
Now, let us show the following theorem:
\begin{thm}\label{b-ify}
As a $\ge^\vee$-crystal, the crystal 
$(\bbB^-_{w_0})_{\Phi^{(+)},\Theta_\bfi}$
is isomorphic to the crystal $B(\ify)$.
\end{thm}
{\sl Proof.} 
To show this theorem, we shall use 
the Kashiwara-Saito characterization of the crystal $B(\ify)$
(Corollary \ref{ks-cor}).
For the purpose, let us define the map $\TY(\Psi,+,i)$
for each $i\in I$,
\begin{equation}
\Psi^{(+)}_i:(\wtil B^-_{w_0})_{\Phi^{(+)},\Theta_{\bfi}}\hookrightarrow 
(\wtil B^-_{w_0})_{\Phi^{(+)},\Theta_{\bfi}}\ot B_i,
\label{embed-}
\end{equation}
by 
\begin{equation}
\Psi^{(+)}_i(x_1,\ldots,x_N):=
\xi^{(i)}_{\bfi}(x_1,\ldots,x_N)\ot
\fit^{\omega_i(x_1,\ldots,x_N)}(0)_i,
\label{embed-def}
\end{equation}
where the map $\omega_i$ and $\xi^{(i)}_{\bfi}$ are given in 
Proposition \ref{optimal}.

\begin{lem}\label{1N=2Nj}
For a reduced longest word $\bfi=i_1i_2\cd i_N$, define $\ist\in I$
to be a unique index satisfying $-\al_\ist=w_0(\al_i)$ $(i=i_1)$.
Then we find that 
\begin{equation}
w_0=s_{i_i}s_{i_2}\cd s_{i_N}=s_{i_2}s_{i_3}\cd s_{i_N}s_\ist.\label{1Nj}
\end{equation}
\end{lem}
{\sl Proof.}
Applying  \cite[Theorem 1.7]{Hum} to $w=s_{i_1}s_{i_2}\cd s_{i_N}s_\ist$ 
we find that the case $i=i_1$ and $j=\ist$ are the unique 
case matching to the claim of the theorem. \qed

\begin{lem}\label{gabe-lem}
Let $\bfi=i_1\cd i_N$ be a reduced word for $w_0\in W$,  
$i=i_1$ and $\ist$ an index satisfying $\Lm_{\ist}=-w_0(\Lm_i)$.
Then, there exists some positive integer $M$ and a Laurent 
polynomial $P_i$ ($i\in I$) in $(c_2,\cd,c_N)$ with positive 
integer coefficients such that 
\begin{eqnarray}
&&\Del_{w_0\Lm_{\ist},s_{\ist}\Lm_{\ist}}(\Theta_{\bfi}^-(c_1,\cd,c_N))
=Mc_1+P_\ist(c_2,\cd,c_N),
\label{gabe-eq}\\
&&\Del_{w_0\Lm_{j},s_{j}\Lm_{j}}(\Theta_{\bfi}^-(c_1,\cd,c_N))
=P_j(c_2,\cd,c_N)\,\,{\rm for}\,j\ne \ist,
\label{gabe-eq2}
\end{eqnarray}
\end{lem}
{\sl Proof.}
For the proof, we shall use Berenstein-Zelevinsky's $\bfi$-trail, 
which will be introduced in Sect.7. In particular Theorem \ref{trail-thm}
claims that 
\begin{equation}
\Del_{w_0\Lm_{\ist},s_{\ist}\Lm_{\ist}}(\Theta_{\bfi}^-(c_1,\cd,c_N))
=\sum_{\pi\in \Pi_\ist}M_\pi c_1^{d_1(\pi)}\cd c_N^{d_n(\pi)}, 
\label{gabe-m}
\end{equation}
where $\Pi_\ist$ is the set of all $\bfi$-trails from 
$\Lm_i$ to $-s_{i^*}(\Lm_{i^*})$ and $M_\pi$ is a positive integer.
For $\pi=\gamma_0,\gamma_1,\cd,\gamma_N)\in \Pi_\ist$, 
it follows from the definition 
$d_r(\pi)=\frac{\gamma_{r-1}+\gamma_r}{2}(h_{i_r})$
that in our setting
$d_r(\pi)=0$ if and only if $\gamma_r=s_{i_r}(\gamma_{r-1})$.
We find that there are two possibilities for $\gamma_1$:
\begin{enumerate}
\item $\gamma_1=\gamma_0-\al_i$.
\item $\gamma_1=\gamma_0$.
\end{enumerate}
In the first case (i), one knows that $d_1(\pi)=0$ and then 
the corresponding term in \eqref{gabe-m} is a Laurent monomial 
in $c_2,\cd,c_N$.

In the second case (ii), namely,  $\gamma_0=\gamma_1=\Lm_i$, 
let $\bfi'=(i_2,\cd,i_N,i^*)$.
By Lemma \ref{1N=2Nj}, this $\bfi'$ is another reduced word for $w_0$.
Then, the sequence $(\gamma_1,\gamma,_2,\cd,\gamma_N,-\Lm_{i^*})$
is an $\bfi'$-trail in $V(\Lm_i)$ from $\Lm_i$ to $-\Lm_{i^*}$,
it contributes a positive multiple of $t_2^{d_2(\pi)}\cd t_N^{d_N(\pi)}$
to $\Del_{w_0\Lm_{\ist},s_{\ist}\Lm_{\ist}}(\Theta_{\bfi'}^-(t_2,\cd,t_N,t_1))$.
But, this generalized minor must be 1 by \eqref{del-1}, which means
$d_2(\pi)=d_3(\pi)=\cd=d_N(\pi)=0$.
Hence, we conclude that there is unique $\bfi$-trail $\pi^*$ satisfying 
$\gamma_0=\gamma_1$, and this trail $\pi^*$ contributes
$M_{\pi^*}c_1^{d_1(\pi^*)}=M_{\pi^*}\cdot c_1$ to 
$\Del_{w_0\Lm_{\ist},s_{\ist}\Lm_{\ist}}(\Theta_{\bfi}^-(c_1,\cd,c_N))$.

Next, we take $j\in I$ such that $j\ne \ist$. It holds
\begin{equation*}
\Del_{w_0\Lm_{j},s_{j}\Lm_{j}}(\Theta_{\bfi}^-(c_1,\cd,c_N))
=\sum_{\pi\in \Pi_j}M_\pi c_1^{d_1(\pi)}\cd c_N^{d_n(\pi)},
\end{equation*}
where $\Pi_j$ is the set of all $\bfi$-trails from 
$-w_0(\Lm_j)$ to $-s_{j}(\Lm_{j})$ and $M_\pi$ is a positive integer.
For any $\pi=(\gamma_0,\cdots,\gamma_N)\in \Pi_j$,
it follows by $-w_0(\Lambda_j)= \Lm_{j^*}$ with some $j^*\in I\setminus\{i_1\}$ that
\[
\gamma_0=\gamma_1=-w_0(\Lambda_j),
\]
which yields $d_1(\pi)=\frac{\gamma_0+\gamma_1}{2}(h_{i_1})=-w_0(\Lm_j)(h_i)=0$. Hence, we get (\ref{gabe-eq2}).\qed


Here we set
$P(c_2,\cd,c_N):=\sum_{i\in I}P_i(c_2,\cd,c_N)$. 
Note that the definition of the map 
$\Psi^{(+)}_i$ is well-defined by Proposition \ref{optimal} 
and Lemma \ref{gabe-lem}. Indeed, to see the well-definedness, we may show
that if $\wtil\Phi^{(+)}_{\bfi}(x_1,x_2,\cd,x_N)\geq 0$, then 
$\wtil\Phi^{(+)}_{\bfi}(\xi^{(i)}_\bfi(x_1,x_2,\cd,x_N))=
\wtil\Phi^{(+)}_{\bfi}(0,x_2,\cd,x_N)\geq 0$. 
By Lemma \ref{gabe-lem}, it is clear that 
\begin{equation}
\wtil\Phi^{(+)}_{\bfi}(x_1,x_2,\cd,x_N)
=\min(x_1,{\mathcal TR}(P)(x_2,\cd,x_N))\geq0,
\label{gabe-star}
\end{equation}
which implies that 
$\wtil\Phi^{(+)}_{\bfi}(0,x_2,\cd,x_N)
={\mathcal TR}(P)(x_2,\cd,x_N)\geq0$.

Let us show that this $\TY(\Psi,+,i)$, 
$(\wtil B^-_{w_0})_{\Phi^{(+)},\Theta_{\bfi}}$ and 
$b_0=(0,0,\cd,0)$ satisfy the conditions 
(i)--(iv) in Theorem \ref{KS-char} and (vii') in 
Corollary \ref{ks-cor}.
Let us denote simply 
$(\wtil B^-_{w_0})_{\Phi^{(+)},\Theta_{\bfi}}$ by $\bbB_{\bfi}$.

Indeed, (iii) and (iv) are trivial by \eqref{vep-ud}. 
Let us see the condition (vi) in Theorem \ref{KS-char}:
To do this, it suffices to show that 
\begin{equation}
\text{For any }(x_1,\cd,x_N)\in \bbB_{\bfi},\,\,\text{we have }
x_1\geq0,
\label{x10>0}
\end{equation}
which is trivial from \eqref{gabe-star}.

As for the condition (i), we refer to the argument in the proof of 
\cite[Lemma 6.23]{BK2} to get 
\begin{equation}
\lan \Lm^\vee_i,\wt(x)\ran\leq0,
\label{wt<0}
\end{equation}
where $\Lm_i^\vee$ is 
the element in $\LL(P^*)$ such that $\lan \Lm^\vee_i,\al_j\ran=
\del_{ij}$, which implies the condition (i).

Let us check the condition (vii') in Corollary \ref{ks-cor}. 
As we have seen that the crystal 
$(\wtil B^-_{w_0})_{\Phi^{(+)},\Theta_{\bfi}}$ is an upper normal
crystal, for any $x\in (\wtil B^-_{w_0})_{\Phi^{(+)},\Theta_{\bfi}}$
there exists some $i$ such that $\eit x\ne0$ iff
$\vep_i(x)>0$. 
Now, we shall introduce the following lemma:
\begin{lem}[\cite{BK2}Lemma 6.19]\label{uniq}
Let $\vep_i(x)$ be the function as in \eqref{vep-ud} and 
$\wtil\Phi_{w_0,\chi^{st}}$ the tropicalization of 
the function $\Phi_{w_0,\chi^{st}}$ as in \eqref{phi-w-chi} for 
the longest element $w_0\in W$ and the standard character $\chi^{st}$.
If an element
$x=(x_1,\cd,x_N)\in {\mathcal TR}(B^-_{w_0})$
satisfies
\begin{equation}
\wtil\Phi_{w_0,\chi^{st}}(x)\geq0,\qq \vep_i(x)\leq 0\q
{\rm for \,any } \,\,i\in I,\label{vep=0}
\end{equation}
then $x=b_0=(0,0,\cd,0)$. 
\end{lem}
By this lemma we know that
if $x\ne(0,0,\cd,0)$ then $\vep_i(x)>0$ for some $i\in I$ and then 
$\eit x\ne0$, which means the condition (vii').

Let us show that the condition (ii) holds.
By Lemma \ref{uniq} and the 
fact $\Phi^{(+)}=\Phi_{w_0,\chi^{st}}$, we find that 
$b_0=(0,0,\cd,0)$ 
is the unique element in $(\wtil B^-_{w_0})_{\Phi^{(+)},\Theta_{\bfi}}$
such that $\wt_i(x)=0$ and $\vep_i(x)=0$ for any $i$.
Here let us assume that there exists 
$x_0\in (\wtil B^-_{w_0})_{\Phi^{(+)},\Theta_{\bfi}}$ such that 
$x_0\ne b_0$ and $\wt(x_0)=0$. 
By the condition (vii'), 
there exists some $j$ such that $\vep_j(x_0)\ne0$ and 
then $\til e_j(x_0)\in
(\wtil B^-_{w_0})_{\Phi^{(+)},\Theta_{\bfi}}$. However, 
one has $\wt(\til e_j(x_0))=\al_j\not\in Q_-$, which 
contradicts (i). Thus, the condition (ii) holds.

Finally, let us show the condition (v). 
Since the braid-type isomorphisms $\phi^{(k)}_{ij}$ 
are isomorphisms of 
crystals, it suffices to show the following (a)--(d)
under the assumption $\bfi=ii_2\cd i_N$, 
\begin{enumerate}
\item[(a)] $\Psi_i^{(+)}$ is injective.
\item[(b)] For any 
$x=(x_1,\cd,x_N)\in (\wtil B^-_{w_0})_{\Phi^{(+)},\Theta_{\bfi}}$, 
one has $\wt((0,x_2,\cd,x_N)\ot (-x_1)_i)=\wt(x).$
\item[(c)] For any $j\in I$ and 
$x=(x_1,\cd,x_N)\in (\wtil B^-_{w_0})_{\Phi^{(+)},\Theta_{\bfi}}$, 
one has $\vep_j((0,x_2,\cd,x_N)\ot (-x_1)_i)=\vep_j(x).$
\item[(d)] For any $j\in I$ and 
$x=(x_1,\cd,x_N)\in (\wtil B^-_{w_0})_{\Phi^{(+)},\Theta_{\bfi}}$, 
one has
$\til e_j(\Psi^{(+)}_i(x))=\Psi^{(+)}_i(\til e_j(x))$ 
and $\til f_j(\Psi^{(+)}_i(x))=\Psi^{(+)}_i(\til f_j(x))$.
\end{enumerate}
Since $\Psi^{(+)}_i(x_1,\cd,x_N)=(0,x_2,\cd,x_N)\ot (-x_1)_i$,
 it is evident that $\Psi^{(+)}_i$ is injective, which shows (a).

We have 
\[
\wt(x_1,\cd,x_N)=-\sum_{k=1}^Nx_k\al_{i_k}=
\wt((0,x_2,\cd,x_N)\ot \til f_{i_1}^{x_1}(0)_{i_1})
\]
which shows (b).

Let us see (c). First suppose that $i=j$. 
Due to \eqref{vep-ud} we have 
\[
\vep_i(x_1,\cd,x_N)=\max_{1\leq m\leq N,i_m=i}\
(x_m+\sum_{l=m+1}^N {\bf a}_{i_l,i}x_l).
\]
On the other-hand, by \eqref{tensor-vep} we find that 
\begin{eqnarray*}
&&\vep_i((0,x_2,\cd, x_N)\ot (-x_1)_i)\\
&&\qq\qq=
\max(\vep_i(0,x_2,\cd, x_N), 
\vep_i((-x_1)_i)-\lan \Lm^\vee_i,\wt((0,x_2,\cd,x_N))
\ran)\\
&&\qq\qq=\max\left(\max\left(\sum_{l=2}^N {\bf a}_{i_l,i}x_l,\max_{2\leq m\leq N,i_m=i}
\left(x_m+\sum_{l=m+1}^N {\bf a}_{i_l,i}x_l\right)\right),
x_1+\sum_{l=2}^N {\bf a}_{i_l,i}x_l\right)\\
&&\qq\qq=
\max_{1\leq m\leq N,i_m=i}\
(x_m+\sum_{l=m+1}^N {\bf a}_{i_l,i}x_l)=\vep_i(x_1,\cd,x_N).
\end{eqnarray*}
Here note that the 3rd equality is derived by 
the fact $x_1\geq0$.

Next, let us show the case $i\ne j$.
Similar to the previous case, we get 
\begin{eqnarray*}
\vep_j((0,x_2,\cd, x_N)\ot (-x_1)_i)
&=&\max(\vep_j(0,x_2,\cd, x_N), 
\vep_j((-x_1)_i)-\lan \Lm^\vee_j,\wt((0,x_2,\cd,x_N))
\ran)\\
&=&\max\left(
\max_{1\leq m\leq k,i_m=j}
\left(x_m+\sum_{l=m+1}^k{\bf a}_{i_l,j}x_l\right),
-\ify\right)\\
&=&
\max_{1\leq m\leq k,i_m=j}\
(x_m+\sum_{l=m+1}^k{\bf a}_{i_l,j}x_l)=\vep_j(x_1,\cd,x_N),
\end{eqnarray*}
where the 2nd and 3rd equalities use the fact that $i\ne j$. Now we get (c).

Finally, let us show (d).
By \eqref{eitn} we obtain the following explicit actions
of $\eit$ and $\fit$ on 
$(\wtil B^-_{w_0})_{\Phi^{(+)},\Theta_{\bfi}}$
(see \cite[Lemma 3.10]{N4}):
\begin{lem}[\cite{N4}]
\label{lem-ef}
For $i\in I$ and $x=(x_1,\cd,x_N)\in 
(\wtil B^-_{w_0})_{\Phi^{(+)},\Theta_{\bfi}}$, 
set $X_m(x):=x_m+\sum_{k=1 }^{m-1}
{\bf a}_{i_k,i_m}x_k$, 
${\mathcal X}^{(i)}(x):=\min\{X_m(x)|1\leq
m\leq N,i_m=i\}$ 
and $M^{(i)}(x):=\{l|1\leq l\leq N,i_l=i,
 X_l(x)={\mathcal X}^{(i)}(x)\}$.
Define $\TY(m,i,e)(x):=\max(M^{(i)}(x))$ and 
$\TY(m,i,f)(x):=\min(M^{(i)}(x))$:
Then we have
\begin{eqnarray}
&&\eit(x)=\begin{cases}(x_1,\cd,x_{\TY(m,i,e)}-1,\cd,x_N)
&\text{ if }\wtil \Phi^{(+)}(x_1,\cd,x_{\TY(m,i,e)}-1,\cd,x_N)\geq0,\\
0&\text{otherwise,}
\end{cases}
\label{th-eaction}\\
&&\fit(x)=\begin{cases}(x_1,\cd,x_{\TY(m,i,f)}+1,\cd,x_N)
&\text{ if }\wtil \Phi^{(+)}(x_1,\cd,x_{\TY(m,i,f)}+1,\cd,x_N)\geq0,\\
0&\text{otherwise.}
\end{cases}
\label{th-faction}
\end{eqnarray}
\end{lem}
\begin{rem}
Note that the definition of $X_m$ in \cite[Lemma 3.10]{N4} is incorrect.
The one here is correct. 
\end{rem}

Consider the case  $i=i_1\ne j$. 
Since $i\ne j$, if $i_m=j$, 
we have $X_m(x)=X_m(0,x_2,\cd,x_N)+{\bf a}_{ij}x_1$ and then 
for $m,l$ with $i_m=i_l=j\ne i$ we obtain 
$l,m>1$ and 
$X_m(x)-X_l(x)=X_m(0,x_2,\cd,x_N)-X_l(0,x_2,\cd,x_N)$. Therefore, we obtain 
both $\TY(m,j,e)(x)=\TY(m,j,e)(0,x_2,\cd,x_N)=:\TY(m,j,e)$ and 
$\TY(m,j,f)(x)=\TY(m,j,f)(0,x_2,\cd,x_N)=:\TY(m,j,f)$.
Assume $\til e_j(x_1,\cd,x_N)=
(x_1,\cd,x_{\TY(m,j,e)}-1,\cd,x_N)\in 
(\wtil B^-_{w_0})_{\Phi^{(+)},\Theta_{\bfi}}$.
And one has also $\vep_j(\fit^{x_1}(0)_i)=-\ify$, which implies
\begin{eqnarray*}
\til e_j((0,x_2,\cd,x_N)\ot \fit^{x_1}(0)_i)
=(0,x_2,\cd,x_{\TY(m,j,e)}-1,\cd,x_N)\ot \fit^{x_1}(0)_i
=\Psi_i^{(+)}(\til e_j(x)).
\end{eqnarray*}
By a similar way, we can show $\til f_j(\Psi^{(+)}_i)(x)=
\Psi^{(+)}_i(\til f_j(x))$.
Next, assume that $\til e_j(x)=0$, namely, 
$\til e_j(x)\not\in (\wtil B^-_{w_0})_{\Phi^{(+)},\Theta_{\bfi}}$. 
It suffices to show that 
$\til e_j(0,x_2,\cd,x_N)\not\in 
(\wtil B^-_{w_0})_{\Phi^{(+)},\Theta_{\bfi}}$.
Since 
$\wtil \Phi^{(+)}(x_1,\cd,x_N)$ is in the form 
\[
\min(x_1,\text{linear functions in }x_2,\cd,x_N),
\]
by \eqref{gabe-star} and the fact that 
$\TY(m,i,e)>1$, if $(x_1,\cd,x_{\TY(m,i,e)}-1,\cd,x_N)\not\in
(\wtil B^-_{w_0})_{\Phi^{(+)},\Theta_{\bfi}}$, then 
one has 
$(0,\cd,x_{\TY(m,i,e)}-1,\cd,x_N)\not\in
(\wtil B^-_{w_0})_{\Phi^{(+)},\Theta_{\bfi}}$.
The case of $\til f_j$ is also done similarly.

Finally, let us assume $i=j$ and show $\eit(\Psi^{(+)}_i(x))
=\Psi^{(+)}_i(\eit(x))$.
Set $x'=(0,x_2,\cd,x_N)$. One has
\begin{eqnarray*}
\vp_i(x')&=&
\wt_i(x')+\vep_i(x')
=-\sum_{j=2}^N{\bf a}_{i_j,i}x_j+
\max\left(\sum_{j=2}^N{\bf a}_{i_j,i}x_j,
\max_{2\leq m\leq N,\,\,i_m=i}
\left(x_m+\sum_{j=m+1}^N{\bf a}_{i_j,i}x_j\right)\right)\\
&=&\max\left(0,\max_{2\leq m\leq N,\,i_m=i}
\left(-x_m-\sum_{j=2}^{m-1}{\bf a}_{i_j,i}x_j\right)\right)\\
&=&-\min(0,\min_{2\leq m\leq N,i_m=i}(X_m(x'))).
\end{eqnarray*}
Note that by this explicit form of $\vp_i(x')$, we find that 
$\vp_i(x')\geq0$.
Let us consider the case 
$\eit(x)\in (\wtil B^-_{w_0})_{\Phi^{(+)},\Theta_{\bfi}}$:
\begin{enumerate}
\item
Assume 
\begin{equation}
\vp_i(x')\geq x_1=\vep_i(\fit^{x_1}(0)_i).\label{vp>vep}
\end{equation}
We have $\eit(x'\ot \fit^{x_1}(0)_i)=\eit(x')\ot \fit^{x_1}(0)_i$.
Thus, it is sufficient to show $\TY(m,i,e)(x)=\TY(m,i,e)(x')$.
Indeed, since $X_j(x)=2x_1+X_j(x')$ for $j>1$, if $\TY(m,i,e)(x)>1$ then 
we obtain $\TY(m,i,e)(x)=\TY(m,i,e)(x')$. If $\TY(m,i,e)(x)=1$, then 
$X_1(x)=x_1<X_j(x)=2x_1+X_j(x')$ for all $j>1$, 
which means $-X_j(x')<x_1$ and then it contradicts \eqref{vp>vep}
since $\vp_i(x')=-X_j(x')$ for some $j$.
Therefore, we get $\TY(m,i,e)(x)>1$ and then $\TY(m,i,e)(x)=\TY(m,i,e)(x')$, which 
implies
$\eit(\Psi^{(+)}_i(x))=\Psi^{(+)}_i(\eit (x))$.
\item
Assume 
\begin{equation}
\vp_i(x')<x_1=\vep_i(\fit^{x_1}(0)_i).\label{vp<vep}
\end{equation}
In this case, we have
\[
\eit(x'\ot \fit^{x_1}(0)_i)=x'\ot \fit^{x_1-1}(0)_i
\]
Here note that by the fact $\vp_i(x')\geq0$, we have $x_1>0$.
We shall show $\TY(m,i,e)(x)=1$. Regarding \eqref{vp<vep}, since
$-x_1<-\vp_i(x')\Leftrightarrow x_1<2x_1-\vp_i(x')$, 
one has
\begin{eqnarray*}
x_1&<&2x_1+\min(0,\min_{2\leq m\leq N,i_m=i}(X_m(x')))\\
&=&\min(2x_1,\min_{2\leq m\leq N,i_m=i}(X_m(x)))),
\end{eqnarray*}
which show that $x_1=X_1(x)<X_j(x)$ for any $j>1$ and then 
$\TY(m,i,e)(x)=1$. Hence, we get 
$\eit(x)=(x_1-1,x_2,\cd,x_N)$ and then 
$\eit(\Psi^{(+)}_i(x))=\Psi^{(+)}_i(\eit (x))$.
\end{enumerate}
Finally, let us consider the case $\eit(x)=0$, namely, 
$\eit(x)\not\in (\wtil B^-_{w_0})_{\Phi^{(+)},\Theta_{\bfi}}
\Leftrightarrow$
\begin{equation}
\wtil \Phi^{(+)}(\eit(x_1,\cd,x_N))<0.
\label{Phi-nega}
\end{equation}
Here we also consider the following two cases, 
\begin{enumerate}
\item
$\vp_i(x')\geq x_1=\vep_i(\fit^{x_1}(0)_i)$.
\item
$\vp_i(x')<x_1=\vep_i(\fit^{x_1}(0)_i)$.
\end{enumerate}
Indeed, $\TY(m,i,e)$ does not depend on $\wtil\Phi^{(+)}$ and then
by the above argument 
\begin{enumerate}
\item
$\TY(m,i,e)(x)=\TY(m,i,e)(x')>1$ .
\item
$\TY(m,i,e)(x)=\TY(m,i,e)(x')=1$. 
\end{enumerate}
Set $(x_1'',x_2'',\cd,x_N''):=\eit(x_1,x_2,\cd,x_N)$.
In the case (i), we obtain 
$\eit((0,x_2,\cd,x_N)\ot \fit^{x_1}(0)_i)=
(0,x_2'',\cd,x_N'')\ot \fit^{x_1}(0)_i$.
Since by $x_1\geq0$ we have
\[
0>\wtil\Phi^{(+)}(\eit x)=
\wtil\Phi^{(+)}(x_1,x_2'',\cd,x_N'')
=\min(x_1,{\mathcal TR}(P)(x_2'',\cd,x_N''))
={\mathcal TR}(P)(x_2'',\cd,x_N''), 
\]
we obtain $\eit(0,x_2,\cd,x_N)\not\in 
(\wtil B^-_{w_0})_{\Phi^{(+)},\Theta_{\bfi}}$. 
In the case (ii), we have
$\eit((0,x_2,\cd,x_N)\ot \fit^{x_1}(0)_i)=
(0,x_2,\cd,x_N)\ot \fit^{x_1-1}(0)_i$.
Since by $x_1>0$ we get 
\[
0>\wtil\Phi^{(+)}(\eit x)=
\wtil\Phi^{(+)}(x_1-1,x_2,\cd,x_N)
=\min(x_1-1,{\mathcal TR}(P)(x_2,\cd,x_N))
={\mathcal TR}(P)(x_2,\cd,x_N), 
\]
we obtain  $(0,x_2,\cd,x_N)\not\in 
(\wtil B^-_{w_0})_{\Phi^{(+)},\Theta_{\bfi}}$. 
Therefore, in both cases we find that
$\eit(\Psi^{(+)}_i(x))=\Psi^{(+)}_i(\eit (x))=0$.
The case for $\fit$ can be done similarly. 
Now we know that the condition (v) holds and 
we have completed the proof of Theorem \ref{b-ify}.\qed

\renewcommand{\thesection}{\arabic{section}}
\section{Explicit form of 
$\Phi^{(+)}(\Theta_{\bf i_0}^-(c))$ for types $\AAN,\BBN,\CCN$, $\DDN$ and 
$\G22$}
\label{fb-ABCDG}
\setcounter{equation}{0}
\renewcommand{\theequation}{\thesection.\arabic{equation}}

The results of this section and the next section will be used  to show the one 
of our main results Theorem \ref{main},  namely, 
the connectedness of cellular crystals $B_\bfi$, where $\bfi$ 
is an arbitrary reduced word of the longest element $w_0\in W$. 
To complete the task, we need to find 
certain set of linear functions $\Xi_\bfi\subset (\bbQ^N)^*$ 
satisfying the condition ${\bf H}_\bfi$ as in Sect.8 and to obtain 
this $\Xi_\bfi$ we need to get 
some explicit forms of the upper half potentials for all types of 
finite-dimensional simple Lie algebras.
Then, in this section we shall see some results given in \cite{N6}
for all classical types $\AAN,\BBN,\CCN$, $\DDN$ and the exceptional types $\G22$. 
Similar results for type 
$\EEE6,\EEE7,\EEE8,\F44$ will be given in the next section.

For those simple Lie algebras, 
let us fix the numbering of Dynkin diagrams. 
For all types except ${\rm E}_7$, we adopt the ones in \cite[Chap.4]{Kac}.
As for the type ${\rm E}_7$, the Dynkin data will be given in the
next section.
We fix the reduced longest word $\bfii0$ as follows:
\begin{equation}
\bfii0
=\begin{cases}\underbrace{1,2,\cd,n}_{},
\underbrace{1,2,\cd,n-1}_{},\cd,\underbrace{1,2,3}_{},1,2,1&\text{ for
  }\AAN,\\
(1,2,\cd,n-1,n)^n&\text{ for }\BBN,\,\,\CCN,\\
(1,2,\cd,n-1,n)^{n-1}&\text{ for }\DDN,\\
121212&\text{ for }\G22.
\label{i0}
\end{cases}
\end{equation}
According to the reduced word $\bfii0$, we consider 
the set of double-indices for type $\AAN,\BBN,\CCN,\DDN$:
\[
I_n:=\begin{cases}
\{(1,1),(1,2),\cd,(1,n),(2,1),(2,2),\cd,(2,n-1)\cd,(n-1,1),(n-1,2),(1,n)\}
&\AAN,\\
\{(1,1),(1,2),\cd,(1,n),(2,1),\cd,(2,n),\cd,(n,1),\cd,(n,n)\}&\BBN,\CCN,\\
\{(1,1),(1,2),\cd,(1,n),(2,1)\cd,(2,n),\cd,(n-1,1),\cd\,(n-1,n)\}&\DDN.
\end{cases}
\]
We also define the total
 order on $I_n$ by the usual lexicographic order, {\it e.g.,}
for type $\BBN$, we have 
\[
(1,1)<(1,2)<\cd<(1,n)<\cd<(n-1,1)<\cd<(n-1,n)<(n,1)\cd<(n,n).
\]
\subsection{Type $\AAN$}
To obtain the explicit form of $\TY(\Phi,+,)(\Theta_{\bf i_0}^-(c))$ for type
$\AAN$, it suffices to know 
$\Del_{w_0\Lm_i,s_i\Lm_i}(\Theta_{\bf i_0}^-(c))$ for
\[
c=(c^{(i)}_j|i+j\leq n+1)=
(c^{(1)}_1,c^{(1)}_2,\cd,c^{(1)}_n,c^{(2)}_1,c^{(2)}_2,\cd,c^{(2)}_{n-1},\cd
c^{(n-1)}_1,c^{(n-1)}_2,c^{(n)}_1)\in (\bbC^\times)^N.
\]
\begin{thm}[\cite{N6}]
\label{thm-a}
For $c\in (\bbC^\times)^N$ as above, 
we have the following explicit forms:
\begin{eqnarray}
&&\Del_{w_0\Lm_j,s_j\Lm_j}(\Theta_\bfii0^-(c))=
c^{(n-j+1)}_1+\frac{c^{(n-j+1)}_2}{c^{(n-j+2)}_1}
+\frac{c^{(n-j+1)}_3}{c^{(n-j+2)}_2}
+\cd+\frac{c^{(n-j+1)}_j}{c^{(n-j+2)}_{j-1}}, 
\label{del1}
\end{eqnarray}
\end{thm}

The relation between the decoration and the 
monomial realizations of crystals of 
type $\AAN$ is explicitly given in \cite{N6,N4}. 
For type $\AAN$ take $p_0:=(p_{i,j})_{i,j\in I,i\ne j}$
such that $p_{i,j}=1$ for $i<j$, $p_{i,j}=0$
for $i>j$, which corresponds to the cyclic sequence
${\mathbf i}=(12\cd n)(12\cd n)\cd$.
The crystal $B(\Lm_1)$ is described as follows (\cite{KN}):
$B(\Lm_1):=\{v_i|1\leq i\leq n+1\}$
and the actions of $\eit$ and $\fit$ are given as
\begin{eqnarray}
&&\fit v_i=v_{i+1},\,\eit v_{i+1}=v_i,
\q(1\leq i\leq n),\q
\til f_nv_{n+1}=\til e_1v_{1}=0,\label{a-f2}
\end{eqnarray}
and the other actions are trivial.
To describe the monomial realizations, we write down the monomials $A_{m,i}$
associated with $p_0$:
\begin{equation}
A_{m,i}=\begin{cases}
\ci(1,m){\ci(2,m)}^{-1}\ci(1,m+1),&\text{ for }i=1,\\
\ci(i,m){\ci(i+1,m)}^{-1}{\ci(i-1,m+1)}^{-1}\ci(i,m+1),
&\text{ for }1<i\leq n-1,\\
\ci(n,m){\ci(n-1,m+1)}^{-1}\ci(n,m+1)&\text{ for }i=n.
\end{cases}
\label{aim-a-1}
\end{equation}

Then we obtain 
\begin{pro}[\cite{N6}]
\label{mono-cry}
The crystal containing the monomial $Y_{n-i+1,1}$ (resp. $Y_{i,1}^{-1}$)
is isomorphic to $B(\Lm_1)$ (resp. $B(\Lm_n)$)
and all basis vectors are given by
\begin{eqnarray*}
&&
\til f_k\cd \til f_2\til
f_1(Y_{n-i+1,1})=\frac{Y_{n-i+1,k+1}}{Y_{n-i+2,k}}\in B(\Lm_1)\q
(k=1,\cd,n).
\end{eqnarray*}
\end{pro}

Applying this results to Theorem \ref{thm-a} and 
changing the variable $Y_{m,l}$ to $c^{(m)}_l$, we find:
\begin{thm}[\cite{N6}]
\label{del-mono}
For $j=1,\cd,n$  we have
\begin{eqnarray*}
&&
\Del_{w_0\Lm_j,s_j\Lm_j}(\Theta_\bfii0^-(c))=
\sum_{k=0}^{j-1}\til f_k\cd \til f_2\til f_1(c^{(n-j+1)}_1).
\end{eqnarray*}
\end{thm}


\subsection{Type $\CCN$}
Here let us see the results for type $\CCN$.
\begin{thm}[\cite{N6}]\label{thm-c-1}
In the case $\CCN$,  for $\bfii0=(12\cd n)^n$ and 
$c=(\ci(i,j))_{1\leq i,j\leq n}\\
=(\ci(1,1),\ci(2,1),\cd,\ci(n-1,n),\ci(n,n))\in
(\bbC^\times)^{n^2}$
for $k=1,2,\cd,n-1$ we have
\begin{eqnarray*}
&&\Del_{w_0\Lm_k,s_k\Lm_k}(\Theta^-_\bfii0(c))
=\ci(1,k)+\frac{\ci(2,k)}{\ci(1,k+1)}+\cd+
\frac{\ci(n,k)}{\ci(n-1,k+1)}+
\frac{\ci(n-1,k+1)}{\ci(n,k+1)}+\frac{\ci(n-2,k+2)}{\ci(n-1,k+2)}+\cd+
\frac{\ci(k,n)}{\ci(k+1,n)},\\
&&\Del_{w_0\Lm_n,s_n\Lm_n}(\Theta^-_\bfii0(c))=\ci(n,n).
\end{eqnarray*}
\end{thm}

We see the monomial realization of $B(\Lm_1)$ associated with 
the cyclic order $\cd(12\cd n)(12\cd n)\cd$, which means that the sign
$p_0=(p_{i,j})$ is given by $p_{i,j}=1$ if $i<j$ and $p_{i,j}=0$ if 
$i>j$. 
The crystal $B(\Lm_1)$ is described as follows (\cite{KN}):
$B(\Lm_1):=\{v_i,v_{\ovl i}|1\leq i\leq n\}$
and the actions of $\eit$ and $\fit$ are given as
\begin{eqnarray}
&&\fit v_i=v_{i+1},\q \fit v_{\ovl i+1}=v_{\ovl i},\q
\eit v_{i+1}=v_i,\q \eit v_{\ovl i}=v_{\ovl{i+1}}
\q(1\leq i<n),\label{c-f1}\\
&&\til f_nv_n=v_{\ovl n},\qq 
\til e_nv_{\ovl n}=v_n,\label{c-f2}
\end{eqnarray}
and the other actions are trivial.
To describe the monomial realizations, we write down the monomials $A_{m,i}$
associated with $p_0$:
\begin{equation}
A_{m,i}=\begin{cases}
\ci(1,m){\ci(2,m)}^{-1}\ci(1,m+1),&\text{ for }i=1,\\
\ci(i,m){\ci(i+1,m)}^{-1}{\ci(i-1,m+1)}^{-1}\ci(i,m+1),
&\text{ for }1<i\leq n-1,\\
\ci(n,m){\ci(n-1,m+1)}^{-2}\ci(n,m+1)&\text{ for }i=n.
\end{cases}
\label{aim-c-1}
\end{equation}
Here the monomial realization of $B(\Lm_1)$ is described explicitly:
\begin{equation}
B(\ci(1,k))(\cong B(\Lm_1))=\left\{\frac{\ci(j,k)}{\ci(j-1,k+1)}=m_1^{(k)}(v_j),
\frac{\ci(j-1,k+n-j+1)}{\ci(j,k+n-j+1)}=m_{1}^{(k)}(v_{\ovl j})\,\,
|\,\,1\leq j\leq n\,\,\right\},
\end{equation}
where $m^{(k)}_i:B(\Lm_i)\hookrightarrow \cY(p)$ ($u_{\Lm_i}\mapsto 
\ci(i,k)$)
is the embedding of crystal as in \eqref{embed-mki}
and we understand $\ci(0,k)=1$.
Now, Theorem \ref{thm-c-1} claims the following:
\begin{thm}\label{thm-mono-c1}
We obtain  $\Del_{w_0\Lm_n,s_n\Lm_n}(\Theta^-_\bfii0(c))
=m_n^{(n)}(u_{\Lm_n})=\ci(n,n)$
 and 
\begin{eqnarray}
&&\Del_{w_0\Lm_k,s_k\Lm_k}(\Theta^-_\bfii0(c))
=\sum_{j=1}^nm^{(k)}_1(v_j)+\sum_{j=k+1}^nm^{(k)}_1(v_{\ovl j}),\qq
(k=1,2\cd,n-1)
\end{eqnarray}
\end{thm}

\subsection{Type $\BBN$}
Here we shall see the results for type $\BBN$.  Let us take  the sequence 
$\bfii0=(12\cd n)^n$ as above.
\begin{thm}\label{thm-Bn1}
For 
$c=(\ci(j,i))_{1\leq,i,j\leq n}=
(\ci(1,1),\ci(2,1),\cd,\ci(n-1,n),\ci(n,n))\in(\bbC^\times)^{n^2}$, we have
\begin{eqnarray*}
&&\hspace{-20pt}\Del_{w_0\Lm_k,s_k\Lm_k}(\Theta^-_\bfii0(c))\\
&&\hspace{-20pt}=
\ci(1,k)+\frac{\ci(2,k)}{\ci(1,k+1)}+\cd+
\frac{\ci(n-1,k)}{\ci(n-2,k+1)}+\frac{{\ci(n,k)}^2}{\ci(n-1,k+1)}
+2\frac{\ci(n,k)}{\ci(n,k+1)}+\frac{\ci(n-1,k+1)}{{\ci(n,k+1)}^2}
+\frac{\ci(n-2,k+2)}{\ci(n-1,k+2)}+\cd+
\frac{\ci(k,n)}{\ci(k+1,n)}\,\,(k=1,\cd,n-1),\\
&&\hspace{-20pt}\Del_{w_0\Lm_n,s_n\Lm_n}(\Theta^-_\bfii0(c))=\ci(n,n).
\end{eqnarray*}
\end{thm}


We see the monomial realization of $B(\Lm_1)$ associated with 
the cyclic order $\cd(12\cd n)(12\cd n)\cd$, which means that the sign
$p_0=(p_{i,j})$ is given by $p_{i,j}=1$ if $i<j$ and $p_{i,j}=0$ if 
$i>j$ as in the previous section.
The crystal $B(\Lm_1)$ is described as follows (\cite{KN}):
$B(\Lm_1):=\{v_i,v_{\ovl i}|1\leq i\leq n\}\sqcup\{v_0\}$ and 
the actions of $\eit$ and $\fit$ $(1\leq i<n)$ are given as 
\begin{eqnarray}
&&\fit v_i=v_{i+1},\q \fit v_{\ovl{i+1}}=v_{\ovl i},\q
\eit v_{i+1}=v_i,\q \eit v_{\ovl i}=v_{\ovl{i+1}}
\q(1\leq i<n),\label{b-f1}\\
&&\til f_nv_n=v_0,\q \til f_nv_0=v_{\ovl n},\q 
\til e_nv_0=v_n,\q \til e_nv_{\ovl n}=v_0,\label{b-f2}
\end{eqnarray}
and the other actions are trivial.

To see the monomial realization $B(\ci(1,k))$, 
we describe the monomials $A_{m,i}$ explicitly:
\begin{equation}
A_{m,i}=\begin{cases}
\ci(1,m){\ci(2,m)}^{-1}\ci(1,m+1)&\text{ for }i=1,\\
\ci(i,m){\ci(i+1,m)}^{-1}{\ci(i-1,m+1)}^{-1}\ci(i,m+1)
&\text{ for }1<i< n-1,\\
\ci(n-1,m){\ci(n,m)}^{-2}{\ci(n-2,m+1)}^{-1}
\ci(n-1,m+1)&\text{ for }i=n-1,\\
\ci(n,m){\ci(n-1,m+1)}^{-1}\ci(n,m+1)&\text{ for }i=n.
\end{cases}
\label{aim-b-1}
\end{equation}
Here the monomial realization $B(\ci(1,k))$ for 
 $B(\Lm_1)$ associated with 
$p_0$ is described explicitly:
\begin{equation}
B(\ci(1,k))=\left\{\frac{{\ci(j,k)}^{\ep_j}}{\ci(j-1,k+1)}=m^{(k)}_1(v_j),
\frac{\ci(n,k)}{\ci(n,k+1)}=m^{(k)}_1(v_0),
\frac{\ci(j-1,k+n-j+1)}{{\ci(j,k+n-j+1)}^{\ep_j}}=m^{(k)}_1(v_{\ovl j})\,\,
|\,\,1\leq j\leq n\,\,\right\},
\end{equation}
where 
$m^{(k)}_i:B(\Lm_i)\hookrightarrow \cY(p)$ ($u_{\Lm_i}\mapsto 
\ci(i,k)$)
is the embedding of crystal as in \eqref{embed-mki}
and we understand $\ci(0,k)=1$.
Now, Theorem \ref{thm-Bn1}  means the following:
\begin{thm}\label{thm-mono-b1}
We obtain $\Del_{w_0\Lm_n,s_n\Lm_n}(\Theta^-_\bfii0(c))
=m_n^{(n)}(u_{\Lm_n})(=\ci(n,n))$
 and 
\begin{eqnarray}
&&\Del_{w_0\Lm_k,s_k\Lm_k}(\Theta^-_\bfii0(c))
=\sum_{j=1}^nm^{(k)}_1(v_j)+2m^{(k)}_1(v_0)
+\sum_{j=k+1}^nm^{(k)}_1(v_{\ovl j}),\qq
(k=1,2\cd,n-1).
\end{eqnarray}
\end{thm}



\subsection{Type $\DDN$}
In case of type $\DDN$, take the cyclic reduced longest word
$\bfii0=(1 2 \cd n-1\,n)^{n-1}$ as above.
\begin{thm}\label{thm-Dn1}
For $k\in\{1,2,\cd,n\}$ and $c=(\ci(i,j))=(\ci(1,1),\ci(2,1),\cd,
\ci(n-1,n-1),\ci(n,n-1))\in(\bbC^\times)^{n(n-1)}$, we have
\begin{eqnarray*}
&&\Del_{w_0\Lm_k,s_k\Lm_k}(\Theta^-_\bfii0(c))\\
&&\hspace{-20pt}=
\ci(1,k)+\frac{\ci(2,k)}{\ci(1,k+1)}+\cd+
\frac{\ci(n-2,k)}{\ci(n-3,k+1)}
+\frac{\ci(n-1,k)\ci(n,k)}{\ci(n-2,k+1)}
+\frac{\ci(n-2,k+1)}{\ci(n-1,k+1)\ci(n,k+1)}+
\frac{\ci(n,k)}{\ci(n-1,k+1)}+\frac{\ci(n-1,k)}{\ci(n,k+1)}
+\frac{\ci(n-3,k+2)}{\ci(n-2,k+2)}+\cd+
\frac{\ci(k,n-1)}{\ci(k+1,n-1)}\\
&&\qq\qq (k=1,2,\cd,n-2),\\
&&\Del_{w_0\Lm_{n-1},s_{n-1}\Lm_{n-1}}(\Theta^-_\bfii0(c))=\ci(n-1,n-1),\qq
\Del_{w_0\Lm_{n},s_n\Lm_{n}}(\Theta^-_\bfii0(c))=\ci(n,n-1).
\end{eqnarray*}
\end{thm}


Let us see the monomial realization of $B(\Lm_1)$ associated with 
the cyclic sequence $\cd(12\cd n)(12\cd n)\cd$, which means that the sign
$p_0=(p_{i,j})$ is given by $p_{i,j}=1$ if $i<j$ and $p_{i,j}=0$ if 
$i>j$ as in the previous sections.
The crystal $B(\Lm_1)$ is described as follows (\cite{KN}):
$B(\Lm_1):=\{v_i,v_{\ovl i}|1\leq i\leq n\}$ and 
the actions of $\eit$ and $\fit$ $(1\leq i\leq n)$ are given as 
\begin{eqnarray}
&&\fit v_i=v_{i+1},\q \fit v_{\ovl{i+1}}=v_{\ovl i},\q
\eit v_{i+1}=v_i,\q \eit v_{\ovl i}=v_{\ovl{i+1}}
\q(1\leq i<n),\label{d-fi}\\
&&\til f_nv_n=v_{\ovl{n-1}},\q \til f_{n-1}v_{\ovl n}=v_{\ovl{n-1}},\q 
\til e_{n-1}v_{\ovl {n-1}}=v_{\ovl n},\q \til e_nv_{\ovl{n-1}}=v_n,
\label{d-fn}
\end{eqnarray}
and the other actions are trivial.
To see the monomial realization $B(\ci(1,k))$, 
we give the explicit forms of the monomials $A_{m,i}$:
\begin{equation}
A_{m,i}=\begin{cases}
\ci(1,m){\ci(2,m)}^{-1}\ci(1,m+1)&\text{ for }i=1,\\
\ci(i,m){\ci(i+1,m)}^{-1}{\ci(i-1,m+1)}^{-1}\ci(i,m+1)
&\text{ for }1<i< n-2,\\
\ci(n-2,m){\ci(n-1,m)}^{-1}{\ci(n,m)}^{-1}{\ci(n-3,m+1)}^{-1}
\ci(n-2,m+1)&\text{ for }i=n-2,\\
\ci(n-1,m){\ci(n-2,m+1)}^{-1}
\ci(n-1,m+1)&\text{ for }i=n-1,\\
\ci(n,m){\ci(n-2,m+1)}^{-1}\ci(n,m+1)&\text{ for }i=n.
\end{cases}
\label{aim-d-1}
\end{equation}
Let $m^{(k)}_i:B(\Lm_i)\hookrightarrow \cY(p)$ ($u_{\Lm_i}\mapsto 
\ci(i,k)$)
is the embedding of crystal as in \eqref{embed-mki}.
Here the monomial realization $B(\ci(1,k))=m^{(k)}_1(B(\Lm_1))=
\{m^{(k)}_1(v_j),m^{(k)}_1(v_{\ovl j})|1\leq j\leq n\}$ associated with 
$p_0$ is described explicitly:
\begin{equation}
m^{(k)}_1(v_j)=\begin{cases}
\frac{\ci(j,k)}{\ci(j-1,k+1)}&1\leq j\leq n-2,\\
\frac{\ci(n-1,k)\ci(n,k)}{\ci(n-2,k+1)}&j=n-1,\\
\frac{\ci(n,k)}{\ci(n-1,k+1)}&j=n,
\end{cases}\qq\qq
m^{(k)}_1(v_{\ovl j})=\begin{cases}
\frac{\ci(j-1,k+n-i)}{\ci(j,k+n-j)}&1\leq j\leq n-2,\\
\frac{\ci(n-2,k+1)}{\ci(n-1,k+1)\ci(n,k+1)}&j=n-1,\\
\frac{\ci(n-1,k)}{\ci(n,k+1)}&j=n,
\end{cases}
\end{equation}
where we understand $\ci(0,k)=1$.
Now, Theorem \ref{thm-Dn1} claims the following:
\begin{thm}\label{thm-mono-d1}
We obtain 
\begin{eqnarray}
&&\Del_{w_0\Lm_k,s_k\Lm_k}(\Theta^-_\bfii0(c))
=\sum_{j=1}^nm^{(k)}_1(v_j)+\sum_{j=k+1}^nm^{(k)}_1(v_{\ovl j}),\qq
(k=1,2\cd,n-2)
\\
&&\Del_{w_0\Lm_{n-1},s_{n-1}\Lm_{n-1}}(\Theta^-_\bfii0(c))
=m^{(n-1)}_{n-1}(u_{\Lm_{n-1}}),\qq
\Del_{w_0\Lm_{n},s_n\Lm_{n}}(\Theta^-_\bfii0(c))=m^{(n-1)}_n(u_{\Lm_n}).
\end{eqnarray}
\end{thm}

\subsection{Type $\G22$}

For $\bfii0=121212$, by direct calculations 
we have $\Del_{w_0{\Lm_2},s_2{\Lm_2}}(\Theta_\bfii0^-(c_1,\cd,c_6))
=c_6 $ and 
\begin{eqnarray*}
&&\Del_{w_0{\Lm_1},s_1{\Lm_1}}(\Theta_\bfii0^-(c_1,\cd,c_6))=
c_1+\frac{c_2^3}{c_3}+\frac{3c_2^2}{c_4}+\frac{3c_2c_3}{c_4^2}\\
&&
+\frac{c_3^2}{c_4^3}+\frac{2c_3}{c_5}+\frac{c_4^3}{c_5^2}
+\frac{3c_2c_4}{c_5}+\frac{3c_2}{c_6}+\frac{3c_3}{c_4c_6}
+\frac{3c_4^2}{c_5c_6}+\frac{3c_4}{c_6^2}+\frac{c_5}{c_6^3}.
\end{eqnarray*}

For type $\G22$ take $p_0:=(p_{1,2},p_{2,1})$ with 
$p_{1,2}=1$ and, $p_{2,1}=0$,  which corresponds to the cyclic sequence
$\bfii0=121212$.
To describe the monomial realizations, we write down the monomials $A_{m,i}$
associated with $p_0$:
\begin{equation}
A_{m,i}=\begin{cases}
c_mc_{m+1}^{-3}c_{m+2}&\text{ for }i=1,\\
c_mc_{m+1}^{-1}c_{m+2}&\text{ for }i=2,\\
\end{cases}
\label{aim-g-1}
\end{equation}
where note that we identify the variable $Y_{m,i}$ with $c_{2m+i}$. Then, 
the crystal graph of $B(\Lm_1)$. 
\[\xymatrix{
\linethickness{30pt}
c_1\ar@{->}[r]^1 &
\frac{c_2^3}{c_3}\ar@{->}[r]^2&
\frac{c_2^2}{c_4}\ar@{->}[r]^2&
\frac{c_2c_3}{c_4^2}
\ar@/^/[drrr]^1\ar@{->}[r]^2&
\frac{c_3^2}{c_4^3}\ar@{->}[r]^1&
\frac{c_3}{c_5}\ar@{->}[r]^1&
\frac{c_4^3}{c_5^2}\ar@/^/[dlll]^2\\
\frac{1}{c_7}&
\frac{c_5}{c_6^3}\ar@{->}[l]^1&
\frac{c_4}{c_6^2}\ar@{->}[l]^2&
\frac{c_4^2}{c_5c_6}\ar@{->}[l]^2&
\frac{c_3}{c_4c_6}\ar@{->}[l]^1&
\frac{c_2}{c_6}\ar@{->}[l]^2&
\frac{c_2c_4}{c_5}\ar@{->}[l]^2
}.
\]
\def\tf{\til f_1}
\def\tff{\til f_2}
Here note that $\frac{1}{c_7}$ does not appear in 
$\Del_{w_0{\Lm_1},s_1{\Lm_1}}(\Theta_\bfii0^-(c_1,\cd,c_6))$ and then 
one has the following:
\begin{thm}\label{mono-g2}
\begin{eqnarray*}
&&\Del_{w_0{\Lm_1},s_1{\Lm_1}}(\Theta_\bfi^-(c_1,\cd,c_6))=
c_1+\tf(c_1)+3\tff\tf(c_1)+3\tff^2\tf(c_1)\\
&&
+\tff^3\tf(c_1)+2\tf\tff^3\tf(c_1)+\tf^2\tff^3\tf(c_1)
+3\tf\tff^2\tf(c_1)+3\tff\tf\tff^2\tf(c_1)+3\tff^2\tf\tff^2\tf(c_1)\\
&&
+3\tf\tff^2\tf\tff^2\tf(c_1)+3\tff\tf\tff^2\tf\tff^2\tf(c_1)
+\tff^2\tf\tff^2\tf\tff^2\tf(c_1).
\end{eqnarray*}
\end{thm}

\subsection{}

Here summarizing the above results, one has the following theorem, 
which will be needed in Sect.8 to show the connectedness of the 
crystal $B_\bfi$:
\begin{thm}
For type $\AAN,\,\BBN,\,\CCN,\,\DDN$ and $\G22$,
 and the specific reduced longest word
$\bfii0$ as in \eqref{i0}, one gets that the generalized minor
$\Del_{w_0\Lm_k,s_k\Lm_k}(\Theta^-_\bfii0(c_1,\cd,c_N))$ ($k\in I$)
is a positive sum of monomials
in the form:
\begin{equation}
c_j\prod_{(m,i)\in T_\bfii0}A_{m,i}^{l_{m,i}},
\end{equation}
for some $j\in\{1,2,\cd,N=L(w_0)\}$ and some integers $l_{m,i}$. 
Note that the monomial $A_{m,i}$ is given in \eqref{aim-a-1},
\eqref{aim-c-1},\eqref{aim-b-1},\eqref{aim-d-1} and \eqref{aim-g-1}.
\end{thm}
{\sl Proof.}
The actions of $\eit$ (resp. $\fit$)
in the monomial realizations are obtained by 
multiplying (resp. dividing) the monomials $A_{m,i}$'s to monomials. 
Thus, it follows from Theorem \ref{del-mono}, Theorem \ref{thm-mono-c1},
Theorem \ref{thm-mono-b1}, 
Theorem \ref{thm-mono-d1} and Theorem \ref{mono-g2}
that the desired results.
\qed

\subsection{Comparison with the polyhedral realizations}

Here we shall see the relation between Theorem \ref{b-ify} and the 
polyhedral realization in \cite{NZ}.

First, let us recall the polyhedral realization of $B(\ify)$ 
for an arbitrary symmetrizable Kac-Moody Lie algebra $\ge$ whose 
Cartan matrix is $(\lan h_i,\al_j\ran)_{i,j\in I}$ (\cite{NZ}).
Consider the infinite dimensional vector space
$\QQ^{\ify}:=\{\textbf{a}=
(\cd,a_k,\cd,a_2,a_1): a_k \in \QQ\,\,{\rm and }\,\,
a_k = 0\,\,{\rm for}\,\, k \gg 0\}$ 
and its dual space $(\QQ^{\ify})^*:={\rm Hom}(\QQ^{\ify},\QQ)$.
Let $x_k\in (\QQ^{\ify})^*$ be the linear function defined as 
$x_k((\cd,a_k,\cd,a_2,a_1)):=a_k$
for $k\in\mathbb{Z}_{\geq1}$.
We will also write a linear form $\vp \in (\QQ^{\ify})^*$ as
$\vp=\sum_{k \geq 1} \vp_k x_k$ ($\vp_j\in \QQ$).

For the fixed infinite sequence
$\io=(i_k)$ and $k\geq1$ we set 
$\kp:={\rm min}\{l:l>k\,\,{\rm and }\,\,i_k=i_l\}$ and
$\km:={\rm max}\{l:l<k\,\,{\rm and }\,\,i_k=i_l\}$ if it exists,
or $\km=0$  otherwise.
We set $\beta_0=0$ and
\begin{equation}
\beta_k:=x_k+\sum_{k<j<\kp}\lan h_{i_k},\al_{i_j}\ran x_j+x_{\kp}\in (\QQ^{\ify})^*
\qq(k\geq1).
\label{betak}
\end{equation}
We define the piecewise-linear operator 
$S_k=S_{k,\io}$ on $(\QQ^{\ify})^*$ by
\begin{equation}
S_k(\vp):=
\begin{cases}
\vp-\vp_k\beta_k & {\mbox{ if }}\vp_k>0,\\
 \vp-\vp_k\beta_{\km} & {\mbox{ if }}\vp_k\leq 0.
\end{cases}
\label{Sk}
\end{equation}
Here we set
\begin{eqnarray}
\Xi_{\io} &:=  &\{S_{j_l}\cd S_{j_2}S_{j_1}x_{j_0}\,|\,
l\geq0,j_0,j_1,\cd,j_l\geq1\},
\label{Xi_io}\\
\Sigma_{\io} & := &
\{\textbf{x}\in \ZZ^{\ify}\subset \QQ^{\ify}\,|\,\vp(\textbf{x})\geq0\,\,{\rm for}\,\,
{\rm any}\,\,\vp\in \Xi_{\io}\}.
\end{eqnarray}
We impose on $\io$ the following {\it positivity condition}:
\begin{equation}
{\hbox{if $\km=0$ then $\vp_k\geq0$ for any 
$\vp=\sum_k\vp_kx_k\in \Xi_{\io}$}}.
\label{posi}
\end{equation}
\begin{thm}[\cite{NZ}]\label{polyhthm}
Let $\io$ be a sequence of indices satisfying $(\ref{seq-con})$ 
and $(\ref{posi})$ and $\Psi_\io:B(\ify)\to \bbZ^\ify_\io$ 
the Kashiwara embedding associated with $\io$. 
 Then we have 
${\rm Im}(\Psi_{\io})(\cong B(\ify))=\Sigma_{\io}$,
which is called 
the {\it polyhedral realization} of the crystal $B(\ify)$ for $\ge$.
\end{thm}
In \cite{N7}, it is shown that for semi-simple $\ge$, the length of 
a sequence $\io$ can be finite, indeed, it is just $N=l(w_0)$.
We obtain the following (similar results were given in \cite{N4,N9}):
\begin{thm}\label{comparison}
Let $\ge$ be a finite-dimensional simple Lie algebra, $\bfi$ 
a reduced longest word  and $\io=\bfi^{-1}$ the opposite ordered longest word 
for the Langlands dual Lie algebra $\ge^\vee$
and $\Psi_\io:B(\ify)\to \bbZ^N_\io$ the Kashiwara embedding 
for $\ge^\vee$ associated with $\io$.
\begin{enumerate}
\item
One has the following identification as crystals in $\bbZ^N$:
\begin{eqnarray}
\Omega\,\,:\,\,(\bbB^-_{w_0})_{\Phi^{(+)},\Theta_{\bfi}}
&\xrightarrow{\q \sim \q}& {\rm Im}(\Psi_\io)
\q (\cong B(\ify)),\label{identif}\\
(x_1,x_2,\cd,x_N)&\longmapsto& 
(x_N,\cd,x_2,x_1)\nn
\end{eqnarray}
\item
Assume that $\io$ satisfies the the positivity condition \eqref{posi} and let
$\Sigma_\io^\vee$ be the polyhedral realization 
for $\ge^\vee$ associated with $\io$. Then the map $\Omega$ gives the 
identification of crystals:
\begin{eqnarray}
\Omega\,\,:\,\,(\bbB^-_{w_0})_{\Phi^{(+)},\Theta_{\bfi}}
&\xrightarrow{\q \sim \q}& \Sigma_{\io}^\vee
\q (\cong B(\ify)),
\end{eqnarray}

\end{enumerate}
\end{thm}
\nd
{\sl Proof.} It is trivial that 
(ii) is derived from (i) and Theorem \ref{polyhthm}. Thus, 
we shall show (i).
First we assume that the map $\Omega$ is well-defined.
Let us show :
\begin{enumerate}
\item[(1)]
$\wt(\Omega(x))=\wt(x)$ 
for any $x\in\bbZ^N$.
\item[(2)]
$\vep_i(\Omega(x))=\vep_i(x)$ for any $i\in I$, 
$x\in (\bbB^-_{w_0})_{\Phi^{(+)},\Theta_{\bfi}}$.
\item[(3)]
$\eit(\Omega(x))=\wtil\Omega(\eit(x))$,
$\fit(\Omega(x))=\wtil\Omega(\fit(x))$ for any $i\in I$,
$x\in (\bbB^-_{w_0})_{\Phi^{(+)},\Theta_{\bfi}}$,.
\end{enumerate}
It is easy to see (1) and (2) by comparing 
\eqref{wt-vep-vp-ify}, \eqref{wt-ud} and \eqref{vep-ud}, where note that 
we take Langlands dual Cartan data for ${\rm Im}(\Psi_\io)$.
Next, let us see (3).
Comparing with $\sigma_k$ in \eqref{sigma-def} and $X_m$ in Lemma \ref{lem-ef} 
we find that $\wtil M^{(i)}(\wtil\Omega(x))=M^{(i)}(x)$ for 
any $i\in I$,
$x\in (\bbB^-_{w_0})_{\Phi^{(+)},\Theta_{\bfi}}$, which implies 
$\til m^{(i)}_f(\Omega(x))=m^{(i)}_f(x)$ and 
$\til m^{(i)}_e(\Omega(x))=m^{(i)}_e(x)$ by \eqref{mfme} and the definitions
in Lemma \ref{lem-ef}. Therefore, it follows from
\eqref{fit-Z-ify}, \eqref{eit-Z-ify},
\eqref{th-eaction} and \eqref{th-faction} that 
$\eit(\Omega(x))=\Omega(\eit(x))$,
$\fit(\Omega(x))=\Omega(\fit(x))$.
 
Finally,  we shall show the well-definedness of the map $\Omega$ by induction
on the height ${\rm ht}(x):=\sum_{k=1}^N|x_k|$ for 
$x\in (\bbB^-_{w_0})_{\Phi^{(+)},\Theta_{\bfi}}$.
Note that both  $(\bbB^-_{w_0})_{\Phi^{(+)},\Theta_{\bfi}}$
and ${\rm Im}(\Psi_\io)$ are generated from $(0,0,\cd,0)$ by applying 
$\fit$'s.
First, for $x_0=(0,0,\cd,0)
\in (\bbB^-_{w_0})_{\Phi^{(+)},\Theta_{\bfi}}$, it is trivial that 
$\Omega(x_0)=(0,0,\cd,0)\in {\rm Im}(\Psi_\io)$. 
Thus, in the case ${\rm ht}(x)=0$, the map $\Omega$ is well-defined.
Here, suppose that for any $x$ with ${\rm ht}(x)=h>0$, the map $\Omega$ is 
well-defined, that is, 
$\Omega(x)\in {\rm Im}(\Psi_\io)$. Then, by (iii) we can easily show that 
$\Omega(\fit(x))=\fit\Omega(x)\in {\rm Im}(\Psi_\io)$, which means 
that $\Omega$ is well-defined for any 
$y\in (\bbB^-_{w_0})_{\Phi^{(+)},\Theta_{\bfi}}$ with ${\rm ht}(y)=h+1$.
Thus, the induction step proceeds and then we 
show that the well-definedness of 
the map $\Omega$.\qed

\renewcommand{\thesection}{\arabic{section}}
\section{Trail and Half Potentials for types $\EEE6,\EEE7,\EEE8$ and $\F44$}
\label{sec-trail}

\setcounter{equation}{0}
\renewcommand{\theequation}{\thesection.\arabic{equation}}

In this section, we review the theory of trail by Berenstein and Zelevinsky
\cite{BZ2} and apply it to obtain 
properly explicit forms of the upper half potential $\Phi^{(+)}$ 
for the exceptional types $\EEE6,\EEE7,\EEE8$ and $\F44$, which 
are not completely explicit as in the previous section, but 
sufficient for later use.
For a reduced word $\bfi$, we also consider the set of 
double-indices and the total order on it like as the ones in the
previous section.

\subsection{Trails}

\begin{df}\label{pretrail}
For a finite dimensional representation $V$ of $\mathfrak{g}$, 
two weights $\gamma$, $\delta$ of $V$ and
a sequence $\bfi=(i_1,\cdots,i_l)$ of indices from $I$,
we say a sequence 
$\pi=(\gamma=\gamma_0,\gamma_1,\cdots,\gamma_l=\delta)$ 
is a {\it pre-$\bfi$-trail} from $\gamma$ to $\delta$ if
for $1\leq k\leq l$, one has $\gamma_k\in P$ , 
$\gamma_{k-1}-\gamma_k=c_k\alpha_{i_k}$ 
for some non-negative integer $c_k$.
\end{df}

\begin{df}\cite{BZ2}
In the setting of Definition \ref{pretrail},
if pre-$\bfi$-trail $\pi$ satisfies the condition (TR):
\begin{itemize}
\item[(TR)] 
$e^{c_1}_{i_1}e^{c_2}_{i_2}\cdots e^{c_l}_{i_l}$ 
is a non-zero linear map from $V_\delta$ to $V_\gamma$,
\end{itemize}
then $\pi$ is called an {\it $\bfi$-trail} from $\gamma$ to $\delta$,
where $V=\oplus_{\mu\in P} V_\mu$ is the weight decomposition of $V$.
\end{df}

\nd
For a pre-$\bfi$-trail $\pi$, we set
\[
d_k(\pi):=\frac{\gamma_{k-1}+\gamma_k}{2}(h_{i_k}).
\]
Note that if $\gamma_{k-1}=s_{i_k}\gamma_k$ then $d_k(\pi)=0$.
For a reduced word $\bfi=(i_1,\cdots,i_l)$ of an element in $W$, 
as before we define
$\Theta^{-}_{\bfi}(c_1,\cdots,c_l)
:=\textbf{y}_{i_1}(c_1)\cdots \textbf{y}_{i_l}(c_l)$.

\begin{thm}[\cite{BZ2}Theorem 5.8]\label{trail-thm}
For $u$, $v\in W$ and $i\in I$, we obtain that
$\Delta_{u\Lambda_i,v\Lambda_i}(\Theta^-_{\bfi}(t_1,\cdots,t_l))$
is a linear combination of the monomials $t_1^{d_1(\pi)}\cdots t_l^{d_l(\pi)}$ with positive coefficients for
all $\bfi$-trail from $-u\Lambda_i$ to $-v\Lambda_i$ in 
$V(-w_0(\Lm_i))$.
\end{thm}

\begin{thm}\label{lowest-term-thm}
For $j\in I$ and a reduced word $\bfi=(i_1,\cdots,i_N)$ of $w_0$, 
we see that
$\Delta_{w_0\Lambda_j,s_j\Lambda_j}(\Theta^-_{\bfi}(t_1,\cdots,t_N))$ has 
a term 
\begin{equation}\label{lowest-term}
t_{J} t_{J+1}^{{\bf a}_{i_{J+1},j}}\cdots t_{N}^{{\bf a}_{i_N,j}},
\end{equation}
where $J:={\rm max}\{1\leq k\leq N | i_k=j\}$. Note that this monomial 
has a positive integer coefficient.
\end{thm}

\nd
{\sl  Proof.}
There is an $\bfi$-trail $\pi_0$ from $-w_0\Lambda_j$ to $-s_j\Lambda_j$ such that
$\gamma_N=\gamma_{N-1}=\cdots=\gamma_{J}=\gamma_{J-1}=-s_j\Lambda_j=-\Lambda_j+\alpha_j$ and
$\gamma_{k-1}=s_{i_k}\gamma_k$ for $k\in[1,J-1]$. Then we get
\[
d_k(\pi_0)=
\begin{cases}
0 & {\rm if}\ k\in[1,J-1], \\
(-\Lambda_j+\alpha_j)(h_{i_k}) & {\rm if}\ k\in[J,N].
\end{cases}
\]
We also get
\[
(-\Lambda_j+\alpha_j)(h_{i_k})=
\begin{cases}
1 & {\rm if}\ k=J,\\
{\bf a}_{i_k,j} & {\rm if}\ k\in[J+1,N].
\end{cases}
\]
By Theorem \ref{trail-thm}, we obtain our claim. \qed

\subsection{Calculations of the generalized minors 
$\Delta_{w_0\Lambda_i,s_i\Lambda_i}(\Theta^-_{\bfi}(\textbf{c}))$}

We suppose $\ge$ is of type $\EEE6,\EEE7,\EEE8$ or $\F44$.
Here, the Dynkin diagram of type $\EEE7$ is different from the one in 
\cite{Kac} but as follows:
The Cartan matrix 
 $({\bf a}_{ij})_{1\leq i,j\leq 7}$ of type $\EEE7$ is given by 
\begin{equation}
{\bf a}_{ij}=\begin{cases}
2&{\rm if }i=j,\\
-1&{\rm if}\ |i-j|=1\ {\rm or }\ (i,j)=(4,7),\,(7,4),\\
0&{\rm otherwise}.
\end{cases}
\end{equation}
%
%
%
%

Let $\bfii0$ be the following longest reduced word 
of the longest element $w_0\in W$:
\begin{equation}
\bfii0=
\begin{cases}
((1,\cdots,6)^4,1,2,3,4,6,1,2,3,6,1,2,1) & {\rm if}\ \EEE6, \\
(1,\cdots,7)^9 & {\rm if}\ {\EEE7},\\
(1,\cdots,8)^{15} & {\rm if}\ {\EEE8},\\
(1,2,3,4)^6 & {\rm if}\ {\F44}.
\end{cases}
\label{EF-i}
\end{equation}
\begin{rem}\label{rem1}
We remark that for each  $\bfii0$ as above, between any two consecutive 
occurrences of the letter $i\in I$,
each index $j$ satisfying ${\bf a}_{i,j}<0$ appears exactly once.
This property will be used later.
\end{rem}
We identify the set $\{1,2,\cdots,N=I(w_0)\}$ 
 with the doubly-index set 
\[
T_\bfii0:=\begin{cases}
\left\{\begin{array}{l}
(1,1),\cdots,(4,6),(5,1),\cd,(5,4),(5,6),\\
(6,1),(6,2),(6,3),(6,6),(7,1),(7,2),(8,1)
\end{array}\right\}&\EEE6,\\
\{(1,1),(1,2),\cd,(1,7),\cd,(9,1),\cd,(9,7)\}&\EEE7,\\
\{(1,1),(1,2),\cd,(1,8),\cd,(15,1),\cd,(15,8)\}&\EEE8,\\
\{(1,1),(1,2),\cdots,(6,1),\cdots,(6,4)\}&\F44
\end{cases}
\]
\begin{rem}\label{rem2}
If we replace the $m$-th occurrence of the letter $i$ in $\bfii0$ 
with the double index $(m,i)$
then we get $T_{\bfii0}$, which is lexicographically sorted.
Then we define the lexicographic  order $<$ on the set 
$T_\bfii0$ by identifying the order $1<2<3<\cdots<N$.
\end{rem}

We put the variable $\bf c$ indexed by $T_\bfii0$:
\def\cc(#1,#2){c^{(#1)}_{#2}}
\[
\textbf{c}=
\begin{cases}
(\cc(1,1),\cdots,\cc(1,6),\cc(2,1),\cdots,\cc(2,6),\cdots,\cc(4,1),\cdots,\cc(4,6),\\
\cc(5,1),\cc(5,2),\cc(5,3),\cc(5,4),\cc(5,6),\cc(6,1),\cc(6,2),\cc(6,3),\cc(6,6),\cc(7,1),\cc(7,2),\cc(8,1))
& {\rm if}\ \EEE6, \\
(\cc(1,1),\cdots,\cc(1,7),\cc(2,1),\cdots,\cc(2,7),
\cdots,\cc(9,1),\cdots,\cc(9,7)) & {\rm if}\ {\EEE7},\\
(\cc(1,1),\cdots,\cc(1,8),\cc(2,1),\cdots,\cc(2,8),
\cdots,\cc(15,1),\cdots,\cc(15,8)) & {\rm if}\ {\EEE8},\\
(\cc(1,1),\cc(1,2),\cc(1,3),\cc(1,4),\cc(2,1),\cc(2,2),\cc(2,3),\cc(2,4),
\cdots,\cc(6,1),\cc(6,2),\cc(6,3),\cc(6,4)) & {\rm if}\ {\F44},\\
\end{cases}
\]
For $s\in\mathbb{Z}$ and $i\in[1,n]$, we set
$p_{j,i}=1$ if $j<i$, $p_{j,i}=0$ if $j>i$ and
\begin{equation}
A_{s,i}:=\cc(s,i)\cc(s+1,i) 
\prod_{j\in[1,n]\setminus i}{\cc(s+p_{j,i},j)}^{{\bf a}_{j,i}}.
\label{asi-ef}
\end{equation}
By direct calculation, we obtain the following:
\begin{lem}\label{lem-ef}
\underline{Type ${\rm E}_6$-case:} 
The monomial of 
$\Delta_{w_0\Lambda_i,s_i\Lambda_i}(\Theta^-_{\bfii0}(\textbf{c}))$
in (\ref{lowest-term}) is equal to
\begin{equation}
\begin{cases}
\cc(8,1) & {\rm if}\ i=1,\\
\frac{\cc(7,2)}{\cc(8,1)}=\cc(7,1)A_{7,1}^{-1} & {\rm if}\ i=2,\\
\frac{\cc(6,3)}{\cc(6,6)\cc(7,2)} =\cc(3,1)
\left(
A_{6, 2} A_{5, 6} A_{6, 1} A_{5, 3} A_{4, 4} A_{3, 5} A_{5, 2} A_{4, 3} 
A_{3, 4} A_{3, 6} A_{3, 3} A_{3, 2} A_{3, 1}\right)^{-1}
& {\rm if}\ i=3,\\
\frac{\cc(5,4)}{\cc(6,3)}=\cc(2,1)
\left(
A_{5, 3} A_{5, 2} A_{4, 6} A_{5, 1} A_{4, 3} A_{3, 4} A_{2, 5} A_{4, 2} 
A_{3, 3} A_{2, 4} A_{2, 6} A_{2, 3} A_{2, 2} A_{2, 1}\right)^{-1}
& {\rm if}\ i=4,\\
\frac{\cc(4,5)}{\cc(5,4)}=\cc(1,1)
\left(
A_{4, 4} A_{4, 3} A_{4, 2} A_{3, 6} A_{4, 1} A_{3, 3} A_{2, 4} A_{1, 5} 
A_{3, 2} A_{2, 3} A_{1, 4} A_{1, 6} A_{1, 3} A_{1, 2} A_{1,1}
\right)^{-1}
 & {\rm if}\ i=5,\\
\cc(6,6) & {\rm if}\ i=6.
\end{cases}
\label{e6-c}
\end{equation}

\nd
\underline{Type ${\rm E}_7$-case:}  
The monomial of 
$\Delta_{w_0\Lambda_i,s_i\Lambda_i}(\Theta^-_{\bfii0}(\textbf{c}))$
in (\ref{lowest-term}) is equal to
\begin{equation}
\left\{
\begin{array}{lll}
&\frac{\cc(9,1)}{\cc(9,2)} =\cc(1,1)
\left(A_{8, 2}A_{7, 3}A_{6, 4}A_{5, 5}A_{4, 6}A_{5, 7}A_{5, 4}A_{4, 5}
A_{5, 3}A_{4, 4}A_{3, 7}A_{5, 2}A_{4, 3}A_{3, 4}A_{2, 5}A_{1, 6}
 A_{5, 1}\right.&\\
& \qq\qq 
\left.A_{4, 2}A_{3, 3}A_{2, 4}A_{1, 5}A_{1, 7}A_{1, 4}A_{1, 3}A_{1, 2}
 A_{1,1}\right)^{-1}&{\rm  if}\ i=1\\
&\frac{\cc(9,2)}{\cc(9,3)} =\cc(2,1)
\left(A_{8, 3}A_{7, 4}A_{6, 5}A_{5, 6}A_{6, 7}A_{6, 4}A_{5, 5}
A_{6, 3}A_{5, 4}A_{4,7}A_{6, 2}A_{5, 3}A_{4, 4}A_{3, 5}A_{2, 6}
 A_{6, 1}\right.&\\
& \qq\qq
\left.A_{5, 2}A_{4, 3}A_{3, 4}A_{2, 5}A_{2, 7}A_{2, 4}A_{2, 3}A_{2, 2}
 A_{2,1}\right)^{-1}&{\rm  if}\ i=2,\\
&\frac{\cc(9,3)}{\cc(9,4)} =\cc(3,1)
\left(A_{8, 4}A_{7, 5}A_{6, 6}A_{7, 7}A_{7, 4}A_{6, 5}
A_{7,3}A_{6, 4}A_{5,7}A_{7, 2}A_{6, 3}A_{5, 4}A_{4, 5}A_{3, 6}
 A_{7, 1}\right.&\\
& \qq\qq
\left.A_{6, 2}A_{5, 3}A_{4, 4}A_{3, 5}A_{3, 7}A_{3, 4}A_{3, 3}A_{3, 2}
 A_{3,1}\right)^{-1}&{\rm  if}\ i=3,\\
&\frac{\cc(9,4)}{\cc(9,5)\cc(9,7)} =\cc(4,1)
\left(A_{8, 5}A_{7, 6}A_{8, 7}A_{8, 4}A_{7, 5}
A_{8,3}A_{7, 4}A_{6,7}A_{8, 2}A_{7, 3}A_{6, 4}A_{5, 5}A_{4, 6}
 A_{8, 1}\right.&\\
& \qq\qq
\left.A_{7, 2}A_{6,3}A_{5, 4}A_{4, 5}A_{4, 7}A_{4, 4}A_{4, 3}A_{4, 2}
 A_{4,1}\right)^{-1}&{\rm  if}\ i=4,\\
&\frac{\cc(9,5)}{\cc(9,6)} =\cc(8,6){A_{8,6}}^{-1}
&{\rm  if}\ i=5,\\
&\cc(9,6) & {\rm if}\ i=6,\\
&\cc(9,7) & {\rm if}\ i=7.
\end{array}
\right.
\label{e7-c}
\end{equation}

\underline{Type ${\rm E}_8$-case:} 
The monomial of 
$\Delta_{w_0\Lambda_i,s_i\Lambda_i}(\Theta^-_{\bfii0}(\textbf{c}))$
in (\ref{lowest-term}) is equal to
\begin{equation}
\left\{
\begin{array}{lll}
\frac{\cc(15,1)}{\cc(15,2)} =&\cc(1,1)
\left(A_{14, 2}A_{13, 3}A_{12, 4}A_{11, 5}A_{10, 6}A_{9, 7}A_{10, 8}\right.&\\
& A_{10, 5}A_{9, 6}A_{10, 4}A_{9, 5}A_{8, 8}A_{10, 3}A_{9, 4}A_{8, 5}
 A_{7, 6}A_{6, 7}A_{10, 2}A_{9, 3}A_{8, 4}A_{7, 5}A_{6, 6}A_{6, 8}&\\
& A_{10, 1}A_{9, 2}A_{8, 3}A_{7, 4} A_{6, 5}^2A_{5, 6}A_{4, 7}A_{5, 8}
A_{6, 4}A_{5, 5}A_{4, 6}A_{6, 3}A_{5, 4}A_{4, 5}A_{3, 8}A_{6, 2}&\\
& \left.A_{5, 3}A_{4, 4}A_{3, 5}A_{2, 6}A_{1, 7}A_{6, 1}A_{5, 2}A_{4, 3}
A_{3, 4}A_{2, 5}A_{1, 6}A_{1, 8}A_{1, 5}A_{1, 4}A_{1, 3}A_{1, 2} A_{1, 1}
\right)^{-1}&{\rm if}\  i=1,\\
\frac{\cc(15,2)}{\cc(15,3)} =&\cc(2,1)
\left(A_{14, 3}A_{13, 4}A_{12, 5}A_{11, 6}A_{10, 7}A_{11, 8}A_{11, 5}\right.&\\
& A_{10, 6}A_{11, 4}A_{10, 5}A_{9, 8}A_{11, 3}A_{10, 4}A_{9, 5}A_{8, 6}
A_{7, 7}A_{11, 2}A_{10, 3}A_{9, 4}A_{8, 5}A_{7, 6}A_{7, 8}A_{11, 1}&\\
& A_{10, 2}A_{9, 3}A_{8, 4} A_{7, 5}^2A_{6, 6}A_{5, 7}A_{6, 8}A_{7, 4}
A_{6, 5}A_{5, 6}A_{7, 3}A_{6, 4}A_{5, 5}A_{4, 8}A_{7, 2}A_{6, 3}&\\
& \left.A_{5, 4}A_{4, 5}A_{3, 6}A_{2, 7}A_{7, 1}A_{6, 2}A_{5, 3}A_{4, 4}
A_{3, 5}A_{2, 6}A_{2, 8}A_{2, 5}A_{2, 4}A_{2, 3}A_{2, 2}A_{2, 1}
\right)^{-1}&{\rm if}\  i=2,\\
\frac{\cc(15,3)}{\cc(15,4)} =&\cc(3,1)
\left(A_{14, 4}A_{13, 5}A_{12, 6}A_{11, 7}A_{12, 8}A_{12, 5}A_{11, 6}\right.&\\
& A_{12, 4}A_{11, 5}A_{10, 8}A_{12, 3}A_{11, 4}A_{10, 5}A_{9, 6}
 A_{8, 7}A_{12, 2}A_{11, 3}A_{10, 4}A_{9, 5}A_{8, 6}A_{8, 8}A_{12, 1}&\\
& A_{11, 2}A_{10, 3}A_{9, 4} A_{8, 5}^2A_{7, 6}A_{6, 7}A_{7, 8}A_{8, 4}
 A_{7, 5}A_{6, 6}A_{8, 3}A_{7, 4}A_{6, 5}A_{5, 8}A_{8, 2}A_{7, 3}&\\
& \left. A_{6, 4}A_{5, 5}A_{4, 6}A_{3, 7}A_{8, 1}A_{7, 2}A_{6, 3}A_{5, 4}
A_{4, 5}A_{3, 6}A_{3, 8}A_{3, 5}A_{3, 4}A_{3, 3}A_{3, 2}A_{3, 1}
\right)^{-1}&{\rm if}\  i=3,\\
\frac{\cc(15,4)}{\cc(15,5)} =&\cc(4,1)
\left(
A_{14, 5}A_{13, 6}A_{12, 7}A_{13, 8}A_{13, 5}A_{12, 6}A_{13, 4}\right.&\\
& A_{12, 5}A_{11, 8}A_{13, 3}A_{12, 4}A_{11, 5}A_{10, 6}A_{9, 7}
A_{13, 2}A_{12, 3}A_{11, 4}A_{10, 5}A_{9, 6}A_{9, 8}A_{13, 1}&\\
& A_{12, 2}A_{11, 3}A_{10, 4} A_{9, 5}^2A_{8, 6}A_{7, 7}A_{8, 8}A_{9, 4}
 A_{8, 5}A_{7, 6}A_{9, 3}A_{8, 4}A_{7, 5}A_{6, 8}A_{9, 2}A_{8, 3}&\\
&\left. A_{7, 4}A_{6, 5}A_{5, 6}A_{4, 7}A_{9, 1}A_{8, 2}A_{7, 3}A_{6, 4}
 A_{5, 5}A_{4, 6}A_{4, 8}A_{4, 5}A_{4, 4}A_{4, 3}A_{4, 2}A_{4, 1}
\right)^{-1},&{\rm if }\ i=4,\\
\frac{\cc(15,5)}{\cc(15,6)\cc(15,8)} =&\cc(5,1)
\left(
A_{14, 6}A_{13, 7}A_{14, 8}A_{14, 5}A_{13, 6}A_{14, 4}A_{13, 5} \right.&\\
& A_{12, 8}A_{14, 3}A_{13, 4}A_{12, 5}A_{11, 6}A_{10, 7}A_{14, 2}
 A_{13, 3}A_{12, 4}A_{11, 5}A_{10, 6}A_{10, 8}A_{14, 1}A_{13, 2}&\\
& A_{12, 3}A_{11, 4} A_{10, 5}^2A_{9, 6}A_{8, 7}A_{9, 8}A_{10, 4}
 A_{9, 5}A_{8, 6}A_{10, 3}A_{9, 4}A_{8, 5}A_{7, 8}A_{10, 2}A_{9, 3}&\\
& \left.A_{8, 4}A_{7, 5}A_{6, 6}A_{5, 7}A_{10, 1}A_{9, 2}A_{8, 3}A_{7, 4}
 A_{6, 5}A_{5, 6}A_{5, 8}A_{5, 5}A_{5, 4}A_{5, 3}A_{5, 2}A_{5, 1}
\right)^{-1}&{\rm if}\  i=5,\\
\frac{\cc(15,6)}{\cc(15,7)} =&\cc(14,7){A_{14,7}}^{-1}
&{\rm if}\  i=6,\\
\cc(15,7) && {\rm if}\ i=7,\\
\cc(15,8) && {\rm if}\ i=8.
\end{array}
\right.
\label{e8-c}
\end{equation}


\underline{Type ${\rm F}_4$-case:}
The monomial of 
$\Delta_{w_0\Lambda_i,s_i\Lambda_i}(\Theta^-_{\bfii0}(\textbf{c}))$
in (\ref{lowest-term}) is equal to
We obtain
\begin{equation}
\left\{
\begin{array}{lll}
&\frac{\cc(6,1)}{\cc(6,2)} =\cc(1,1)
\left(
A_{5, 2} A_{4, 3}^2 A_{3, 4}^2 A_{4, 2} A_{3, 3}^2 A_{4, 1} A_{3, 2}^2 
A_{2, 3}^2 A_{1,4}^2 A_{3,1}A_{2, 2}A_{1,3}^2 A_{1, 2} A_{1,1}
\right)^{-1}& {\rm if}\ i=1,\\
&\frac{\cc(6,2)}{{\cc(6,3)}^2} =\cc(2,1)
\left(
A_{5, 3}^2 A_{4, 4}^2 A_{5, 2} A_{4, 3}^2 A_{5, 1} A_{4, 2}^2 A_{3, 3}^2 
A_{2,4}^2 A_{4,1} A_{3, 2} A_{2, 3}^2 A_{2, 2} A_{2,1}
\right)^{-1}& {\rm if}\ i=2,\\
&\frac{\cc(6,3)}{\cc(6,4)} =\cc(5,4)A_{5,4}^{-1}
& {\rm if}\ i=3,\\
&\cc(6,4)& {\rm if}\ i=4.
\end{array}
\label{f4-c}
\right.
\end{equation}
Here note that in the last term of each equation
we rewrite the monomial using $\TY(c,i,j)$ and $A_{m,k}^{-1}$ 
for the later use.
\end{lem}
\begin{pro}\label{pro-ef}
Let $Y$ be the Laurent monomial (\ref{lowest-term}) in
$\Delta_{w_0\Lambda_i,s_i\Lambda_i}(\Theta^-_{\bfii0}(\textbf{c}))$.
Then, we find that each Laurent monomial in 
$\Delta_{w_0\Lambda_i,s_i\Lambda_i}(\Theta^-_{\bfii0}(\textbf{c}))$
is in the form
\[
Y
\times \prod_{(s,j)\in T_\bfii0}A^{l_{s,j}}_{s,j}
\]
with some $l_{s,j}\in\mathbb{Z}$.
\end{pro}
{\sl  Proof.}
Let $\pi=(\gamma_0,\gamma_{s,j})_{(s,j)\in T_\bfii0}$ be a pre-$\bfii0$-trail, 
where  $\gamma_0=-w_0\Lm_i$.
If a pre-$\bfi$-trail $\pi'=(\gamma_0,\gamma'_{s,j})_{(s,j)\in T_\bfi}$ satisfies
\begin{equation}\label{pre-trans}
\begin{split}
& \gamma'_{s,j}=\gamma_{s,j}-\alpha_{k}, \q
{\rm if}\ (t,k)\leq (s,j)< (t+1,k),\\
& \gamma'_{s,j}=\gamma_{s,j},\qq\q\,\, {\rm otherwise},
\end{split}
\end{equation}
then we obtain
\[
d_{t,k}(\pi')=\frac{\gamma'_{(t,k)-1}+\gamma'_{t,k}}{2}(h_{k})=
\frac{\gamma_{(t,k)-1}+(\gamma_{t,k}-\alpha_k)}{2}(h_{k})
=d_{t,k}(\pi)-1.
\]
where the notation $(t,k)-1$ is the entry in $T_\bfi$ just before $(t,k)$.
We also get if $l>k$ then
\[
d_{t,l}(\pi')=\frac{\gamma'_{(t,l)-1}+\gamma'_{t,l}}{2}(h_{l})=
\frac{(\gamma_{(t,l)-1}-\alpha_k)+(\gamma_{t,l}-\alpha_k)}{2}(h_{l})
=d_{t,l}(\pi)-{\bf a}_{l,k},
\]
if $l<k$ then
\[
d_{t+1,l}(\pi')=\frac{\gamma'_{(t+1,l)-1}+\gamma'_{t+1,l}}{2}(h_{l})=
\frac{(\gamma_{(t+1,l)-1}-\alpha_k)+(\gamma_{t+1,l}-\alpha_k)}{2}(h_{l})
=d_{t+1,l}(\pi)-{\bf a}_{l,k},
\]
and
\[
d_{t+1,k}(\pi')=\frac{\gamma'_{(t+1,k)-1}+\gamma'_{t+1,k}}{2}(h_{k})=
\frac{(\gamma_{(t+1,k)-1}-\alpha_k)+\gamma_{t+1,k}}{2}(h_{k})
=d_{t+1,k}(\pi)-1.
\]
Otherwise,  $d_{s,j}(\pi')=d_{s,j}(\pi)$. Thus, we have
\begin{eqnarray*}
\prod_{s,j} {\cc(s,j)}^{d_{s,j}(\pi')}
&=&
(\prod_{s,j} {\cc(s,j)}^{d_{s,j}(\pi)})
\times
{\cc(t,k)}^{-1}{\cc(t+1,k)}^{-1}
\prod_{(t,k)< (m,l)<(t+1,k)} {\cc(m,l)}^{-{\bf a}_{l,k}}\\
&=&
(\prod_{s,j} {\cc(s,j)}^{d_{s,j}(\pi)}) \times A_{t,k}^{-1}.
\end{eqnarray*}

Following Definition \ref{pretrail}, we take non-negative integers $\{c_{s,j}\}_{(s,j)\in T_{\textbf{i}_0}}$,
$\{c'_{s,j}\}_{(s,j)\in T_{\textbf{i}_0}}$ such that $\gamma_{(s,j)-1}-\gamma_{s,j}=c_{s,j}\alpha_j$,
$\gamma'_{(s,j)-1}-\gamma'_{s,j}=c'_{s,j}\alpha_j$. Then it holds
\[
c_{s,j}=c'_{s,j}\ \text{if }(s,j)\neq (t,k),\ (t+1,k),\qquad c'_{t,k}=c_{t,k}+1,\quad c'_{t+1,k}=c_{t+1,k}-1.
\]
Each pre-$\textbf{i}_0$-trail $(\gamma''_{0},\gamma''_{s,j})_{(s,j)\in T_{\textbf{i}_0}}$
from $-w_0\Lambda_i$ to $-s_i\Lambda_i$ is uniquely determined from 
non-negative integers $\{c''_{s,j}\}_{(s,j)\in T_{\textbf{i}_0}}$ such that
$\gamma''_{(s,j)-1}-\gamma''_{s,j}=c''_{s,j}\alpha_j$.
Note that for $m\in I$,
 the value $\sum_{s\in\mathbb{Z}_{\geq1};(s,m)\in T_{\textbf{i}_0}}c''_{s,m}$ is
constant not according to the choice of pre-$\textbf{i}_0$-trail.
Therefore, for two pre-$\textbf{i}_0$-trails, 
one is obtained from the other via iteration of the transformations
$(\gamma_{s,j})\leftrightarrow (\gamma'_{s,j})$ in (\ref{pre-trans}).
In particular, each $\textbf{i}_0$-trail is obtained from
the $\textbf{i}_0$-trail $\pi_0$ in the proof of Theorem \ref{lowest-term-thm},
which means our claim.
\qed
%

Here summarizing the results in this section, 
one has the following theorem, 
which is similar to 
Theorem \ref{mono-g2} and 
will be needed in Sect.8 to show the connectedness of the 
crystal $B_\bfi$:
\begin{thm}\label{mono-ce}
For type ${\rm E}_6,\ {\rm E}_7,\ {\rm E}_8$ and ${\rm F}_4$,
 and the specific reduced longest word
$\bfii0$ as in \eqref{EF-i}, one gets that the generalized minor
$\Del_{w_0\Lm_k,s_k\Lm_k}(\Theta^-_\bfii0(c_1,\cd,c_N))$ ($k\in I$)
is a positive sum of monomials
in the form:
\begin{equation}
c_j\prod_{(m,i)\in T_\bfii0}A_{m,i}^{l_{m,i}},
\end{equation}
for some $j\in\{1,2,\cd,N=L(w_0)\}$ and some integers $l_{m,i}$. 
Note that the monomial $A_{m,i}$ is given as \eqref{asi-ef}.
\end{thm}
{\sl Proof.}
It is immediate from the explicit forms of monomials 
in Lemma \ref{lem-ef} and Proposition \ref{pro-ef}.\qed

\renewcommand{\thesection}{\arabic{section}}
\section{Connectedness of the cellular crystal $B_\bfi$}
\label{conn-cell}
\setcounter{equation}{0}
\renewcommand{\theequation}{\thesection.\arabic{equation}}

In this section, we shall show the connectedness of the cellular crystal 
$B_\bfi=B_{i_1}\ot\cd\ot  B_{i_l}$ associated with a reduced word 
$\bfi=i_1i_2\cd i_l$ for a 
Weyl group element $w\in W$.

\subsection{Subset  $\HI\subset B_{\bfi}$}

Let $\bfi=i_1i_2\cd i_N$ be the reduced longest word 
of the Weyl group $W$ and the linear function $\beta_k(x)$ is given as 
in \eqref{betak}:
\[
\beta_k(x)=x_k+\sum_{k<j<\kp}{\bf a}_{i_k,i_j}x_j+x_{\kp}.
\]
Now, define the set $\HI\subset \bbZ^N(=B_{\bfi})$ as follows:
\begin{equation}
\HI:=\{x\in\bbZ^N(=B_{\bfi})\,|\,
\beta_k(x)=0 \,\text{for any }k\in [1,N]\,\,\text{such that }
\kp\leq N\},
\label{cH}
\end{equation}
where by the correspondence $(x_1,\cd,x_N)\leftrightarrow
(-x_1)_{i_1}\ot\cd\ot(-x_N)_{i_N}$ we identify $\bbZ^N$ with $B_\bfi$.

This $\HI$ is presented in the following simple form.
First, for $\bfi=i_1i_2\cd i_N$, define
$\TY(\al,k,)=s_{i_N}s_{i_{N-1}}\cd s_{i_{k+1}}(\al_{i_k})$ for 
the simple reflections $s_1,\cd,s_n$ and simple roots $\al_1,\cd,\al_n$ 
of the Langlands dual Lie algebra 
$\LL(\ge)$.
It is well-known that $\LL(\Del)=\{\TY(\al,k,)\,|\,1\leq k\leq N\}$
just coincides with the set of all positive roots of $\LL(\ge)$.
Thus, we can write 
\begin{equation}
\TY(\al,k,)=\sum_{i\in I}\TY(m,k,i)\al_i
\label{alk}
\end{equation}
with certain non-negative integral coefficients $\TY(m,k,i)$. 
Then,  we obtain 
\begin{lem}\label{HI-lem}
\begin{equation}
\HI=\left\{\left.\left(\sum_i\TY(m,1,i)h_i,\sum_i\TY(m,2,i)h_i,\cd,
\sum_i\TY(m,N,i)h_i\right)\,\right|\,(h_1,h_2,\cd,h_n)\in\bbZ^n\right\}
\label{Hm}
\end{equation}
\end{lem}
{\sl Proof.}
For $k$ such that $\kp\leq N$, we have
\begin{eqnarray*}
\TY(\al,k,)&=&s_{i_N}s_{i_{N-1}}\cd s_{i_{k+1}}(\al_{i_k})=
s_{i_N}s_{i_{N-1}}\cd s_{i_{k+2}}(\al_{i_k}-{\bf a}'_{i_{k+1},i_k}\al_{i_{k+1}})\\
&=&
s_{i_N}s_{i_{N-1}}\cd s_{i_{k+2}}(\al_{i_k})
-{\bf a}'_{i_{k+1},i_k}\TY(\al,k+1,)=\cd\\
&=&
s_{i_N}s_{i_{N-1}}\cd s_{i_{\kp}}(\al_{i_{\kp}})-
\sum_{k+1\leq j\leq \kp-1}{\bf a}'_{i_j,i_k}\TY(\al,j,)\\
&=&
-s_{i_N}s_{i_{N-1}}\cd s_{i_{\kp+1}}(\al_{i_{\kp}})-
\sum_{k< j< \kp}{\bf a}'_{i_j,i_k}\TY(\al,j,)\\
&=&-\TY(\al,\kp,)-\sum_{k< j< \kp}{\bf a}'_{i_j,i_k}\TY(\al,j,)\\
&=&-\TY(\al,\kp,)-\sum_{k< j< \kp}{\bf a}_{i_k,i_j}\TY(\al,j,),
\end{eqnarray*}
where ${\bf a}'_{ij}$ is an entry of the Cartan matrix of $\LL(\ge)$ and then 
${\bf a}'_{ij}={\bf a}_{ji}$.
Thus, substituting \eqref{alk} to the above, we obtain 
\begin{equation}
\sum_{i\in I}\TY(m,k,i)\al_i=-\sum_{i\in I}\TY(m,\kp,i)\al_i
-\sum_{k< j <\kp}{\bf a}_{i_k,i_j}(\sum_{i\in I}\TY(m,j,i)\al_i).
\end{equation}
Furthermore, since $\{\al_1,\cd,\al_n\}$ is linearly independent, 
for any $i\in  I$, we have
\begin{equation}
\TY(m,k,i)=-\TY(m,\kp,i)-\sum_{k<j<\kp}{\bf a}_{i_k,i_j}\TY(m,j,i).
\label{mki}
\end{equation}
For $k=1,2,\cd,N$ and $(h_1,\cd,h_n)\in \bbZ^n$, set
$H_k:=\sum_{i\in I}\TY(m,k,i)h_i$. By \eqref{mki} we get
\begin{equation}
H_k=\sum_{i\in I}(-\TY(m,\kp,i)-\sum_{k<j<\kp}{\bf a}_{i_k,i_j}\TY(m,j,i))h_i
=-H_{\kp}-\sum_{k<j<\kp}{\bf a}_{i_k,i_j}H_j,
\end{equation}
which implies that $H_k,H_{k+1},\cd,H_{\kp}$ satisfies 
the equation $\beta_k(H_1,H_2,\cd,H_N)=0$. Thus, we find that 
$(H_1,H_2,\cd,H_N)$ is a solution of the equations
$\beta_k(x)=0$ for $k$ such that $\kp\leq N$. Therefore, we obtain
L.H.S. of \eqref{Hm}$\supset$ R.H.S. of \eqref{Hm}.
On the other hand, the solution space $\cH(\bbQ)$ 
in $\bbQ^N$ of the linear equalities $\beta_k(x)=0$ for 
$k$ such that $\kp\leq N$ has the dimension $n$ since 
the number of the equations is $N-n$ and $\{\beta_k\}_k$ are all independent.
Thus, we find 
\[
\HI\subset\cH(\bbQ)\cap \bbZ^N,
\]
and 
\[
\cH(\bbQ)=\left\{\left.\left(\sum_i\TY(m,1,i)h_i,\sum_i\TY(m,2,i)h_i,\cd,
\sum_i\TY(m,N,i)h_i\right)\,\right|\,(h_1,h_2,\cd,h_n)\in\bbQ^n\right\}.
\]
Since $\{\TY(\al,1,),\cd,\TY(\al,N,)\}$ includes all simple roots
$\{\al_1,\cd,\al_n\}$, there exist $k_1,\cd,k_n\subset \{1,\cd,N\}$ 
such that 
$\sum_{i\in I}\TY(m,k_l,i)h_i=h_l$ ($1\leq l\leq n)$, which means
R.H.S. of \eqref{Hm}$=\cH(\bbQ)\cap \bbZ^N$ and then we obtain 
$\HI\subset $R.H.S.of \eqref{Hm}.\qed

We have the following result: 
\begin{pro}\label{eitfit-pro}
For any $i\in I$, $x=(x_1,\cd,x_N)\in B_{\bfi}$ and 
$H=(H_1,\cd,H_N)\in\HI$, we get
\begin{equation}
\eit(x+H)=\eit(x)+H,\qq
\fit(x+H)=\fit(x)+H.
\label{eihfih}
\end{equation}
\end{pro}
{\sl Proof.}
Now let us recall Theorem \ref{thm-KN}. For 
$x=(x_1,\cd,x_N)\in B_\bfi$ the function $a_k=a_k(x)$ in the
theorem can be rewritten as 
\[
a_k(x)=\begin{cases}
-x_k-\sum_{1\leq \nu<k}{\bf a}_{i_k,i_\nu}x_\nu&\text{ if }i_k=i,\\
-\ify&\text{ if }i_k\ne i,
\end{cases}
\]
where $({\bf a}_{ij})$ is the Cartan matrix for $\ge$.
If $i_k=i$ and $\kp\leq N$, then we find that 
\[
a_k(x)-a_{\kp}(x)=x_k+x_{\kp}+\sum_{k<\nu<\kp}{\bf a}_{i_k,i_\nu}x_j=\beta_k(x).
\]
Thus, we obtain 
\begin{equation}
a_k(x+H)-a_{\kp}(x+H)=\beta_k(x+H)=\beta_k(x)+\beta_k(H)=\beta_k(x)
=a_k(x)-a_{\kp}(x).
\label{xHx}
\end{equation}
Indeed, $k_e$ and $k_f$ in Theorem \ref{thm-KN} are 
determined by the differences $a_k-a_\kp=\beta_k$'s and by 
\eqref{xHx} we obtain the desired results.\qed

\subsection{Condition $\CHI$}

For any reduced longest word $\bfi=i_1i_2\cd i_N$, the condition 
$\CHI$ is defined as follows:
\begin{df}[The Condition $\CHI$]\label{cond-h}
We say the condition $\CHI$ holds if 
there exists the finite set of linear functions $\Xi_{\bfi}$:
\[
\begin{array}{l}\displaystyle
\Xi_{\bfi}\subset
\left\{
\vp(x)\in (\bbQ^N)^*\,\left|\,\vp(x)=x_j-\sum_{1\leq k\leq N, \TY(k,+,)
\leq N}c_k\beta_k(x)\,\,(c_k\in\bbZ,\,\,1\leq j\leq N)\right.
\right\}\\
\text{such that  }\,\, B(\ify)=\{x\in B_{\bfi}(=\bbZ^N)\,|\,
\vp(x)\geq0\,\,(\vp\in\Xi_\bfi)\},
\end{array}
\]
where $\beta_k(x)$ is the function defined as in \eqref{betak}.
\end{df}

Suppose that the condition $\CHI$ holds. 
The crystal $B(\ify)$ is given by
\begin{equation}
B(\ify)=\{x\in B_\bfi\,|\,\varphi(x)\geq0(\forall \vp\in\Xi_\bfi)\}.
\label{bify}
\end{equation}
For $H\in\cH_\bfi$ set
\begin{equation}
B^H(\ify):=\{x+H\in B_\bfi\,|\,x\in B(\ify)\}.
\label{BH}
\end{equation}
Here if we define $\eit(x)=0$ (resp. $\fit(x)=0$) for 
$\eit(x)\not\in B^H(\ify)$ (resp. $\fit(x)\not\in B^H(\ify)$), then 
$B^H(\ify)$ becomes a crystal. 
The following lemma is immediate from Proposition \ref{eitfit-pro} and the
connectedness of $B(\ify)$.
\begin{lem}
The crystal 
$B^H(\ify)$ is connected.
\end{lem}

As in Sect.\ref{fb-ABCDG}, we have seen that 
the explicit forms of $\Phi^{(+)}=\Phi^{(+)}_{\bfi_0}$ 
for the types $\AAN,\BBN,\CCN,\DDN$ and $\G22$, where 
$\bfi_0$ is the specific longest reduced word as in \eqref{i0}.
Here, let us define the set of linear functions $\Xi_{\bfi}$ as follows:
For the tropicalization $\wtil\Phi^{(+)}_{\bfi}$, if we have 
\begin{equation}
\wtil\Phi^{(+)}_{\bfi}(x)=\min(f_1(x),f_2(x),\cd,f_M(x)),
\label{min-f}
\end{equation}
then $\Xi_{\bfi}:=\{f_1(x),f_2(x),\cd,f_M(x)\}$. Namely,
if $\Phi^{(+)}_{\bfi}$ is a Laurent polynomial and 
each Laurent monomial is $F_j$, then $f_j(x)={\mathcal TR}(F_j)(x)$.

In the previous section, we have also  seen that 
certain explicit forms of $\Phi^{(+)}=\Phi^{(+)}_{\bfi_0}$ 
for the types $\EEE6$, $\EEE7$, $\EEE8$ and $\F44$.
Thus, we get
\begin{pro}\label{H0-ABCDG}
For the all types of $\AAN,\BBN,\CCN,\DDN$, $\EEE6$, $\EEE7$, $\EEE8$, $\F44$  and $\G22$, 
the set of linear functions $\Xi_{\bfi_0}$ obtained form 
$\wtil\Phi^{(+)}_{\bfi_0}$ satisfies the condition $\CHZ$.
\end{pro}
{\sl Proof.}
First let us show for $\ge=\AAN,\CCN,\BBN,\DDN$  and $\G22$. 
For those types, by Theorem \ref{del-mono}, Theorem \ref{thm-mono-c1}
Theorem \ref{thm-mono-b1}, Theorem \ref{thm-mono-d1} and 
Theorem \ref{mono-g2}, each monomial in $\Phi^{(+)}_{\bfi_0}(c)$ 
is in the form 
\[
c^{(j)}_i\prod_{p=1}^rA_{i_p,m_p}^{-1},
\]
which means that each $f_j(x)$ in \eqref{min-f} is in the form 
\[
x^{(j)}_i-\sum_{p=1}^r \beta_{i_p,m_p}^\vee,
\]
since one has that ${\mathcal TR}(A_{i,m})(x)=\beta_{i,m}^\vee(x)
:=\TY(x,m,i)+\TY(x,m+1,i)+\sum_{(i,m)<(j,r)<(i,m+1)}
{\bf a}_{j,i}\TY(x,r,j)$. 
This implies that 
the set of linear functions $\Xi_{\bfi_0}$ holds the 
condition $\CHZ$ for Langlands dual $\LL(\ge)$. 

Let us show for $\EEE6$, $\EEE7$, $\EEE8$ and $\F44$.
By the similar argument  to the types $\AAN,\CCN,\BBN,\DDN$  and $\G22$, 
and the results in the previous section,
we know that each monomial in $\Phi^{(+)}_{\bfi_0}(c)$ is in the form
\[
c^{(j)}_i\prod_{p=1}^rA_{i_p,m_p}^{l_{i_p,m_p}}\,\qq
(l_{i_p,m_p}\in \bbZ),
\]
and then its tropicalization is also in the form 
\[
x^{(j)}_i+\sum_{p=1}^r l_{i_p,m_p}\beta_{i_p,m_p}^\vee.
\]
Thus, 
we get that 
the set of linear functions $\Xi_{\bfi_0}$ holds the 
condition $\CHZ$ for Langlands dual $\LL(\ge)$. 
\qed

\subsection{Connectedness of the cellular crystals}

First let us see the following lemma:
\begin{lem}\label{Hi-cnn}
For a reduced longest word $\bfi=i_1i_2\cd i_N$, if 
the condition $\cH_\bfi$ holds, then the cellular crystal 
$B_\bfi=B_{i_1}\ot\cd\ot B_{i_N}$ is connected as a crystal graph.
\end{lem}
{\sl Proof.}
Let us show 
\begin{equation}
B(\ify)\cap B^H(\ify)\ne\emptyset \q\text{for any }H\in\HI.
\label{nonempty}
\end{equation}
Since the set of linear functions $\Xi_\bfi$ is finite and all elements 
in $\Xi_\bfi$ have no constant term, 
there exists $\til x=(x_1,\cd,x_N)\in B(\ify)$ such that 
$\vp(\til x)\geq1$ for any $\vp\in \Xi_\bfi$.
By the condition $\CHI$, each $\vp(x)$ can be written in the form:
\begin{equation}
\vp(x)=x_k-\sum_j m^{(\vp)}_j\beta_j(x).\label{vp-exp}
\end{equation}
Here, for $H=(H_1,H_2,\cd, H_N)\in\HI$ we set 
\[
\wtil H:=\max(|H_1|,|H_2|,\cd,|H_N|).
\]
By \eqref{vp-exp}  we get 
\[
\vp(\wtil H\cdot \til x)=\wtil H\cdot \vp(\til x)=
\wtil H\left(x_k-\sum_j m^{(\vp)}_j\beta_j(\til x)\right)
\geq0,
\]
which implies that $\wtil H\cdot \til x\in B(\ify)$ and then 
$\wtil H\cdot \til x+H\in B^H(\ify)$. 
Regarding the facts $\vp(\wtil x)\geq1$  and $\beta_j(H)=0$, 
we have
\begin{eqnarray*}
\vp(\wtil H\cdot \til x+H)&=&
\wtil H\cdot \til x_k+H_k
-\sum_jm^{(\vp)}_j\beta_j(\wtil H\cdot\til x+H)
=\wtil H\cdot\vp(\wtil x)+H_k
\geq\wtil H+H_k\geq0,
\end{eqnarray*}
which means that $\wtil H\cdot\til x+H\in B(\ify)$ and then we obtain 
\eqref{nonempty} by 
\[
\wtil H\cdot\til x+H\in B(\ify)\cap B^H(\ify).
\]
It follows from \eqref{nonempty} and connectedness of 
$B(\ify)$ that $\cup_{H\in\HI}B^H(\ify)$
is connected. Here the actions of $\eit$ and $\fit$ on $\cup_{H\in\HI}B^H(\ify)$
are considered as follows:
for some $x\in B^H(\ify)$, suppose $\eit(x)\not\in B^H(\ify)$
(resp. $\fit(x)\not\in B^H(\ify)$). If $\eit(x)\in B^{H'}(\ify)$
(resp. $\fit(x)\in B^{H'}(\ify)$) for some $H'\in \HI$, we understand 
$\eit(x)\ne0$ (resp. $\fit(x)\ne0$).
Indeed, 
for any $x, y\in \cup_{H\in\HI}B^H(\ify)$ 
there exist $H,H'\in\HI$ such that $x\in B^H(\ify)$ and $y\in B^{H'}(\ify)$.
Then, by \eqref{nonempty} we find that there exist $x'\in B(\ify)\cap B^H(\ify)$
and $y'\in B(\ify)\cap B^{H'}(\ify)$. Thus, we know that $(x,x')$, $(x',y')$ and 
$(y',y)$ are connected each other and then $(x,y)$ is connected.

Therefore, for the purpose 
it suffices to show that
\begin{equation}
B_\bfi=\bigcup_{H\in\HI}B^H(\ify).
\label{Bi=union}
\end{equation}
For $x\in B_\bfi$, $\vp\in \Xi_{\bfi}$ and 
$H=(H_1,\cd,H_N)\in \HI$, by a similar way as above we have 
\[
\vp(x+H)=\vp(x)+H_k,
\]
for some $k\in [1,N]$, 
Thus, 
by Lemma \ref{HI-lem}, $H_k$ can be written in 
\[
H_k=\sum_i m^{(k)}_ih_i,
\]
where $(h_1,\cd,h_n)\in\bbZ^n$
and $(m^{(k)}_1,\cd,m^{(k)}_n)\in\bbZ_{\geq0}^n\setminus\{(0,0,\cd,0)\}$. Thus, 
taking all $h_i$'s  to be sufficiently large such that
$\vp(x)+H_k\geq0$ for any $k$ and $\vp\in \Xi_\bfi$, then we 
know that $x+H\in B(\ify)$ and then $x\in \cup_{H\in\HI}B^H(\ify)$. 
Therefore, \eqref{Bi=union} is obtained and then 
we completed the proof of the lemma.\qed

The following is the main theorem of this section.
\begin{thm}\label{main}
For an arbitrary finite-dimensional simple Lie algebra $\ge$ 
and any associated 
reduced word $\io=i_1i_2\cd i_k$ of any Weyl group element 
$w\in W$, the cellular crystal 
$B_\io=B_{i_1}\ot B_{i_2}\ot\cd\ot B_{i_k}$ is connected as 
a crystal graph.
\end{thm}
\begin{rem}\label{rem-ten}
\begin{enumerate}
\item If a tensor product of crystals $B\ot B'$ is connected, 
then both the crystals $B$ and $B'$ are connected. Indeed, 
suppose that  $B$ is not connected and $b_1,b_2\in B$ are in the 
different connected components. Then, $b_1\ot b'$ and $b_2\ot b'$ 
are never connected for any $b'\in B'$. By the actions of $\eit$ and $\fit$ 
on the tensor product we have
\[
\til X_{j_1}\til X_{j_2}\cd \til X_{j_L}(b_1\ot b')=
\til X_{l_1}\til X_{l_2}\cd \til X_{l_m}(b_1)\ot  
\til X_{l_{m_1}}\til X_{l_{m+2}}\cd 
\til X_{l_L}(b'),
\]
where $X=e$ or $f$ and $\{j_1,\cd,j_L\}=\{l_1.\cd,l_L\}$. 
We know that $\til X_{l_1}\til X_{l_2}\cd \til X_{l_m}(b_1)$ never be $b_2$.
\item
For any reduced word $i_1i_2\cd i_k$, there exists
$i_{k+1},i_{k+2},\cd ,i_N\in I$
such that $i_1i_2\cd i_ki_{k+1}\cd i_N$ becomes a reduced longest word.
\item
For the proof of the theorem, by the above (i) (ii)
it suffices to show that $B_{i_1}\ot \cd\ot B_{i_N}$ is connected for 
any reduced longest word $i_1\cd i_N$. Furthermore, 
by the braid-type isomorphisms, we may show the connectedness for 
some specific reduced longest word.
\item
If a word ${\bf  j}=j_1\cd j_k$ is not reduced, then the crystal
$B_{\bf j}=B_{j_1}\ot \cd\ot B_{j_k}$ is not connected. Indeed, 
by applying the braid-type isomorphisms to the crystal  $B_{\bf  j}$, one has 
a crystal $B_{\bf  m}=B_{m_1}\ot \cd\ot B_{m_k}$ such that 
$m_p=m_{p+1}$ for some $p\in\{1,2,\cd,k-1\}$. 
And then it is known by \cite[Example 1.3.4]{K3} 
that the crystal 
$B_m\ot B_m=\oplus_{l\in\bbZ}B_m\ot T_{l\al_m}$ $(m\in I)$ 
is not  connected. 

\end{enumerate}
\end{rem}
{\sl Proof of Theorem \ref{main}.}
By the remark above, to prove the theorem, it 
suffices to show the connectedness of the cellular 
crystal $B_{\bfi_0}=B_{i_1}\ot\cd\ot B_{i_N}$ for some 
specific reduced longest word $\bfi_0$. Therefore,  to complete the proof, 
we may show  
that the condition ${\bf H}_{\bfi_0}$ holds, which has been already 
shown in Proposition \ref{H0-ABCDG}.\qed

\begin{df}
A geometric crystal $\bbX=(X,\{e_i\},\{\gamma_i\},\{\vep_i\})$ is 
{\it prehomogeneous} if there exists a Zariski open dense subset 
$\Omega\subset X$ which is an orbit by the actions of $e_i$'s.
\end{df}
\begin{cor}\label{prehom}
The geometric crystal $B^-_{\bfi}$ for a reduced word $\bfi=i_1i_2\cd i_k$ 
is prehomogeneous.
\end{cor}
{\sl Proof.}
In \cite[Theorem 3.3]{KNO2}, we know that
for a positive geometric crystal $X$, if 
its tropicalization $B={\mathcal TR}(X)$
is connected in the sense of crystal graph, 
then $X$ is prehomogeneous. 
Thus, by Theorem \ref{main} we find that 
$B_\bfi={\mathcal TR}(B^-_\bfi)$ is connected and then 
we obtain the desired result.\qed
\bibliographystyle{amsalpha}

\end{document}